\documentclass[a4paper,11pt]{article}
\usepackage{amsmath, amssymb}
\usepackage{amsthm}% for \proof,\newtheorem
\usepackage{latexsym}

\usepackage{subfiles}
\usepackage{comment}
\usepackage{setspace} % for gyoukan chousei
\usepackage{authblk} %\author関連

\usepackage{enumerate} %\begin{enumerate}[{(}i{)}]
\usepackage{cases}
  
\usepackage{ascmac} % for \screen,\itembox,\shadebox,highlight,\yen
\usepackage{framed}
\usepackage{fancybox}

\usepackage{cancel}
\usepackage{url}

\usepackage[left=19.91truemm,right=19.91truemm,top=26.4truemm,bottom=42truemm]{geometry}

\usepackage[colorlinks=true]{hyperref}
\hypersetup{
    colorlinks=true,
    citecolor=blue,
    linkcolor=blue,
}
\usepackage{cleveref}

\usepackage[backend=bibtex,style=alphabetic,isbn=false,url=false,eprint=false,doi=false]{biblatex}
\usepackage[normalem]{ulem}

\usepackage{xcolor}
\usepackage{layout} % for layout 
\usepackage{tikz} 
\usepackage{graphicx}
\usepackage[labelformat=parens, font=footnotesize]{caption}

\usepackage{subcaption}

\usepackage{here}

\usepackage{dsfont}
\usepackage{bbm} %\mathbbm{ABCabc\nabla}
\usepackage[bbgreekl]{mathbbol} %\mathbb{abcABC\nabla\bbalpha\bbbeta\Delta  }
\usepackage{upgreek}
\usepackage{mathrsfs}

\newtheorem{theorem}{Theorem}[section]
\newtheorem{ass}[theorem]{Assumption}

\newtheorem{lemma}[theorem]{Lemma}

\newtheorem{proposition}[theorem]{Proposition}

\newtheorem{remark}[theorem]{Remark}

\numberwithin{equation}{section}
\newtheorem{theorem*}{Theorem}
\newtheorem{ass*}[theorem*]{Assumption}
\newtheorem{note*}[theorem*]{Note}
\newtheorem{lemma*}[theorem*]{Lemma}
\newtheorem{definition*}[theorem*]{Definition}
\newtheorem{proposition*}[theorem*]{Proposition}
\newtheorem{corollary*}[theorem*]{Corollary}
\newtheorem{remark*}[theorem*]{Remark}
\newtheorem{example*}[theorem*]{Example}
\numberwithin{equation}{section}

\usepackage{comment}

\makeatletter
\newsavebox{\@brx}
\newcommand{\llangle}[1][]{\savebox{\@brx}{\(\m@th{#1\langle}\)}%
  \mathopen{\copy\@brx\kern-0.5\wd\@brx\usebox{\@brx}}}
\newcommand{\rrangle}[1][]{\savebox{\@brx}{\(\m@th{#1\rangle}\)}%
  \mathclose{\copy\@brx\kern-0.5\wd\@brx\usebox{\@brx}}}
\makeatother

\newcommand\vspsm{\vspace{5mm}}

\renewcommand{\emptyset}{\varnothing}

\newcommand\sfz{{\sf z}}

\newcommand\sfi{{\tt i}}
\newcommand\tti{{\tt i}}

\newcommand{\stout}[1]{\ifmmode\text{\sout{\ensuremath{#1}}}\else\sout{#1}\fi}

\newcommand{\myred}{\color[rgb]{0.8,0,0}}
\newcommand{\delc}[1]{}
\newcommand{\delb}[1]{[]}

\newcommand\vpr{{v'}}

\newcommand{\bnorm}[1]{{\bigl\| #1 \bigr\|}}

\newcommand{\snorm}[1]{{\left\| #1 \right\|}}
\newcommand{\norm}[1]{{\left\| #1 \right\|}}

\newcommand{\abs}[1]{{\left| #1 \right|}}

\newcommand{\babs}[1]{{\bigl| #1 \bigr|}}

\newcommand{\cbr}[1]{{\left\{ #1 \right\} }}
\newcommand{\ncbr}[1]{{\{ #1 \} }}
\newcommand{\bcbr}[1]{{\bigl\{ #1 \bigr\} }}
\newcommand{\Bcbr}[1]{{\Bigl\{ #1 \Bigr\} }}

\newcommand{\rbr}[1]{{\left( #1 \right) }}
\newcommand{\nrbr}[1]{{( #1 ) }}
\newcommand{\brbr}[1]{{\bigl( #1 \bigr) }}
\newcommand{\Brbr}[1]{{\Bigl( #1 \Bigr) }}

\newcommand{\sbr}[1]{{\left[ #1 \right] }}

\newcommand{\nabr}[1]{{\langle #1 \rangle }}
\newcommand{\abr}[1]{{\left\langle #1 \right\rangle }}
\newcommand{\babr}[1]{{\bigl\langle #1 \bigr\rangle }}
\newcommand{\Babr}[1]{{\Bigl\langle #1 \Bigr\rangle }}

\newcommand{\subotimes}[1]{{\underset{ #1 }{\otimes}}}

\newcommand\qtor{{\tt qTor}}
\newcommand\mqtor{{\tt mqTor}}
\newcommand\qtan{{\tt qTan}}

\newcommand{\mR}{\mathfrak R}
\newcommand{\mS}{\mathfrak S}
\newcommand{\mT}{\mathfrak T}

\newcommand{\scrC}{{\mathscr C}}

\newcommand{\bbf}{{\mathbb f}}

\newcommand{\bbone}{{\mathbb 1}}

\newcommand{\tfor}{\quad\text{for}\ }
\newcommand{\tforsm}{\ \text{for}\ }

\newcommand{\tifsm}{\ \text{if}\ }

\newcommand{\tand}{\quad\text{and}\quad}

\newcommand{\tandsm}{\ \text{and}\ }

\newcommand{\torsm}{\ \text{or}\ }
\def\mH{{\mathfrak H}}

\def\bd{\begin{description}}
\def\ed{\end{description}}

\def\D2{\bbD_{2,\infty-}}

\def\n{{\bf n}}

\def\C{{\bf C}}
\def\D{{\bf D}}

\def\cala{{\cal A}}

\def\cale{{\cal E}}
\def\calf{{\cal F}}

\def\calh{{\cal H}}
\def\cali{{\cal I}}

\def\yeq{\>=\>}

\def\simleq{\ \raisebox{-.7ex}{$\stackrel{{\textstyle <}}{\sim}$}\ }

\def\half{\frac{1}{2}}

\def\nn{\nonumber}
\def\be{\begin{equation}}
\def\ee{\end{equation}}
\def\bea{\begin{eqnarray}}
\def\eea{\end{eqnarray}}
\def\beas{\begin{eqnarray*}}
\def\eeas{\end{eqnarray*}}
\def\bi{\begin{itemize}}
\def\ei{\end{itemize}}

\def\bd{\begin{description}}
\def\ed{\end{description}}
\newcommand{\bbD}{{\mathbb D}}

\newcommand{\bbN}{{\mathbb N}}

\newcommand{\bbR}{{\mathbb R}}

\newcommand{\bbT}{{\mathbb T}}

\newcommand{\bbZ}{{\mathbb Z}}

\newcommand{\barq}{{{\bar  q}}}
\newcommand{\bartheta}{{ {\bar  \theta}}}
\newcommand{\Comp}{\scrC}

\usepackage{stmaryrd}
\newcommand{\secDiff}[2]{\Delta^{(2)}_{#1,#2}}

\newcommand{\ntwo}{{\bbN_2}}

\newcommand{\constInProd}{c_H}

\newcommand{\calii}{\cali^{(i)}}

\newcommand{\bbd}{{\mathbb d}}

\newcommand{\wt}[1]{\widehat{#1}}
\newcommand{\wtv}{\wt{V}}

\newcommand{\whv}{\widehat{V}}
\newcommand{\extv}{\widetilde{V}}

\newcommand{\tp}{\hat{p}}
\newcommand{\hp}{\hat{p}}
\newcommand{\tv}{\hat{v}}
\newcommand{\hv}{\hat{v}}

\newcommand{\diffker}{\bbd}

\newcommand{\vwq}{{\mathbf{q}}}

\newcommand{\ewt}{{\mathbf{\uptheta}}}
\newcommand{\cewt}{\check{\ewt}}
\newcommand{\bff}{{\mathbf{f}}}

\newcommand{\kap}{\kappa}
\newcommand{\kapr}{{\kappa'}}

\newcommand{\tnj}[2]{{t^{#1}_{#2}}}

\usepackage{accents}

\newcommand{\eq}{e_{q}}
\newcommand{\et}{e_{\theta}}

\newcommand{\newComp}{{\tt C}}

\DeclareMathOperator*{\Dom}{Dom}

\newrobustcmd*{\citeallauthors}{%
  \AtNextCite{\AtEachCitekey{\defcounter{maxnames}{999}}}%
  \citeauthor}

\renewcommand{\tnj}{{t^n_j}}
\newcommand{\tnjp}{{t^n_{j+1}}}
\newcommand{\tnjm}{{t^n_{j-1}}}
\newcommand{\diff}[2]{\Delta_{#1,#2}}

\newcommand{\cnstVar}[2]{C^{#1}_{#2}}

\newcommand{\rsdOne}[1]{R_n^{(1,#1)}}

\newcommand{\tnjo}{{t^n_{j_1}}}

\newcommand{\tnjt}{{t^n_{j_2}}}

\newcommand{\tnjs}{{t^n_{j_3}}}

\newcommand{\fvo}{V^{[1]}}
\newcommand{\fvt}{V^{[2]}}

\newcommand{\caliqtor}{\cali_n^{(\ell,m;\sharp)\dagger\dagger}}
\newcommand{\graphQtorSharp}{G^{(\ell,m;\sharp)}}

\newcommand{\graphNnThree}{G_{N^{(3,\ell)}}}
\newcommand{\graphNnOne}{G_{N^{(1)}}}

\newcommand{\GnTwoMain}{G_n^{(\ell_1,\ell_1;2\ell_1-1)\dagger\dagger}}
\newcommand{\dfvo}{V^{[1,1]}}
\newcommand{\dfvt}{V^{[2,1]}}
\newcommand{\fvoo}{V^{[1,1;1]}}
\newcommand{\fvto}{V^{[2,1;1]}}
\newcommand{\fvot}{V^{[1,1;2]}}

\newcommand{\unlconst}[1]{C_{u_n^{(#1)}}}

\newcommand{\pfqtanconst}[2]{C_{#1,#2}^{\eqref{eq:240403.1648}}}
\newcommand{\pfqtorConstOne}[1][\pi]{C_{\ell,m,#1}^{(1)}}
\newcommand{\pfqtorConstTwo}[1][\pi]{C_{\ell,m,#1}^{(2)}}

\newcommand{\gtwoconst}[1]{C_{(#1)}^{\eqref{eq:240607.1332}}}
\newcommand{\qtrhoconst}[1]{C^{(#1)}_{\widehat\rho}}%C^{(\ell_1,\ell_2,m)}

\newcommand{\qtorconst}{C^{(\ell,m;\sharp)}}

\newcommand{\qtorconstTotal}{C_{\qtor}}
\newcommand{\qtorconstRho}[1][\ell,m]{C_{\qtor}^{(#1)}}

\newcommand{\qtanresid}{R^\qtan_n}

\newcommand{\weifuncNOne}{A^{(N_n^{(1)})}_{n,j}}
\newcommand{\kernelNOne}{\diffker^{N^{(1)}}_{n,j}}

\newcommand{\innerProdDiffker}{c_0}

\renewcommand{\calii}{\cali}

\newcommand{\newCompProd}{\widetilde{\tt C}^{(\hv_1,\pi)}}
\newcommand{\newGraph}{{\widetilde{\tt G}}}
\newcommand{\extVertex}{{\widetilde{V}}}
\newcommand{\extKer}{{\tilde{\bbf}}}
\newcommand{\extIndex}{{\tilde{j}}}
\newcommand{\newWeight}{{\tilde{A}}}

\newcommand{\prodFormulaIndex}{\Pi}

\excludecomment{en-text}
\includecomment{scrap}
\excludecomment{scrap}

\includecomment{todelete}
\excludecomment{todelete}

\includecomment{furthertodo}
\excludecomment{furthertodo}

\includecomment{memo}
\allowdisplaybreaks
\title{Asymptotic expansion of 
the weighted power variation 
with second order differences of
a stochastic differential equation driven by fBm%
\footnote{
This work was in part supported by 
Japan Science and Technology Agency CREST JPMJCR14D7, JPMJCR2115; 
Japan Society for the Promotion of Science Grants-in-Aid for Scientific Research 
No. 17H01702 (Scientific Research);
and by a Cooperative Research Program of the Institute of Statistical Mathematics. %%
}
}

\author{Hayate Yamagishi}
\affil{Graduate School of Mathematical Sciences, University of Tokyo
\footnote{Graduate School of Mathematical Sciences, 
University of Tokyo: 3-8-1 Komaba, Meguro-ku, Tokyo 153-8914, Japan. 
e-mail: 
yhayate@ms.u-tokyo.ac.jp
}}
\affil{CREST, Japan Science and Technology Agency%\footnote{}
}
\begin{document}
\maketitle
\begin{abstract}
We study a process satisfying a one-dimensional stochastic differential equation 
driven by fractional Brownian motion with Hurst index $H>1/2$, 
and consider the weighted power variation based on the second order differences of the process.
We derive the asymptotic expansion formula of its distribution based on the 
theory of expansion of Skorohod integrals by \cite{nualart2019asymptotic}.
The formula includes the rate of convergence as a corollary. 
To facilitate the application of the general expansion theory, 
% To implement the general theory of expansion, 
we employ the theory of exponents from \cite{2024Yamagishi-asymptotic} to obtain estimates of functionals.
\end{abstract}

% \keywords{%
{\bf Keywords:}
fractional Brownian motion, 
stochastic differential equation,
weighted power variation,
second-order difference, 
mixed central limit theorem,
% mixed normal distribution,
asymptotic expansion,
random symbol,
Malliavin calculus,
Skorohod integral,
% multiple stochastic integral,
exponent, 
weighted graph,
asymptotic order
% }

\tableofcontents

%%%%%%%%%%%%%%%%%%%%%%%%%%%%%%%%%%%%%%%%%
%%%%%%%%%%%%%%%%%%%%%%%%%%%%%%%%%%%%%%%%%

% \newpage
% \subfile{0-total-intro.tex}

\section{Introduction}
% \paragraph*{{\rb modelの説明}}
% \paragraph*{{(1) SDEの説明}}
% {\myblue 
In this paper,
we consider the following one-dimensional SDE driven by a fractional Brownian motion (fBm)
with its Hurst parameter $H>1/2$:
\begin{align}\label{eq:240626.1939}%{eq:230925.1452}
  dX_t &= V^{[2]}(X_t) dt + V^{[1]}(X_t) dB^H_t
  \\
  X_0&=x_0\in\bbR,
  \nn
\end{align}
where $V^{[i]}$ ($i=1,2$) are functions from $\bbR$ to $\bbR$,
$t\in[0,T]$ and $T>0$ is a fixed terminal time.
To make the presentation simple, we set $T=1$.
The stochastic integral in \eqref{eq:240626.1939} is a pathwise Riemann-Stieltjes
integral (Young integral), 
and it is known that 
there exists a unique solution under some regularity conditions on $V^{[i]}$ ($i=1,2$), 
as Nualart and Rascanu \cite{rascanu2002differential} detailed.
% }

% \paragraph*{{ (2) estimatorの説明}}
Weighted realized power variations of 
diffusion-related processes are a central tool for 
estimation problems about the diffusion coefficient.
When the driving noise process is a classical Brownian motion
(i.e. Wiener process),
the problem has been well studied since 90's,
for example \cite{florens1993estimating},
\cite{genon1993estimation} 
and \cite{2000Jacod-NonparametricKE} to name a few.
In the case where the driving noise is a fractional Brownian motion, 
the error %of convergence 
of the estimator was studied in 
\cite{leon2004stable} and \cite{leon2007limits}.
The asymptotic mixed normality of the weighted power variation 
for the solution to SDE \eqref{eq:240626.1939} was 
obtained in \cite{leon2007limits},
but the Hurst parameter was restricted to $(\half,\frac34)$
since the asymptotic variance diverges to infinity
due to the covariance structure of fBm.
To overcome this problem, \cite{leon2007limits} considered 
the weighted power variation based on the second order differences of the process $X_t$.
Formally, 
for $n\geq2$ and $j=1,...,n-1$,
we denote the second order difference of $X_t$ by
\begin{align*}
  \secDiff{n}{j}X = (X_{t^n_{j+1}} - X_{t^n_{j}}) - (X_{t^n_{j}} - X_{t^n_{j-1}}), 
\end{align*}
with $t^n_{j}= j/n$, and  %for $j=1,..,n$.
define the $2k$-th $(k\geq1)$ weighted power variation $S_n$ as follows:
\begin{align}
  S_n&=
  n^{2kH-1}\sum_{j=1}^{n-1}
  f(X_\tnj) (\secDiff{n}{j}X)^{2k},
  \label{eq:240628.1935}
\end{align}
where the function $f:\bbR\to\bbR$ has enough regularity.
It is known that the functional $S_n$ converges in $L^p$ to 
$S_\infty$ defined by 
\begin{align}
  S_\infty&= \cnstVar{2k}{0}
  \int^1_0 f(X_t) (V^{[1]}(X_t))^{2k} dt
  % \int^1_0 f_t (V^{[1]}_t)^{2k} dt.
  \label{eq:240628.1936}
\end{align}
with the constant $\cnstVar{2k}{0}$ defined at \eqref{eq:240605.1731}.
Since the covariance function of the second order difference of 
a fBm with $H>\half$ decays fast enough,
we can obtain the asymptotic mixed normality of 
the convergence error of $S_n$ to $S_\infty$ 
for all $H>\half$ as \cite{leon2007limits} proved,
though they did not show the rate of convergence.

% \paragraph*{{\rb 漸近展開の説明}}
% \paragraph*{{ (3) 漸近展開の一般論, (4) 漸近展開の具体的な研究について}}
% \paragraph*{{\rb (4) 漸近展開の具体的な研究について}}
% \paragraph*{{\rb 漸近展開の歴史について}}　(位置未定)
In the current paper, we consider the asymptotic expansion
of the distribution of the error $\sqrt{n}(S_n-S_\infty)$
beyond the rate of convergence.
Asymptotic expansion of the distribution of estimators are important in the context of statistical applications,
since they can yield more an accurate approximation than 
that based on the central limit theorem (i.e. normal approximation).
The theory in the case of i.i.d. dates back to the beginning of twentieth century (cf. \cite{edgeworth1905law} and \cite{cramer1928composition}) and 
the extension to stochastic processes with more complicated dependent structure has been considered since 1980's 
(see \cite{gotze1983asymptotic} and \cite{gotze1994asymptotic}).
Recently, 
a great progress on general theories of asymptotic expansion for Wiener functionals
has been achieved by Yoshida and his coauthors, e.g. 
\cite{tudor2019asymptotic}, \cite{nualart2019asymptotic} and 
\cite{tudor2023high}.
They can be applied to the target in this paper,
since indeed it can be written as a functional on a Wiener space.
Examples of the application of \cite{nualart2019asymptotic} are 
weighted quadratic variations of a (classical) Brownian motion
and fractional Brownian motion in \cite{nualart2019asymptotic},
a quadratic variation with anticipative weights and 
robust realized volatility in \cite{yoshida2023asymptotic},
a quadratic variation of SDE driven by fBm with Hurst index $H\in(\half,\frac34)$ in \cite{2024Yamagishi-AsymptoticEO},
a Hurst estimator of fractional SDE with $H>1/2$ in \cite{2024Yamagishi-asymptotic}.
The expansion method of \cite{tudor2019asymptotic} and \cite{tudor2023high} based on the so-called gamma factors are applied to 
quadratic variations of mixed fractional Brownian motions in \cite{tudor2020asymptotic} and 
a Hurst estimator of fBm in \cite{mishura2023asymptotic}.
In \cite{mishura2023asymptotic},
the second-order modification of the estimator are also considered.

% \paragraph*{{\rb (5) exponentについてのremark}}
The general theories of asymptotic expansion mentioned above are so powerful that we can obtain the expansion formulas
only by checking a set of several sufficient conditions.
However, while doing so, we need to
separate functionals of a greater magnitude from negligible ones,
and it is a highly non-trivial task to
estimate the asymptotic order of norms of these functionals,
making difficult the derivation of asymptotic expansion formulas:
As we will see in the following argument,
the number of the functionals appearing there are so large 
that it is virtually impossible 
to estimate the order of them one by one.
Thus it is essential to find 
a method for systematic evaluation of orders of functionals having a certain form to make use of the general theory.
Indeed, similar theories have been developed for 
different targets;
\cite{yoshida2023asymptotic}
for functionals related to Wiener processes (i.e. $H=1/2$),
\cite{yamagishi2023order} in
the context of variation processes of fBm with $H\in(\half,\frac34)$, 
which was used in \cite{2024Yamagishi-AsymptoticEO},
and 
\cite{2024Yamagishi-asymptotic} for functionals related to 
variations using second-order difference of fBm with $H>1/2$.
In this paper, we shall employ that from \cite{2024Yamagishi-asymptotic} for the problem of weighted power variations of SDE driven by fBm.
% As we will see in the following argument,
% the number of the functionals appearing there are so large 
% that it is virtually impossible 
% to estimate the order of them one by one,
% and we can argue that without this tool we could not
% sort out non-negligible functionals.

% \newpage
\subsection{Main result and high-level plan}
% \subsubsection*{{\rb (1) $G_\infty$の定義}}
We denote the error of convergence by 
\begin{align}
  Z_n=\sqrt{n}(S_n-S_\infty)
  \label{eq:240628.1948}
\end{align}
with $S_n$ and $S_\infty$ defined at \eqref{eq:240628.1935} and 
\eqref{eq:240628.1936}, respectively.
We define the functional $G_\infty$ at \eqref{eq:240630.1538}, namely
% We define the functional $G_\infty^{\dagger}$ by
% \begin{align*}
%   % C_{\Delta}&:=
%   % \sum_{\ell_1=1}^k
%   % 2\ell_1\times 
%   % \gtwoconst{\ell_1,\ell_1,2\ell_1-1}\times
%   % \qtrhoconst{\ell_1},\tand 
%   G_\infty^{\dagger}=
%   \int_0^1 a(X_t)^2dt
% \end{align*}
\begin{align}
  % C_{G_\infty} &=   
  % \sum_{\ell_1=1}^k 
  % \gtwoconst{\ell_1,\ell_1,2\ell_1-1} \times 
  % \qtrhoconst{\ell_1}
  % \label{eq:240628.1937}
  % \\
  G_\infty&= 
  C_{G_\infty} \int_0^1 f(X_t)^2 (V^{[1]}(X_t))^{4k} dt
  % C_{G_\infty} \int_0^1 a(X_t)^2 dt
  \label{eq:240628.1951}
\end{align}
with a positive constant $C_{G_\infty}$ 
defined at \eqref{eq:240628.1937}.
% $a(x)=f(x) (V^{[1]}(x))^{2k}$ 
This functional plays the role of conditional asymptotic variance of $Z_n$.

% \paragraph*{{(2) 仮定の条件}}
To derive the asymptotic expansion of the distribution of $Z_n$, 
we assume the following conditions: % on SDE \eqref{eq:240626.1939}. 
\begin{ass}\label{ass:240626.1518}%{ass:230927.1617}
  (i) $V^{[i]}\in C^\infty_b(\bbR)$ for $i=1,2$,
  where $C^k_b(\bbR)$ is the set of 
  $k$ times continuously differentiable bounded functions $\bbR\to\bbR$
  with bounded derivatives of order up to $k$.
  % where $C^\infty_b(\bbR)$ is the set of the smooth functions whose derivatives of any order are bounded together with itself.
  \item[(ii)]  
  $f\in C^\infty_b(\bbR)$.

  \item[(iii)] 
  The functional $G_\infty$ satisfies 
  $(G_\infty)^{-1}\in L^{\infty-}=\cap_{p>1}L^p$.
\end{ass}
\noindent
The third condition is assumed to ensure the nondegeneracy of the distribution of $Z_n$.

% \paragraph*{{\rb 主結果}}
% \paragraph*{{(3) 必要なnotationを説明}}
We introduce some notations to state the main result.
We write $\phi(z;\mu,v)$ for the density of the normal distribution 
with mean $\mu\in\bbR$ and variance $v>0$.
% {\redmy 
For a polynomial random symbol 
$\varsigma(\xi) = \sum_{\alpha\in\bbZ_{\geq0}} \varsigma_\alpha\,\xi^\alpha$ %(finite sum)
% $\varsigma(\xi) = \sum_{\alpha\in\bbZ_{\geq0}} c_\alpha\,\xi^\alpha$ %(finite sum)
with coefficient random variables $\varsigma_\alpha$ and a dummy variable $\xi$,
% with coefficient random variables $c_\alpha$ and a dummy variable $\xi$,
the operation of the adjoint $\varsigma(\partial_z)^*$ on $\phi(z;0,G_\infty)$ under the expectation
is  defined by 
\begin{align*}
  E\sbr{\varsigma(\partial_z)^* \phi(z;0,G_\infty)}
  =\sum_\alpha (-\partial_z)^\alpha E\sbr{\varsigma_\alpha\;\phi(z;0,G_\infty)}.
  % =\sum_\alpha (-\partial_z)^\alpha E\sbr{c_\alpha\phi(z;0,G_\infty)}.
\end{align*}%}
See Section \ref{sec:240626.1819} for details of the definition of random symbols.
We denote by $\cale(M,\gamma)$ 
the set of measurable functions $g:\bbR\to\bbR$
such that 
$\abs{g(z)}\leq M(1+\abs{z})^\gamma$. % for all $z\in\bbR$.

% \paragraph*{{(4) 主結果}}
Our main result is the following asymptotic expansion formula of 
the distribution of the functional $Z_n$.
\begin{theorem}\label{thm:240620.2214} %{prop:240620.2214}
  Let $M,\gamma>0$.
  Suppose that 
  % $(X_t)_{t\in[0,T]}$ is the solution to SDE \eqref{eq:240626.1939}
  % satisfying Assumption \ref{ass:240626.1518} and
  Assumption \ref{ass:240626.1518} is satisfied and
  the functional $Z_n$ is defined at \eqref{eq:240628.1948}.
  %
  % Then, the distribution of the functional $Z_n$ 
  % has the following asymptotic expansion as $n$ goes to $\infty$:
  Then the following estimate holds as $n$ goes to $\infty$:
  \begin{align*}
    \sup_{g\in\cale(M,\gamma)}\abs{E[g(Z_n)]-
    \int_\bbR g(z) p_n(z) dz}
    =o(n^{-\half}),
  \end{align*}
  with the approximate density $p_n$ written as 
  \begin{align*}
    p_n(z)=
    E[\phi(z;0,G_\infty)]
    +n^{-\half}
    E[\mS(\partial_z)^* \phi(z;0,G_\infty)].
    % E\sbr{\mS(\partial_z)^* \phi(z;0,G_\infty)}
  \end{align*}
  Here %{\redmy the asymptotic conditional variance $G_\infty$ of $Z_n$ 
  $G_\infty$ is defined at \eqref{eq:240628.1951}, 
  % {\rb[$G_\infty$はどこで定義される？]} 
  and the random symbol $\mS=\mS^{(3,0)}+\mS^{(1,0)}$ is written 
  with the random symbols 
  $\mS^{(3,0)}$ and $\mS^{(1,0)}$ defined at
  \eqref{eq:240626.1954} and \eqref{eq:240626.1955}, respectively.
\end{theorem}

% \paragraph*{{\rb (5) 証明のhigh-level-picture}}
Here we provide an outline of the proof.
In Section \ref{sec:240626.1840},
we expand the functional $Z_n$ as $Z_n=\delta(u_n)+r_nN_n$,
where $\delta(u_n)$ is a Skorohod integral, and %of $u_n$, 
$r_nN_n$ is the residual term (Proposition \ref{prop:240522.2303}).
In the context of the present paper, we set $r_n=n^{-\half}$.
Next, in Section \ref{sec:240628.1844}, we decompose the functionals 
$\abr{DM_n,u_n}$, $\abr{D\abr{DM_n,u_n},u_n}$, $N_n$ and $D_{u_n} N_n$,
which appear in the subsequent arguments several times,
% repeatedly,
and also identify the asymptotic variance $G_\infty$.
In Section \ref{sec:240628.1845},
we determine
the summands of the limit random symbol $\mS$,
namely $\mS^{(2,0)}_0,\mS^{(3,0)},\mS^{(1,0)}$ and 
$\mS^{(2,0)}_1$.
In Section \ref{sec:240630.2339},
for $u_n,r_n,N_n,G_\infty$ and $\mS$ specified by then, %in the present context,
we will verify the sufficient condition of Theorem \ref{thm:240702.2413}
to apply the general theory of asymptotic expansion of \cite{nualart2019asymptotic}.
The condition consists of four parts related to
differentiability in Malliavin's sense,
asymptotic orders of functionals, random symbols, and 
non-degeneracy of the distributions.
(They are detailed in Section \ref{sec:240626.1819}.)
The condition on regularity is easily verified by the assumption on SDE \eqref{eq:240626.1939}.
(See \cite{nualart2009malliavin} and \cite{hu2016rate}
for the Malliavin differentiability of SDE driven by fBm with $H>1/2$.)
Also the conditions on random symbols are essentially checked 
when finding the limit random symbol in Section \ref{sec:240628.1845}.
For the rest, 
we check the order conditions in Section \ref{sec:240628.1902},
and the non-degeneracy in Section \ref{sec:240628.1903}.
% 

% \paragraph*{{\rb (6) exponentについてのfurther remark}}
The theory of exponent from \cite{2024Yamagishi-asymptotic}
will be heavily used in the above argument as hinted in the introduction.
Although the method was originally developed for a weighted quadratic variation of fBm based on second-order differences,
it was written for much more general cases than needed there, allowing us to 
apply the main theorem of order estimation
(Theorem \ref{thm:240628.1910})
to the context of this paper.
On the other hand, we need to estimate the effect of 
the operator $D_{v_n}$ on the order of functionals,
where $v_n$ is a summand of $u_n$,
and in Proposition \ref{prop:240617.2009},
we extend a related result from \cite{2024Yamagishi-asymptotic} 
to the case in this paper.
%  with $v_n^{(q)}$ written as 
% \begin{align*}
%   v_n^{(q)}&=
%   % u_n^{(q)}
%   n^{q H-\half}
%   \sum_{{\redmy j_0\in[n-1]}}
%   A'_{n,j_0}
%   I_{q-1}({\diffker^n_{j_0}}^{\otimes q-1})
%   \diffker^n_{j_0}
%   % \label{eq:240617.1534}
% \end{align*}
% 
% To highlight the applicability of the method of exponent, 

% \paragraph*{{(7) organization}}
The organization of the paper is as follows.
% The rest of this paper is organized in the following way:
% 
Section \ref{sec:240626.1456} provides some preliminary results used in this paper.
Namely, 
we provide an overview of the general theory by \Citeauthor{nualart2019asymptotic} 
% \cite{nualart2019asymptotic}
of the asymptotic expansion of
the distributions of perturbed Skorohod integrals 
in Section \ref{sec:240626.1819}, 
and the basic facts of fractional Brownian motions 
in Section \ref{sec:240626.1825}.
In Section \ref{sec:240626.1834},
we review the theory of exponent developed in \cite{2024Yamagishi-asymptotic}.
Particularly, 
the aforementioned extension of a result from \cite{2024Yamagishi-asymptotic}
will be proved in Section \ref{sec:240626.1835}.
Section \ref{sec:240617.1526} treats the main argument of the proof of Theorem \ref{thm:240620.2214},
% In Section \ref{sec:240626.1840},
while we defer technical lemmas to Section \ref{sec:240626.1842};
% Precisely,
we prove lemmas concerned with the stochastic expansion of $Z_n$ 
(Proposition \ref{prop:240522.2303}) in Section \ref{sec:240626.1850},
and those to identify the limit of the functionals $(D_{u_n})^iM_n$ ($i=2,3$)
in Section \ref{sec:240607.1650}.

% \newpage
\subsection{Notations}
The following notations are repeatedly used in the following sections.
\begin{itemize}
  \setlength{\parskip}{0cm} \setlength{\itemsep}{5pt}

  \item We denote by $\bbZ$, $\bbZ_{\geq0}$ and $\bbN$
  the set of integers, that of nonnegative integers and
  that of positive integers.
  We write $\ntwo=\cbr{j\in\bbZ \mid j\geq2}$ for notational convenience.
  For $n\in\bbN$, we write $[n]=\cbr{1,...,n}$ for short.

  \item For a (finite) set $S$, we denote 
  the cardinality of $S$ by $\abs{S}$.

  \item 
  % $C^\infty_b(\bbR)$: the set of the smooth functions 
  % whose derivatives of any order are bounded together with itself.
  For $k\in\bbN\cup\cbr{\infty}$, $C^k_b(\bbR)$ denotes
  the space of $k$ times continuously differentiable functions $\bbR\to\bbR$
  which are bounded together with its derivatives of order up to $k$.

  \item The functions $V^{[i]}:\bbR\to\bbR$ ($i=1,2$) are
  the diffusion coefficient and drift coefficient
  of SDE \eqref{eq:240626.1939}, respectively.
  We denote the $k$-th derivative of $V^{[i]}$ by $V^{[i;k]}$.
  
  % \item The $\bbR$-valued functions $V^{[1]}$ and $V^{[2]}$ 
  % defined on $\bbR$ appears in SDE \eqref{eq:230925.1452} as 
  % the diffusion coefficient and drift coefficient, respectively, and 
  % $V^{[i;k]}$ is the $k$-th derivative of $V^{[i]}$.

  \item For a function $g:\bbR\to\bbR$, we write $g_t=g(X_t)$ for brevity,
  where $X_t$ is the solution to SDE \eqref{eq:240626.1939}.
  In particular, we often write $V^{[i]}_t$ for the coefficient $V^{[i]}(X_t)$ of SDE \eqref{eq:240626.1939}.

  \item 
  We write $t^n_j=j/n$ for $j=0,...,n$ and $n\in\bbN$.
  % However, when there is no risk of confusion, we use the notation $t_j$ to represent $t^n_j=j/n$.
  % We denote the indicator function of the interval $\sbr{\tjm,\tj}$
  % by $1_{n,j}$ or $1_j$.
  
  \item 
  We denote by $L^p$ the $L^p$-space of
  random variables on the probability space fixed in the following arguments.
  % When we simply write $L^p$, it means the $L^p$-space of
  % random variables on the probability space.
  We also write
  $L^{1+}=\cup_{p>1} L^p$ and $L^{\infty-}=\cap_{p>1}L^p$.
  The $L^p$-norm of a random variable is denoted by 
  $\norm{.}_p$ or $\norm{.}_{L^p}$.
  % We also write $\norm{.}_{L^p(P)}$ 
  % to distinguish it from the $L^p$-norm of functions on $[0,T]^k$, 
  % which is denoted by $\norm{.}_{L^p([0,T]^k)}$.
  We denote the Sobolev norm in Malliavin calculus by 
  $\norm{.}_{i,p}$ for $i\in\bbN$ and $p\geq1$.
  
  \item 
  Let  $\alpha\in\bbR$ and 
  $(F_n)_{n\in\bbN}$ be a sequence of random variables.
  We use the following notation to express asymptotic orders.
  \begin{itemize}
    \item [\labelitemi]
    $F_n=O_{L^{p}}(n^\alpha)$ 
    [resp. $F_n=o_{L^{p}}(n^\alpha)$],\quad
    if $\norm{F_n}_{p}=O(n^\alpha)$
    [resp. $\norm{F_n}_{p}=o(n^\alpha)$]\; for $p\geq1$.
    \item [\labelitemi]
    $F_n=O_{L^{\infty-}}(n^\alpha)$ 
    [resp. $F_n=o_{L^{\infty-}}(n^\alpha)$],\quad
    if $F_n=O_{L^{p}}(n^\alpha)$
    [resp. $F_n=o_{L^{p}}(n^\alpha)$]\; for every $p>1$.
  \end{itemize}
  For $(F_n)_{n\in\bbN}$ with $F_n\in\bbD^\infty$, we write
  \begin{itemize}
    \item [\labelitemi]
    $F_n=O_M(n^\alpha)$ [resp. $F_n=o_M(n^\alpha)]$, \quad
    if $\norm{F_n}_{i,p}=O(n^\alpha)$
    [resp. $\norm{F_n}_{i,p}=o(n^\alpha)$]
    for every $i\in\bbN$ and $p>1$.
  \end{itemize}
  When $F_n=O_{L^{p}}(n^{\alpha+\epsilon})$ for any $\epsilon>0$, we write 
  $F_n=\hat O_{L^p}(n^\alpha)$.
  For the other norms, we use the similar notations.

  % \item We use the notation $\norm{.}_\beta$
  % for the $\beta$-H\"{o}lder seminorm of 
  % $\beta$-H\"{o}lder continuous functions on $[0,T]$.
  % We use the notation $\norm{.}_\infty$ for the uniform norm.
  % For the seminorms restricted on an interval $[t,t']$,
  % we write $\norm{.}_{t,t',\beta}$ and $\norm{.}_{t,t',\infty}$.
  % See Section \ref{220405.1158} for their formal definitions.
\end{itemize}

% \newpage
\section{Preliminaries}\label{sec:240626.1456}

% \newpage
\subsection{Asymptotic expansion of Skorohod integrals}
\label{sec:240626.1819}% \label{230901.1140}
We review the theory of asymptotic expansion of Skorohod integrals 
described in \cite{nualart2019asymptotic} in the one-dimensional case.
Let $(\Omega, \calf, P)$ be a complete probability space equipped with an isonormal Gaussian process 
$W=\cbr{W(h)_{h\in\calh}}$ 
on a separable real Hilbert space $\calh$.
We denote the Malliavin derivative operator by $D$ and 
its adjoint operator, namely the divergence operator or the Skorohod integral, by $\delta$.
For $p\geq1$, $k\in\bbN$ and a real separable Hilbert space $V$, 
we write $\bbD^{k,p}(V)$ for the Sobolev space of $V$-valued random variables 
that have the Malliavin derivatives up to $k$-th order which have finite moments of order $p$.
We write $\bbD^{k,p}=\bbD^{k,p}(\bbR)$, 
$\bbD^{k,\infty}(V)=\cap_{p>1}\bbD^{k,p}(V)$ and
$\bbD^{\infty}(V)=\cap_{p>1,k\geq1}\bbD^{k,p}(V)$.
The $(k,p)$-Sobolev norm is denoted by $\norm{\cdot}_{k,p}$.
We refer to the monograph \cite{nualart2006malliavin} for a detailed account on this subject.

Consider a sequence of random variables $Z_n$ defined on the probability space 
$(\Omega, \calf, P)$ written as 
\begin{align}\label{eq:231004.1810}
  Z_n=M_n+r_nN_n,
\end{align}
where
$M_n=\delta(u_n)$ is the Skorohod integral of an $\calh$-valued random variable
$u_n\in\Dom(\delta)$,
$N_n$ is a random variable and 
$(r_n)_{n\in\bbN}$ is a sequence of positive numbers such that 
$\lim_{n\to\infty}r_n=0$.
%
% {\mygreen
The variable $Z_n$ of \eqref{eq:231004.1810} is a perturbation of $M_n$ when $N_n=O_p(1)$ as $n\to\infty$.
Such a perturbation always appears when one derives 
a stochastic expansion of a statistic $Z_n$ around a principal part $M_n$ that is easy to handle by a limit theorem. 
%
% In the case of the quadratic variation 
% $\bbV_n =n^{2H-1}\sum_{j=1}^n (\Delta_jX)^2$ of (\ref{220420.1130}), 
% the scaled variable $Z_n =n^{1/2}\big(\bbV_n-\bbV_\infty)$ admits the stochastic expansion (\ref{202308260454}) 
% with 
% $r_n = n^{2H-3/2}$ and 
% $M_n=\delta(u_n)$ for $u_n$ of (\ref{220404.1632}): 
% $u_n = n^{2H-1/2} \sum_{j\in[n]} a_{t_{j-1}}  I_1(1_j) 1_j$. 
% The variable $N_n$ has an expression given in (\ref{220404.1633}). 
% }
%

We consider an asymptotic expansion of the distribution of $Z_n$ in the situation where 
$M_n$ stably converges to 
a mixed normal distribution, that is, 
$M_n\overset{d_s}{\to}M_\infty=G_\infty^{1/2}\zeta$ as $n\to\infty$,
where $G_\infty$ is a positive random variable 
and $\zeta$ is a standard Gaussian random variable independent of $\calf$. 
The variable $\zeta$ is given on an extension of the probability space $(\Omega,\calf,P)$. 
Here the stable convergence means that the convergence $(M_n,Y)\to^d(M_\infty,Y)$ holds for any 
random variable $Y$ measurable with respect to $\sigma[W]$, the $\sigma$-field generated 
by the isonormal Gaussian process $W$.
%
% We consider a stable convergence of $Z_n$ to 
% a mixed normal distribution $G_\infty^{1/2}\zeta$,
% where $G_\infty$ is a positive random variable 
% and $\zeta$ is a standard normal distribution independent of $\calf$.}

Nualart and Yoshida \cite{nualart2019asymptotic} introduced the following random symbols.
To write random symbols, we use a simplified notation  
due to the one-dimensional setting.
%We write $\mS(\xi)=\sum_{k}\chi_k\,\xi^k$ (a finite sum) for a polynomial random symbol,
% where $\chi_k$ are coefficient random variables
% %\redb{$\chi_k\,\xi^k$ is the random symbol of degree $k$ with its coefficient random variable $\chi_k$}
%  and $\xi$ stands for a dummy variable.
We denote
$D_{u_n}F=\abr{DF,u_n}_\calh$ for a random variable $F$ regular enough.
The quasi-tangent is defined by 
\begin{align*}
  \qtan_n{[\sfi\sfz]^2}=
  r_n^{-1}\brbr{\babr{DM_n{[\sfi\sfz]}, u_n{[\sfi\sfz]}}_\calh
  -G_\infty{[\sfi\sfz]^2}}
  =r_n^{-1}\brbr{D_{u_n}M_n-G_\infty}{[\sfi\sfz]^2}.
\end{align*}
The quasi-torsion and modified quasi-torsion are defined by 
\begin{align*}
  \qtor_n{[\sfi\sfz]^3}&=
  r_n^{-1}\Babr{D\babr{DM_n{[\sfi\sfz]}, u_n{[\sfi\sfz]}}_\calh,
  u_n{[\sfi\sfz]}}_\calh
  =r_n^{-1}(D_{u_n})^2M_n{[\sfi\sfz]^3}
  \\
  \mqtor_n{[\sfi\sfz]^3}&=
  r_n^{-1}\babr{DG_\infty{[\sfi\sfz]^2},u_n{[\sfi\sfz]}}_\calh
  =r_n^{-1}D_{u_n}G_\infty{[\sfi\sfz]^3},
\end{align*}
respectively.

We write
\begin{align*}
  G^{(2)}_n
  %G^{(2)}_n(\sfz) 
  &%= r_n\qtan_n
  = D_{u_n}M_n - G_\infty
  \tand% \\
  G^{(3)}_n
  %G^{(3)}_n(\sfz) 
  % &%= r_n\mqtor_n 
  = D_{u_n}G_\infty,
\end{align*}
and define the following random symbols
\begin{align}
  \mS^{(3,0)}_n(\tti\sfz)&=
  \frac13\qtor_n[\tti\sfz]^3
  =\frac13r_n^{-1}(D_{u_n})^2 M_n[\tti\sfz]^3
  \label{eq:240629.2341}
  \\
  \mS^{(2,0)}_{0,n}(\tti\sfz)&
  % \frac12 r_n^{-1}G^{(2)}_n(\sfz)=
  =\frac12 \qtan_n[\tti\sfz]^2
  =\frac12 r_n^{-1}\brbr{D_{u_n} M_n - G_\infty}[\tti\sfz]^2
  % =\frac12 r_n^{-1}\brbr{D_{u_n[\tti\sfz]} M_n[\tti\sfz] - G_\infty[\tti\sfz]^2}
  \label{eq:240629.2342}
  \\
  \mS^{(1,0)}_n(\tti\sfz)&=N_n[\tti\sfz]
  \label{eq:240629.2343}
  \\
  \mS^{(2,0)}_{1,n}(\tti\sfz)&=
  D_{u_n} N_n[\tti\sfz]^2
  \label{eq:240629.2344}
  % D_{u_n[\tti\sfz]} N_n[\tti\sfz]
\end{align}
for $\tti\sfz\in\tti\bbR$.
%
% We consider random symbols 
We denote by
$\mS^{(3,0)}$,
$\mS^{(2,0)}_{0}$,
$\mS^{(1,0)}$ and
$\mS^{(2,0)}_{1}$
% These work as 
the limits of the above random symbols
$\mS^{(3,0)}_n$,
$\mS^{(2,0)}_{0,n}$,
$\mS^{(1,0)}_n$ and
$\mS^{(2,0)}_{1,n}$, respectively.
(The meaning of the limit is explained in the condition {\bf[D]} (iii).)
Let
\begin{align*}%\label{eq:230810.1905}
  \Psi(\sfz)=
  \exp\rbr{2^{-1} G_\infty[\tti\sfz]^2}
  \yeq \exp\rbr{-2^{-1} G_\infty\>\sfz^2}
\end{align*}
for $\sfz\in\bbR$.
The following set of conditions {\bf [D]} from Nualart and Yoshida \cite{nualart2019asymptotic}
is a sufficient condition to validate the asymptotic expansion of $Z_n$ of \eqref{eq:231004.1810}.
The parameter $l$ about differentiability below is $l=9$ in this case.
For a one-dimensional functional $F$, we write 
$\Delta_F %=\sigma_F
=\det(\abr{DF,DF}_\calh)=\abr{DF,DF}_\calh=\norm{DF}_\calh^2$
%for the Malliavin covariance (matrix) of $F$,
for the Malliavin covariance of $F$.

\begin{itemize}
  \item [{\bf [D]}]
  \begin{itemize}
    \item [(i)]
    $u_n\in\bbD^{l+1,\infty}(\mH)$,
    %$u_n\in\bbD^{l+1,\infty}(\mH\otimes\bbR^\sfd)$,
    $G_\infty \in \bbD^{l+1,\infty}(\bbR_+)$,
    %$G_\infty \in \bbD^{l+1,\infty}(\bbR^\sfd \otimes_+ \bbR^\sfd)$,
    %$W_n, 
    and
    $N_n\in\bbD^{l,\infty}$.
    %$N_n\in\bbD^{l,\infty}(\bbR^\sfd)$,
    %$W_\infty\in\bbD^{l\vee\sfd_2,\infty}(\bbR^\sfd)$,
    %$X_n\in\bbD^{l,\infty}(\bbR^{\sfd_1})$,
    %$X_\infty\in\bbD^{l\vee(\sfd_2+1),\infty}(\bbR^{\sfd_1})$,
  
    \item [(ii)]
    There exists a positive constant $\kappa$ such that the following estimates hold for every $p>1$:
    \begin{align}
      &\norm{u_n}_{l,p}=O(1)
      \label{220215.1241}\\
      &\bnorm{G^{(2)}_n}_{l-2,p}=O(r_n)
      \label{220215.1242}\\
      &\bnorm{G^{(3)}_n}_{l-2,p}=O(r_n)
      \label{220215.1243}\\
      &\bnorm{D_{u_n}G^{(3)}_n}_{l-1,p}=O(r_n^{1+\kappa})
      \label{220215.1244}\\
      &\bnorm{D^2_{u_n}G^{(2)}_n}_{l-3,p}=O(r_n^{1+\kappa})
      \label{220215.1245}\\
      &\norm{N_n}_{l-1,p}=O(1)
      \label{220215.1246}\\
      &\bnorm{D^2_{u_n} N_n}_{l-2,p}=O(r_n^{\kappa})
      \label{220215.1247}
    \end{align}
  
    \item [(iii)]
    For each pair 
    $(\mT_n,\mT)= (\mS_n^{(3,0)},\mS^{(3,0)})$,
    $(\mS_{0,n}^{(2,0)},\mS_0^{(2,0)})$,
    $(\mS_n^{(1,0)},\mS^{(1,0)})$ and 
    $(\mS_{1,n}^{(2,0)},\mS_1^{(2,0)})$,
    the following conditions are satisfied.
    \begin{itemize}
      \item [(a)]
      $\mT$ is a polynomial random symbol the coefficients of which are in
      $L^{1+}=\cup_{p>1} L^p$.
  
      \item [(b)]
      For some $p>1$, there exists a polynomial random symbol $\bar\mT_n$ such that
      the coefficients of $\bar\mT_n$ belongs to $L^p$,
      the equation
      $E \sbr{\Psi (\sfz)\mT_n(\tti\sfz)} = E[\Psi (\sfz)\bar\mT_n(\tti\sfz)]$
      holds for $\sfz\in\bbR$, and
      the convergence $\bar\mT_n\to\mT$ in $L^p$ holds.
    \end{itemize}
  
    \item [(iv)]
    \begin{itemize}
      \item [(a)] $G^{-1}_{\infty}\in L^{\infty-}$
      \item [(b)]
      There exist 
      $\kappa'>0$ such that
      \begin{align*}
      P[\Delta_{M_n}<s_n]=O(r_n^{1+\kappa'})
      \end{align*}
      for some positive random variables $s_n\in\bbD^{l-2,\infty}$ satisfying 
      $\sup _{n\in\bbN}(\norm{s_n^{-1}}_p + \norm{s_n}_{l-2,p})<\infty$
      for every $p>1$.
    \end{itemize}
  \end{itemize}
\end{itemize}
Writing 
$\bar\mT_n(\sfi\sfz)=\sum_{k\geq0}c_n^k[\sfi\sfz]^k$ and 
$\mT(\sfi\sfz)=\sum_{k\geq0}c^k[\sfi\sfz]^k$,
the convergence $\bar\mT_n\to\mT$ in $L^p$ means that
there exists $k_0\in\bbZ_{\geq0}$ such that 
$c_n^k$ and $c^k$ are all zero for $k>k_0$, and
$c_n^k\to c^k$ in $L^p$ for every $k\leq k_0$.
%
\begin{comment}
  {\myred [なにか説明を入れるか？ quad-varの文章を参照せよとjournal投稿版に書くかな？]
  The functional $s_n$ in (iv) (b) is used to make
  a truncation functional to gain local asymptotic non-degeneracy.
  See Section 7 of \cite{nualart2019asymptotic} for a construction of a truncation function.}
\end{comment}

We write $\phi(z;\mu,v)$ for the density function of the normal distibution 
with mean $\mu\in\bbR$ and variance $v>0$.
For a (polynomial) random symbol 
$\varsigma(\sfi\sfz) = \sum_{\alpha\in\bbZ_{\geq0}} \varsigma_\alpha[\sfi\sfz]^\alpha$
% $\varsigma(\sfi\sfz) = \sum_\alpha c_\alpha[\sfi\sfz]^\alpha$
% $\varsigma(\xi) = \sum_\alpha c_\alpha\,\xi^\alpha$
with random variables $\varsigma_\alpha$,
% with a random variable $c_\alpha$ and $\alpha\in\bbZ_{\geq0}$,
% where $\xi$ is a dummy variable, 
the action of the adjoint $\varsigma(\partial_z)^*$ to $\phi(z;0,G_\infty)$ under the expectation
is defined by 
\begin{align*}
  E\sbr{\varsigma(\partial_z)^* \phi(z;0,G_\infty)}
  =\sum_\alpha (-\partial_z)^\alpha E\sbr{\varsigma_\alpha\phi(z;0,G_\infty)}.
  % =\sum_\alpha (-\partial_z)^\alpha E\sbr{c_\alpha\phi(z;0,G_\infty)}.
  % \label{eq:230927.1623}
  %{220401.1952}
\end{align*} 
Define the random symbol 
$\mS_n = 1 +r_n\mS$ with
\begin{align*}
  \mS(\tti\sfz) = 
  \mS^{(3,0)}(\tti\sfz) + \mS^{(2,0)}_{0}(\tti\sfz) +
  \mS^{(1,0)}(\tti\sfz) + \mS^{(2,0)}_{1}(\tti\sfz),
\end{align*}
and the approximate density 
$p_n(z)$ by
% $\hat p_n(z)$ by
\begin{align*}
  p_n(z) = 
  % \hat p_n(z) = 
  E\sbr{\mS_n(\partial_z)^* \phi(z;0,G_\infty)}.
\end{align*}
% where the action of the adjoint of a random symbol is defined at (\ref{220401.1952}).
%
For $M,\gamma>0$, we denote by
$\cale(M,\gamma)$ 
% $\hat\cale(M,\gamma)$
the set of measurable functions $g:\bbR\to\bbR$
such that 
$\abs{g(z)}\leq M(1+\abs{z})^\gamma$. % for all $z\in\bbR$.
The following theorem rephrases Theorem 7.7 of \cite{nualart2019asymptotic}.
\begin{theorem}
  [Theorem 7.7 of \cite{nualart2019asymptotic}]\label{thm:240702.2413}
  Suppose that Condition {\bf [D]} is satisfied.
  Then, for each $M,\gamma>0$,
  it holds that
  \begin{align*}
    \sup_{g\in\cale(M,\gamma)}
    % \sup_{f\in\hat\cale(M,\gamma)}
    \abs{E\sbr{g(Z_n)}- \int_{\bbR}g(z) p_n(z)dz}
    % \abs{E\sbr{f(Z_n)}- \int_{\bbR}f(z)\hat p_n(z)dz}
    =o(r_n)
  \end{align*}
  as $n\to\infty$.
\end{theorem}

% \newpage
\subsection{Fractional Brownian motion}%\label{sec:231004.1910}
\label{sec:240626.1825}
We review the definitions and notations about a fractional Brownian motion and the related Malliavin calculus.
Although we only consider the case $T=1$ in the following sections,
we give an exposition for a general $T>0$.
Let $B=(B_t)_{t\in[0,T]}$ a fractional Brownian motion 
% Let $B=\cbr{B_t\mid t\in[0,T]}$ a fractional Brownian motion 
with Hurst parameter $H\in(1/2,1)$ %$H\in\rbr{1/2,3/4}$
defined on some complete probability space $(\Omega,\calf, P)$.
An inner product on the set $\cale$ of step functions on $[0,T]$ is defined by 
\begin{align*}
  \abr{1_{[0,t]},1_{[0,s]}}_\calh = E\sbr{B_s\,B_t}
  = \frac12\rbr{\abs{t}^{2H}+\abs{s}^{2H}-\abs{t-s}^{2H}}.
\end{align*}
and $\calh$ is the closure of $\cale$ with respect to
$\snorm{\cdot}_\calh=\abr{\cdot,\cdot}_\calh^{1/2}$.
The map $\cale\ni1_{[0,t]}\mapsto B_t\in L^2(\Omega, \calf, P)$
can be extended linear-isometrically to $\calh$.
We denote this map by $\phi\mapsto B(\phi)$ and 
the process $\cbr{B(\phi), \phi\in \calh}$ is an isonormal Gaussian process. 
In the following sections, tools from Malliavin calculus are based on 
this isonormal Gaussian process.
We write $I_q(h)$ for the $q$-th multiple stochastic integral of $h\in\calh^{\otimes q}$.

It is known that the Hilbert space $\calh$ contains 
not only measurable functions on $[0,T]$ but also distributions as its elements;
see \cite{pipiras2000integration} and \cite{pipiras2001classes} for detailed accounts.
As a subspace of the Hilbert space $\calh$,
we have the linear space $\abs\calh$ of measurable functions $\phi:[0,T]\to\bbR$ such that
\begin{align*}
  \int_0^T\int_0^T \abs{\phi(s)}\abs{\phi(t)}\abs{t-s}^{2H-2}dsdt<\infty,
\end{align*}
and for $\phi$ and $\psi$ in $\abs\calh$, the inner product of $\calh$ is written as 
\begin{align*}
  \abr{\phi,\psi}_\calh=
  \constInProd
  \int^T_0\int^T_0\phi(s)\psi(t)\abs{t-s}^{2H-2}dsdt,
\end{align*}
where $\constInProd=H(2H-1)$.
For a measurable function $\phi$ on $[0,T]^l$, we define
\begin{align*}
  \norm\phi_{\abs\calh^{\otimes l}}^2
  =\constInProd^l \int_{[0,T]^{l}} \int_{[0,T]^{l}} \abs{\phi(u)} \abs{\phi(v)}
  \abs{u_1-v_1}^{2H-2}...\abs{u_{l}-v_{l}}^{2H-2}dudv,
\end{align*}
with $u=(u_1,...,u_l)$ and $v=(v_1,...,v_l)$.
The space
$\abs\calh^{\otimes l}=\bcbr{\phi:[0,T]^l\to\bbR\mid \norm\phi_{\abs\calh^{\otimes l}}<\infty}$ 
of measurable functions
forms a subspace of the $l$-fold tensor product space $\calh^{\otimes l}$ of $\calh$.
We often drop $\calh$ from the notation $\abr{\cdot,\cdot}_\calh$, and 
instead write $\abr{\cdot,\cdot}$ for brevity, when there is no risk of confusion.
We refer to \cite{nualart2006malliavin} or 
\cite{nourdin2012selected} for a detailed account on fBm. 
We collect some basic estimates related to SDE's driven by a fBm with Hurst index $H>1/2$
in Section \ref{sec:240626.1842}.

Recall that we write 
$t^n_j=j/n$ for $n\in\bbN$ and $j\in\cbr{0,...,n}$.
We denote by $\bbone^n_j$ the indicator function $\bbone_{[t^n_{j-1},t^n_j]}$
of the interval $[t^n_{j-1},t^n_j]$,
and we write $\diffker^n_j=\bbone^n_{j+1}-\bbone^n_j$.
We regard $\bbone^n_j$ and $\diffker^n_j$ as an element of $\abs{\calh}\subset\calh$.
%%%
%%%
% {\mygreen 
% We write $[n]=\cbr{1,..,n}$.
% Define $\rho_H(k)$, $c_H$ and 
% $\beta_n\rbr{j_1,j_2}$ by %=\beta_{j_1,j_2}
%  %:= \abr{1_{\jon}, 1_\jtw}_\calh$  for $j_1,j_2\in[n]$ by
% \begin{align}
%   \rho_H(k) &= \frac12 \rbr{\abs{k+1}^{2H} +\abs{k-1}^{2H} -2\abs{k}^{2H}}
%   =\abr{1_{[0,1]}, 1_{[k,k+1]}}_\calh
%   &&\tfornsp k\in\bbZ,
%   \nn\\
%   c_H^2 &%= c_{H,(\ref{210417.1805})}^2
%   =\sum_{k\in\bbZ} \rho_H(k)^2 ,
%   \label{210417.1805}
%   \\
%   \beta_n\rbr{j_1,j_2}=\beta_{j_1,j_2}%=\beta_{j_1,j_2;T}
%   &= \abr{1_{\jon}, 1_\jtw}_\calh
%   = T^{2H} n^{-2H} \rho_H(\jon-\jtw)
%   &&\tfornsp \jon,\jtw\in[n],
%   \nn%\label{220422.1911}
% \end{align}
% respectively.}
%%%
%%%
We define the constant $\innerProdDiffker$ by 
% {\redmy We define a constant $c_{2,H}$ by 
\begin{align}
  % \bbone^n_j&=\bbone_{[t^n_{j-1},t^n_j]}
  % \label{230925.1608}
  % \\
  % \diffker^n_j&=\bbone^n_\jp - \bbone^n_j
  % \label{230925.1609}
  % \\
  \innerProdDiffker&=4-2^{2H}.
  % =\abr{\diffker^1_j,\diffker^1_j}
  \label{230925.1607}
\end{align}
Note that the relation 
$n^{-2H}\innerProdDiffker=\babr{\diffker^n_j,\diffker^n_j}_\calh$
% $\innerProdDiffker=n^{2H}\babr{\diffker^n_j,\diffker^n_j}_\calh$
holds.

\subsection{Theory of exponent}\label{sec:240626.1834}

Here we recap the theory of exponent developed in \cite{2024Yamagishi-asymptotic}
to estimate the order of norms of functionals appearing in Section \ref{sec:240617.1526},
such as \eqref{eq:240607.1334}, \eqref{eq:240610.1145} and \eqref{eq:240629.2414}.
% 
% In contrast to the context of {\rb\cite{???}},
% the functionals treated in this paper is written as the sum on $[n-1]^k$
First we review the definitions and notations of the theory,
and then give the main proposition in \cite{2024Yamagishi-asymptotic}
% for the estimation of order of functionals 
(Theorem \ref{thm:240628.1910}).
Then Proposition \ref{prop:240617.2009} explains in terms of exponents
the change of the order of a functional by the operator $D_{v_n}$,
where $v_n$ is a $\calh$-valued random variable of a certain form.
This proposition is a generalization of Proposition 4.6 of \cite{2024Yamagishi-asymptotic}.

\subsubsection{Order estimate using the theory of exponent from \cite{2024Yamagishi-asymptotic}}
% {\rb Definitions of weighted graph, exponent and related functional [タイトル直す]}}
%%%%%%%%%%%%%%%%%%%%%%%%%%%%%%
%%%%% underlying sets %%%%%%%%
%%%%%%%%%%%%%%%%%%%%%%%%%%%%%%
% \paragraph{{\rb Underlying sets}}
First we introduce some underlying sets on which 
the weighted graphs and the theory of exponent are built.
%
%
%%%%% set of pairs %%%%%
For a nonempty finite set $V$, 
we denote by $p(V)$
the set of subsets of $V$ having two distinct elements.
We will write $[v,v']$ for 
the element $\cbr{v,v'}$ of $p(V)$. 
%%%%% set of `expanded' vertices $\whv$ %%%%%
We write
$\whv=V\times\cbr{1,2}$ and 
denote elements of $\whv$ by $\hv, \hv', \hv_0, \hv_1...$.
\begin{comment}
  {\myred [1]ここにもう少しfunctionalとの対応についてのコメントを書いても良い．}
\end{comment}
%
%%%%% expanded pairs %%%%%
\begin{comment}
  {\myred [2]ここにもう少しfunctionalとの対応についてのコメントを書いても良い．}
\end{comment}
We define $\hp(V)$ to be
the set of subsets 
$\cbr{(v,\kap),(v',\kap')}$
of $\whv=V\times \cbr{1,2}$ such that 
$v\neq v'(\in V)$.
Again we write $[(v,\kap),(v',\kap')]$ or $[\hv,\hv']$ for elements of $\hp(V)$.
When $V$ is a singleton, i.e. $V=\cbr{v_0}$ with some $v_0$,
we define $p(V)=\emptyset$ and also $\hp(V)=\emptyset$.

%%%%%%%%%%%%%%%%%%%%%%%%%%%%%%
%%%%% weighted graph %%%%%%%%%
%%%%%%%%%%%%%%%%%%%%%%%%%%%%%%
% \paragraph{{\rb Weighted graph}}
Let $V$ be a nonempty finite set.
Given functions 
$\ewt:\hp(V)\to\bbZ_{\geq0}$ and 
$\vwq:\whv\to\bbZ_{\geq0}$,
we call the triplet $(V,\ewt,\vwq)$ {\it a weighted graph}.
When $V$ is a singleton, the domain of $\ewt$ is $\hp(V)=\emptyset$, and 
$\ewt$ is the function from the empty set to $\bbZ_{\geq0}$,
which we denote by $0$.
For a weighted graph $G=(V,\ewt,\vwq)$, 
we call $V$, $\ewt$ and $\vwq$ 
{\it the set of vertices, 
the weights on edges and % the edge weight (function) and 
the weights on vertices}, % the vertex weight (function), 
respectively.
We denote the set of vertices of $G$ (i.e. $V$) by $V(G)$, and
write $\whv(G)=V(G)\times\cbr{1,2}$.
%
%%%%% concatenation of multiple weighted graphs %%%%%
For two weighted graphs $G=(V,\ewt,\vwq)$ and $G'=(V',\ewt',\vwq')$
such that $V\cap V'=\emptyset$,
we define the combined weighted graph 
% we define a weighted graph 
$\widetilde{G}=(\extv,\tilde\ewt,\tilde\vwq)$ 
by
$\extv=V\sqcup V'$,
\begin{align*}
  &\tilde\ewt|_{\hp(V)} = \ewt, \quad
  \tilde\ewt|_{\hp(V')} = \ewt'\tand
  \tilde\ewt|_{\hp(V\sqcup V')\setminus (\hp(V)\sqcup\hp(V'))} = 0
  % \tilde\ewt|_{\hp(V,V')} = 0\tand
  \\&
  \tilde\vwq|_{\whv} = \vwq,\quad
  \tilde\vwq|_{\whv'} = \vwq',
\end{align*}
and denote this weighted graph $\widetilde{G}$ by $G\vee G'$.
In general, we write 
\begin{align*}
  \mathop\vee_{k'\in\cbr{1,...,k}} G_{k'} = G_1 \vee ... \vee G_{k}
\end{align*}
for $k\in\bbN$ and weighted graphs $G_1,...,G_{k}$ whose sets of vertices are disjoint.

%%%%%%%%%%%%%%%%%%%%%%%%%%%%%%%%%%%%%%%
%%%%% connectedness, component %%%%%%%%
%%%%%%%%%%%%%%%%%%%%%%%%%%%%%%%%%%%%%%%
% \paragraph{{\rb Connectedness and component}}
%%%%% projected edge weight %%%%%
For a nonempty finite set $V$ and 
$\ewt:\hp(V)\to\bbZ_{\geq0}$,
define the function $\cewt:p(V)\to\bbZ_{\geq0}$ by 
\begin{align*}
  \cewt([v,v']) = \sum_{\kap,\kapr\in\cbr{1,2}}\ewt([(v,\kap),(\vpr,\kapr)])
\end{align*}
for $[v,v']\in p(V)$.
We call this function $\cewt$ {\it the projected edge weight (function)}.
%
%%%%% Set of edges, projected graph %%%%%
For a weighted graph $G=(V,\ewt,\vwq)$,
we denote 
\begin{align*}
  E(G) = \cbr{[v,\vpr]\in p(V)\mid \cewt([v,\vpr])>0},
\end{align*}
which consists of the pairs of $V$ which has positive projected edge weight,
and we call 
$(V(G),E(G))(=(V,E(G)))$
{\it the projected graph} of the weighted graph $G$.
\begin{comment}
  {\myred projected graphに$\check{G}$という記号を与えては？}
\end{comment}
% 
%%%%% connectedness %%%%%
If the projected graph $(V(G), E(G))$ of a weighted graph $G$ 
is a connected graph in a usual sense,
we say $G$ is {\it connected}.
%
%%%%% component %%%%%
Writing $\check{\Comp}(G)$ for
the set of connected components of the projected graph $(V(G),E(G))$
of $G$,
we define $\Comp(G)$ by 
\begin{align*}
  \Comp(G) = \cbr{(V',\ewt|_{\hp(V')}, \vwq|_{\whv'})\mid 
  %\Comp(G) = \cbr{(V',\ewt\bbone_{\hp(V')}, \vwq\bbone_{\whv'})\mid 
  (V',E')\in\check\Comp(G)},
\end{align*}
and call an element of $\Comp(G)$ a {\it component} of $G$.

%%%%%%%%%%%%%%%%%%%%%%%%%%%%%%%%%%%%%%%%%%%%%
%%%%% summarizing quantity of graph %%%%%%%%%
%%%%%%%%%%%%%%%%%%%%%%%%%%%%%%%%%%%%%%%%%%%%%
% \paragraph{{\rb Summarizing quantity}}
To define the exponent for weighted graphs, we introduce some summarizing quantities.
We start with those related to weights on vertices. 
Let $G=(V,\ewt,\vwq)$ be a weighted graph.
%
%%%%% Total weight on vertices $\barq(G)$, $\barq_1(G)$, $\barq_2(G)$ %%%%%
We denote
\begin{align*}
  \barq(G)=\sum_{\hv\in\whv(G)}\vwq(\hv),\tand
    % =\sum_{v\in V(G)}\sum_{\kap=1,2}\vwq(v,\kap)
  \barq_\kap(G)=\sum_{v\in V(G)}\vwq(v,\kap)
  \tfor \kap=1,2.
\end{align*}
Since 
$\barq(G)=\sum_{v\in V(G)}\sum_{\kap=1,2}\vwq(v,\kap)$,
it holds that 
$\barq(G)=\barq_1(G)+\barq_2(G)$.
%
%%%%% Classification of components %%%%%
% $\Comp_0(G)$, $\Comp_+(G)$, $\Comp_1(G)$, $\Comp_2(G)$;
With the above quantities, we classify the components of $G$ as follows:
\begin{itemize}
  \item $\Comp_0(G)=\cbr{C\in\Comp(G)\mid \barq(C)=0},$
  \item $\Comp_+(G)=\cbr{C\in\Comp(G)\mid \barq(C)>0},$
  \begin{itemize}
    \item $\Comp_1(G)=\bcbr{C\in\Comp(G)\mid \barq_1(C)>0},$
    \item $\Comp_2(G)=\bcbr{C\in\Comp(G)\mid \barq_1(C)=0 \tandsm \barq_2(C)>0}.$
  \end{itemize}
\end{itemize}
% 
% \begin{align*}
%   \Comp_0(G)&=\cbr{C\in\Comp(G)\mid \barq(C)=0},&
%   \Comp_+(G)&=\cbr{C\in\Comp(G)\mid \barq(C)>0},
%   \\*
%   \Comp_1(G)&=\bcbr{C\in\Comp(G)\mid \barq_1(C)>0},&
%   \Comp_2(G)&=\bcbr{C\in\Comp(G)\mid \barq_1(C)=0 \tandsm \barq_2(C)>0}.
% \end{align*}
Notice that 
$\Comp_+(G) = \Comp_1(G)\sqcup\Comp_2(G)$ and 
$\Comp(G) = \Comp_+(G)\sqcup\Comp_0(G)$.
%
%%%%% Total weight on edges $\bartheta(G)$ %%%%%
As for edges, we define %for a weighted graph $G=(V,\ewt,\vwq)$
\begin{align*}
  \bartheta(G)=\sum_{[\hv,\hv']\in\hp(V(G))}\ewt([\hv,\hv']).
\end{align*}
To obtain a sharper estimate, we classify the edges of 
the projected graph of a weighted graph.
For a weighted graph $G=(V,\ewt,\vwq)$,
we denote 
\begin{align*}
  % E_1(V',\ewt) &= \cbr{[v,\vpr]\in p(V')\mid \cewt_1([v,\vpr])>0 \tandsm 
  % \cewt_2([v,\vpr])=0}
  % \\
  E_2(G) &= \cbr{[v,\vpr]\in p(V)\mid \cewt_2([v,\vpr])>0},
\end{align*}
where we define 
$\cewt_2:p(V)\to\bbZ_{\geq0}$ by 
\begin{align*}
  \cewt_2([v,\vpr]) = 
  \ewt([(v,2),(\vpr,2)])+
  \ewt([(v,1),(\vpr,2)])+
  \ewt([(v,2),(\vpr,1)]).
  % \sum_{\substack{\kap,\kapr\in\cbr{1,2}\\(\kap,\kapr)\neq(1,1)}}
  %   \ewt([(v,\kap),(\vpr,\kapr)]).
\end{align*}
For a connected weighted graph $C$,
we say that $\bbT\subset E(C)$ is a {\it spanning tree} of $C$,
if the subgraph $(V(C), \bbT)$ of the $(V(C), E(C))$ 
is a tree in a usual sense, that is a connected graph without cycles.
When $\abs{V(C)}=1$ and hence $p(V(C))=\emptyset$, 
we regard $\emptyset$ as the only spanning tree of $C$.
%
%%%%% exponent for connected weighted graphs, $\ell_2(C)$ %%%%%
Then 
the exponent $e(C)$ of a connected weighted graph $C$ is defined by
\begin{align*}%\label{def:231005.1440}
  e(C) = \et(C) + \eq(C)
\end{align*}
with
\begin{align*}
  \eq(C) &=
  \begin{cases}
    -\half -H\barq(C)
    % \rbr{-H-\half} -H(\barq(C)-1) {\myred(=-\half -H\barq(C))}
    &\tifsm \barq_2(C)>0 \tandsm \barq_1(C)=0 
    \\
    -1 -H(\barq(C)-1)               
    &\tifsm \barq_1(C)>0 
    \\
    0 
    & \tifsm \barq(C)=0
  \end{cases}
  \\
  \et(C) &=
  % 1 - (\abs{V(C)}-1) + (1-2H) \ell_2(C) 
  % - 2H (\bartheta(C)-(\abs{V(C)}-1))
  % \\&=
  1 - 2H\bartheta(C) + (2H-1) \brbr{\abs{V(C)}-1-\ell_2(C)},
\end{align*}
and 
\begin{align*}
  \ell_2(C) = \max_\bbT\abs{E_2(C)\cap \bbT}
\end{align*}
where the maximum is taken over the set of all spanning trees of $C$. 
Notice that 
when $\abs{V(C)}=1$,
% for a connected weighted graph $C$ such that $\abs{V(C)}=1$,
it holds that $e_\theta(C)=1$, since $\bartheta(C)=0$ and $\ell_2(C)=0$.
%
%%%%% Exponent for a weighted graph %%%%%
The exponent of a weighted graph $G$ is defined as 
the sum of the exponent of its components, namely 
\begin{align*}
  e(G)=\sum_{C\in\Comp(G)} e(C).
\end{align*}

% The estimates of the next lemma is basic in the following arguments.
% \begin{lemma}
% (i) For any connected graph $C$,
%   \begin{align}
%     \eq(C)\leq-H\barq(C). \label{221206.1156}  
%   \end{align}
% (ii) For any connected graph $C$ such that $\barq(C)>0$,
%   \begin{align}
%     \eq(C)&\geq-\half-H\barq(C), \label{221206.1202}
%     \\
%     \eq(C)&\leq-1-H(\barq(C)-1). \label{eq:230818.1803}
%   \end{align}
% \end{lemma}

%%%%%%%%%%%%%%%%%%%%%%%%%%%%%%
%%%%% factor B %%%%%%%%%%%%%%%
%%%%%%%%%%%%%%%%%%%%%%%%%%%%%%
% \paragraph{{\rb Factor $B$}}
% \subsubsection*{Realization of functional}
Now we introduce 
% the functionals represented by aforementioned weighted graphs 
the correspondence between the weighted graphs and functionals.
%
%%%%% \beta %%%%%
Suppose that 
$V$, $V_\bff$ and $V_\ewt$ are nonempty finite sets
satisfying $V\subset V_\bff, V_\ewt$.
Let
$\bff=(\bff^{(\hv)})_{\hv\in\wtv_\bff}\in\calh^{\wtv_\bff}$.
% $\bff\in\abs{\calh}^{\wtv_\bff}$
For 
$\ewt:\hp(V_\ewt)\to\bbZ_{\geq0}$,
% such that $\supp(\ewt)\subset\hp(V)$,
we denote
\begin{align*}
  \beta_{V}(\bff, \ewt) = 
  \prod_{[\tv, \tv']\in\hp(V)} \babr{\bff^{(\tv)}, \bff^{(\tv')}}^{\ewt[\tv,\tv']}.
\end{align*}
If $\abs{V}=1$, and therefore $\hp(V)=\emptyset$, we define 
$\beta_{V}(\bff, \ewt)=1$.
% Notice that $\beta_{V}(\bff, \ewt)=1$
% when $\sum_{[\tv,\tv']\in\tp(V)}\ewt[\tv,\tv']=0$.
%
%%%%% \delta %%%%%
For a nonempty finite set $V_\vwq$ satisfying $V\subset V_\vwq$ and 
$\vwq:\whv_\vwq\to\bbZ_{\geq0}$, % such that $\supp(\vwq)\subset\whv$,
we denote 
\begin{align*}
  \delta_V(\bff, \vwq)
  = 
  \delta^{\bar\vwq}
  \Brbr{
    \subotimes{\hv\in\whv}
    \brbr{\nrbr{\bff^{(\hv)}}^{\otimes \vwq(\hv)}}},
\end{align*}
which is a $\bar\vwq$-th Skorohod integral with
$\bar\vwq = \sum_{\hv\in\whv}\vwq(\hv)$.
When $\bar\vwq=0$ (i.e. $\vwq=0$), 
we define $\delta_V(\bff, \vwq)=1$.
%
%%%%% B %%%%%
For a connected weighted graph $C=(V, \ewt, \vwq)$,
% $\bff\in{\calh}^{\wtv_\bff}$ with some $V_\bff(\supset V)$,
we define 
\begin{align*}
  \beta(C, \bff) = \beta_{V}(\bff, \ewt),\qquad
  \delta(C,\bff) = \delta_{V}(\bff, \vwq)\tand
  B(C,\bff) = \beta(C, \bff)\;\delta(C, \bff).
\end{align*}

%%%%%%%%%%%%%%%%%%%%%%%%%%%%%%
%%%%% conditions on bbf %%%%%%
%%%%%%%%%%%%%%%%%%%%%%%%%%%%%%
% \paragraph{{\rb Conditions on $\bbf$}}
%%%%% functionals %%%%%
%
% %%%%% $i:V\to\cbr{1,2},j\in\bbJ_n(V,i)$;  %%%%%
% For nonempty finite sets $V$ and $V'$ with $V\subset V'$
% and $i:V'\to\cbr{1,2}$,
% we write 
% $\bbJ_n(V,i)=\prod_{v\in V} [i(v)n-1]$.
% %
% 
%
% \item $\calf(V), \bbf, \cala(V,i), \cali_n^{(i)}$
%%%%% \calf(V) %%%%%
For a nonempty finite set $V$,
consider the following set of conditions for 
$\bbf=(\bbf^{(\hv)})_{\hv\in\whv}\in{L^\infty(\bbR)}^{\whv}$:
% $\bbf\in\rbr{L^\infty(\bbR)}^{\whv}$:
\begin{itemize}
  \item For $v\in V$, 
  the support of the function $\bbf^{(v,1)}$ 
  % the support of the function $\bbf^{(v,1)}\in L^\infty(\bbR)$ 
  is included in $[a,a+1]$ with some $a\in[-1,0]$. 

  \item For $v\in V$, 
  the support of the function $\bbf^{(v,2)}$ 
  % the support of the function $\bbf^{(v,2)}\in L^\infty(\bbR)$ 
  is included in $[-1,0]$. 
\end{itemize}
Denote by $\calf(V)$ the subset of 
${L^\infty(\bbR)}^{\whv}$ satisfying the above conditions.
% 
%%%%% $T_{m,j}, T^{(1)}_{m,j}, T^{(2)}_{m,j}$: 関数のtranslation %%%%%
For $f:\bbR\to\bbR$, $m\in\bbN$ and $j\in\bbZ$, 
we define the function $T_{m,j}(f):\bbR\to\bbR$ by 
\begin{align*}
  T_{m,j}(f)(x) = f(mx-j),
\end{align*}
and we denote 
\begin{align*}
  T_{m,j}^{(1)}(f) = T_{m,j}(f)\tand
  T_{m,j}^{(2)}(f) = T_{m,j+1}(f) - T_{m,j}(f).
\end{align*}
% 
%
%%%%% \bbf %%%%%
For $\bbf\in\calf(V)$, 
% $(v,\kappa)\in\whv$,
$n\in\ntwo=\bbZ_{\geq2}$ and
% $n\in\ntwo=\cbr{j\in\bbZ\mid j\geq2}$ and
$j\in[n-1]^V$,
we define $\bbf_{n,j}=\brbr{\bbf^{(v,\kap)}_{n,j}}_{(v,\kap)\in\wtv}
\in{L^\infty(\bbR)}^{\wtv}$
% \in\abs{\calh}^{\wtv}$
by
\begin{align*}
  \bbf^{(v,\kappa)}_{n,j}  %= \bbf^{(v,\kappa)}_{n,i,j}
  &= T_{n,j_v}^{(\kappa)}(\bbf^{(v,\kappa)}),
\end{align*}
for each $(v,\kappa)\in\whv$.
Notice that $\bbf_{n,j}\in\abs{\calh}^{\wtv}$.

%%%%%%%%%%%%%%%%%%%%%%%%%%%%%%%%%%%%%%%%%%%%%%%%%%%%%%%
%%%%% conditions on weight functional A %%%%%%%%%%%%%%%
%%%%%%%%%%%%%%%%%%%%%%%%%%%%%%%%%%%%%%%%%%%%%%%%%%%%%%%
% \paragraph{{\rb Conditions on weight functional $A$}}
For a nonempty finite set $V$, %. {\redmy and $i:V\to\cbr{1,2}$.}
we consider a family of functionals
$A=(A_{n,j})_{n\in\bbN_2, j\in[n-1]^V}$ with
$A_{n,j}\in\bbD^{\infty}$
satisfying the following conditions:
% Writing $\bbJ(V,i) = \bigsqcup_{n\geq2} \cbr{n}\times\bbJ_n(V,i)$,
% we consider the following conditions
% for $A=(A_{n,j})_{(n,j)\in\bbJ(V,i)}\in(\bbD^\infty)^{\bbJ(V,i)}$:
\begin{itemize}
  \item For any $p\geq1$,
  $\sup_{n\geq2, j\in[n-1]^V}\norm{A_{n,j}}_{L^p(P)}<\infty$.
  \item For $k\geq1$, %$j\in\bbJ_n(V,i)$, 
  $D^kA_{n,j}$ is $\abs{\calh}^{\otimes k}$-valued,
  % $D^kA_{n,j}$ is an $\abs{\calh}^{\otimes k}$-valued random variable.
  % \item $\sup_{n\geq2, j\in\bbJ_n(V,i)}\norm{A_{n,j}}_{L^p(P)}<\infty$ for any $p\geq1$.
  % \item For $k\in\bbN$, $D^kA_{n,j}$ is 
  and represented by 
  $(D^k_{s_1,...,s_k}A_{n,j})_{s_1,...,s_k\in[0,1]}$ 
  such that
  \begin{align*}
    \sup_{n\geq2, j\in[n-1]^V}
    % \sup_{n\geq2, j\in\bbJ_n(V,i)}
    \sup_{s_1,...,s_k\in[0,1]}\norm{D^k_{s_1,...,s_k}A_{n,j}}_{L^p(P)}
    <\infty
    \tfor p\geq1.
  \end{align*}
\end{itemize}
We denote by $\cala(V)$ %$\cala(V,i)$
the set of the families of functionals
% the subset of $(\bbD^\infty)^{\bbJ(V,i)}$ 
satisfying the above conditions.

%%%%%%%%%%%%%%%%%%%%%%%%%%%%%%
%%%%% def of functional %%%%%%
%%%%%%%%%%%%%%%%%%%%%%%%%%%%%%
% \paragraph{{\rb Def. of functional}}
Having prepared the definitions above, we define the functional to consider in the following arguments.
For $n\in\ntwo$, a weighted graph $G=(V,\ewt, \vwq)$, %$i:V\to\cbr{1,2}$, 
$A\in\cala(V)$ 
% $A\in\cala(V,i)$ 
and $\bbf\in\calf(V)$,
we define the functional $\cali_n\rbr{G, A,\bbf}$ by
% we define the functional $\calii_n\rbr{G, A,\bbf}$ by
\begin{align}
  \cali_n\rbr{G, A,\bbf} &=
  % \calii_n\rbr{G, A,\bbf} &=
  \sum_{j\in[n-1]^V}
  % \sum_{j\in\bbJ_n(V,i)}
  A_{n,j} %A_n(j)
  \prod_{C\in\Comp(G)} B\rbr{C, \bbf_{n,j}}.
  \label{eq:240622.2118}
\end{align}

Then we have the following result from
Propositions 4.4 and 4.5 in \cite{2024Yamagishi-asymptotic}.
\begin{theorem}\label{thm:240628.1910}
  Let $G=(V,\ewt,\vwq)$ be a weighted graph,
  $A\in\cala(V)$ and 
  % $A=(A_{n,j})_{n\in\bbN_2, j\in[n-1]^V}\in\cala(V)$ and
  $\bbf\in\calf(V)$.
  Consider the functional 
  $\cali_n=\cali_n(G,A,\bbf)$
  defined by \eqref{eq:240622.2118}.
  Then the Sobolev norm of the functional $\cali_n$ is bounded as 
  \begin{align*}
    \norm{\cali_n}_{k,p}=O(n^{e(G)})
  \end{align*}
  for any $k\geq1$ and $p\geq1$.
  In particular, it holds that 
  $\norm{\cali_n}_{L^p}=O(n^{e(G)})$.

\end{theorem}

% \subfile{expo-Hurst-fSDE/sec_exponent_def.tex}
% \subfile{expo-Hurst-fSDE/sec_exponent_main.tex}

\subsubsection{The change of exponent by the operator $D_{v_n}$}
\label{sec:240626.1835}
For $q\geq2$,
consider the following $\calh$-valued random variable $v_n^{(q)}$
\begin{align}
  v_n^{(q)}&=
  % u_n^{(q)}
  n^{q H-\half}
  \sum_{{j_0\in[n-1]}}
  A'_{n,j_0}
  I_{q-1}({\diffker^n_{j_0}}^{\otimes q-1})
  \diffker^n_{j_0}.
  \label{eq:240617.1534}
\end{align}
Here $A'=(A'_{n,j_0})_{n\in\bbN_2,j_0\in[n-1]}$ satisfies 
$A'\in\cala(V')$ with some singleton $V'$.
Typical examples of $v_n^{(q)}$ will appear in Section \ref{sec:240617.1526}
as $u_n^{(\ell)}$ {with $q=2\ell$}.
For a functional $\cali_n$ of the form \eqref{eq:240622.2118},
the following proposition shows that
$D_{v_n^{(q)}}\cali_n$ again can be written as %decomposes into 
a sum of functionals of the form \eqref{eq:240622.2118}
and describes the behavior of the maximum of the asymptotic order of them
in terms of exponents.

In the following argument, given $V$,
we take $V'=\cbr{v_0}$ such that $V\cap V'=\emptyset$ and 
write $\extVertex=V\sqcup V'$.
For $\bbf\in\calf(V)$,
we define an element $\extKer$ of 
$\calf(\extVertex)$ by
% ${L^\infty(\bbR)}^{\widehat\extVertex}$ by
\begin{align*}
  \extKer|_{\whv}=\bbf,\tand
  \extKer^{(v_0,2)}=\bbone_{[-1,0]}.
\end{align*}
Since $\extKer^{(v_0,1)}$ is not used in the following argument, 
$\extKer^{(v_0,1)}$ can be set arbitrarily, say $\bbone_{[0,1]}$.

\begin{proposition}\label{prop:240617.2009}
  Let $G=(V,\ewt,\vwq)$ be a weighted graph,
  % Suppose that 
  $A\in\cala(V)$ and $\bbf\in\calf(V)$.
  Let $q\geq2$, and 
  consider $v_n^{(q)}$ of \eqref{eq:240617.1534} with 
  some $A'\in\cala(V')$.
  % $V'=\cbr{v_0}$ be a singleton such that $v_0\not\in V$.

  Then the functional 
  $\babr{D\calii_n(G,A,\bbf), v_n^{(q)}}$
  % , which we shall denote by $D_{v_n^{(q)}}\cali_n$,
  decomposes as 
  \begin{align*}
    \babr{D\calii_n(G,A,\bbf), v_n^{(q)}}=
    \sum_{\gamma\in\Gamma} 
    n^{\alpha^{(\gamma)}}
    \calii_n(G^{(\gamma)},A^{(\gamma)},\extKer)
  \end{align*}
  with some finite set $\Gamma$, and 
  $\alpha^{(\gamma)}\in\bbR$, 
  weighted graphs $G^{(\gamma)}$ and
  $A^{(\gamma)}\in\cala(\extVertex)$
  for $\gamma\in\Gamma$
  satisfying 
  \begin{align*}
    \alpha^{(\gamma)} + e(G^{(\gamma)}) \leq e(G).
  \end{align*}
  % for all $\gamma\in\Gamma$.

  Furthermore, %{\rb the following assertions hold:}
  \item [(i)] if $\Comp_+(G)=\emptyset$, i.e. $\barq(G)=0$,
  then $\Gamma$ is a singleton (denoted by $\Gamma=\cbr{\gamma_0}$) and 
  $\alpha^{(\gamma_0)} + e(G^{(\gamma_0)}) = e(G)-H$.

  \item [(ii)] If $\barq(C)\neq q$ for all $C\in\Comp_+(G)$, then
  \begin{align*}
    \alpha^{(\gamma)} + e(G^{(\gamma)}) \leq e(G)-\half
  \end{align*}
  holds for any $\gamma\in\Gamma$.
\end{proposition}

\begin{proof}
The functional {$\babr{D\calii_n(G,A,\bbf), v_n^{(q)}}$} decomposes as
\begin{align*}
  &\babr{D\calii_n(G,A,\bbf), v_n^{(q)}}
  % &\abr{D\cali_n, v_n^{(q)}}
  % \\&=
  % \Big\langle
  %   D\Bcbr{\sum_{j\in[n-1]^V}
  %   A_{n,j} 
  %   \prod_{C\in\Comp(G)} B\rbr{C, \bbf_{n,j}}},
  % % \\&\qquad
  % n^{q H-\half}
  % \sum_{j_0\in[n-1]}
  % A'_{n,j_0}
  % I_{q-1}({\diffker^n_{j_0}}^{\otimes q-1})
  % \diffker^n_{j_0}
  % \Big\rangle
  % % \\&=
  % % n^{q H-\half}
  % % \sum_{j\in[n-1]^V}
  % % \sum_{j_0\in[n-1]}
  % % \abr{D\cbr{
  % %   A_{n,j} 
  % %   \prod_{C\in\Comp(G)} B\rbr{C, \bbf_{n,j}}},
  % % \diffker^n_{j_0}}
  % % A'_{n,j_0}
  % % I_{q-1}({\diffker^n_{j_0}}^{\otimes q-1})
  \\&=
  n^{q H-\half}
  \sum_{\substack{j\in[n-1]^V\\j_0\in[n-1]}}
  \abr{D\Bcbr{
    A_{n,j} 
    \prod_{C\in\Comp(G)} B\rbr{C, \bbf_{n,j}}},
  \diffker^n_{j_0}}
  A'_{n,j_0}
  I_{q-1}({\diffker^n_{j_0}}^{\otimes q-1})
  % \quad\text{\redmy どっちか消す．}
  \\&=
  \cali_n^{(A)}+
  \sum_{C_1\in\Comp_+(G)}
  \cali_n^{(C_1)},
\end{align*}
where we define 
\begin{align*}
  \cali_n^{(A)}&=
  n^{q H-\half}
  \sum_{\substack{j\in[n-1]^V\\j_0\in[n-1]}}
  \abr{DA_{n,j}, \diffker^n_{j_0}}
  % \abr{D\cbr{
  %   A_{n,j} 
  %   \prod_{C\in\Comp(G)} B\rbr{C, \bbf_{n,j}}},
  % \diffker^n_{j_0}}
  A'_{n,j_0}
  I_{q-1}({\diffker^n_{j_0}}^{\otimes q-1})
  \prod_{C\in\Comp(G)} B\rbr{C, \bbf_{n,j}}
  \\
  \cali_n^{(C_1)}&=
  n^{q H-\half}
  \sum_{\substack{j\in[n-1]^V\\j_0\in[n-1]}}
  A_{n,j}
  \abr{DB\rbr{C_1, \bbf_{n,j}},\diffker^n_{j_0}}
  A'_{n,j_0}
  I_{q-1}({\diffker^n_{j_0}}^{\otimes q-1})
  \prod_{\substack{C\in\Comp(G)\\C\neq C_1}} B\rbr{C, \bbf_{n,j}}
\end{align*}

Assume that $\Comp_+(G)$ is not empty for a while and 
fix $C_1\in\Comp_+(G)$.
Let us write $V_1:=V(C_1)$ and $\whv_1:=\whv(C_1)$.
We can write
\begin{align*}
  \abr{DB\rbr{C_1, \bbf_{n,j}},\diffker^n_{j_0}}
  &=
  % \abr{D\rbr{\beta\rbr{C_1, \bbf_{n,j}}
  % \delta\rbr{C_1, \bbf_{n,j}}},\diffker^n_{j_0}}
  % =
  \beta\rbr{C_1, \bbf_{n,j}}
  \abr{D\delta\rbr{C_1, \bbf_{n,j}},\diffker^n_{j_0}}
\end{align*}
and
\begin{align*}
  D\delta\rbr{C_1, \bbf_{n,j}}&=
  \sum_{\substack{\hv_1\in\whv_1\\\vwq(\hv_1)>0}}
  \vwq(\hv_1)\times
  \delta^{\bar\vwq_{C_1}-1}
  % \delta^{{\redmy \bar\vwq}-1}
  \Brbr{\subotimes{\hv\in\whv_1}
  \brbr{\nrbr{\bbf_{n,j}^{(\hv)}}^{\otimes \vwq(\hv)-\bbone_\cbr{\hv_1}(\hv)}}}
  \bbf_{n,j}^{(\hv_1)},
\end{align*}
% {\redmy[上か下かどちらかを消す]
% \begin{align*}
%   \abr{D\delta\rbr{C_1, \bbf_{n,j}},\diffker^n_{j_0}}&=
%   \sum_{\substack{\hv_1\in\whv_1\\\vwq(\hv_1)>0}}
%   \vwq(\hv_1)\times
%   \delta^{{\redmy \bar\vwq}-1}
%   \Brbr{\subotimes{\hv\in\whv_1}
%   \brbr{\nrbr{\bbf_{n,j}^{(\hv)}}^{\otimes \vwq(\hv)-\bbone_\cbr{\hv_1}(\hv)}}}
%   \abr{\bbf_{n,j}^{(\hv_1)},\diffker^n_{j_0}}
% \end{align*}}  
where we denote 
$\bar\vwq_{C_1}=\sum_{\hv\in\whv_1}\vwq(\hv)$.
Hence we have
\begin{align*}
  \cali_n^{(C_1)}&=
  \sum_{\substack{\hv_1\in\whv_1\\\vwq(\hv_1)>0}}
  \cali_n^{(C_1,\hv_1)}
\end{align*}
with 
\begin{align*}
  \cali_n^{(C_1,\hv_1)}
  &=% \\&=
  \vwq(\hv_1)\times
  n^{q H-\half}
  \sum_{\substack{j\in[n-1]^V\\j_0\in[n-1]}}
  A_{n,j} A'_{n,j_0}
  % \abr{DB\rbr{C_1, \bbf_{n,j}},\diffker^n_{j_0}}
  \times
  \prod_{\substack{C\in\Comp(G)\\C\neq C_1}} B\rbr{C, \bbf_{n,j}}
  \\&\hspace{50pt}\times
  \beta\rbr{C_1, \bbf_{n,j}}
  \abr{\bbf_{n,j}^{(\hv_1)},\diffker^n_{j_0}}
  \delta^{\bar\vwq_{C_1}-1}
  % \delta^{{\redmy \bar\vwq}-1}
  \Brbr{\subotimes{\hv\in\whv_1}
  \brbr{\nrbr{\bbf_{n,j}^{(\hv)}}^{\otimes \vwq(\hv)-\bbone_\cbr{\hv_1}(\hv)}}}
  I_{q-1}({\diffker^n_{j_0}}^{\otimes q-1})
  % \\&\hspace{50pt}
\end{align*}

By the product formula, we have
\begin{align*}
  &\delta^{\bar\vwq_{C_1}-1}
  \Brbr{\subotimes{\hv\in\whv_1}
  \brbr{\nrbr{\bbf_{n,j}^{(\hv)}}^{\otimes \vwq(\hv)-\bbone_\cbr{\hv_1}(\hv)}}}
  I_{q-1}({\diffker^n_{j_0}}^{\otimes q-1})
  \\&=
  \sum_{\pi\in\prodFormulaIndex}
  c(\pi)
  \delta^{\bar\vwq_{C_1}-2 -2\bar\pi}
  \Brbr{\subotimes{\hv\in\whv_1}
  \brbr{\nrbr{\bbf_{n,j}^{(\hv)}}
  ^{\otimes \vwq(\hv)-\bbone_\cbr{\hv_1}(\hv)-\pi(\hv)}}
  \otimes {\diffker^n_{j_0}}^{\otimes q-1-\sum_{\hv\in\whv_1}\pi(\hv)}}
  \\&\hspace{30pt}\times% \times% 
  \prod_{\hv\in\whv_1}\abr{\bbf_{n,j}^{(\hv)}, \diffker^n_{j_0}}^{\pi(\hv)}
\end{align*}
with 
\begin{align*}
  \prodFormulaIndex=
  \Bcbr{\pi:\whv_1\to\bbZ_{\geq0}:~
  \pi(\hv)\leq \vwq(\hv)-\bbone_\cbr{\hv_1}(\hv)
  \tforsm \hv\in\whv_1,
  \sum_{\hv\in\whv_1}\pi(\hv)\leq q-1}
\end{align*}
and 
$\bar\pi:=\sum_{\hv\in\whv_1}\pi(\hv)$.
Hence we have 
\begin{align*}
  &\cali_n^{(C_1,\hv_1)}=
  \sum_{\pi\in\prodFormulaIndex}
  \cali_n^{(C_1,\hv_1,\pi)}
\end{align*}
with
\begin{align*}
  \cali_n^{(C_1,\hv_1,\pi)}
  &=% \\*&=
  \vwq(\hv_1)\times c(\pi)\times
  n^{q H-\half}
  \sum_{\substack{j\in[n-1]^V\\j_0\in[n-1]}}
  A_{n,j} A'_{n,j_0}
  \times% \\&\hspace{50pt}\times
  \prod_{\substack{C\in\Comp(G)\\C\neq C_1}} B\rbr{C, \bbf_{n,j}}
  \\&\hspace{50pt}\times
  \beta\rbr{C_1, \bbf_{n,j}}
  \prod_{\hv\in\whv_1}\abr{\bbf_{n,j}^{(\hv)}, \diffker^n_{j_0}}
  ^{\pi(\hv)+\bbone_\cbr{\hv_1}(\hv)}
  \\&\hspace{50pt}\times
  \delta^{\bar\vwq_{C_1}-2 -2\bar\pi}
  \Brbr{\subotimes{\hv\in\whv_1}
  \brbr{\nrbr{\bbf_{n,j}^{(\hv)}}
  ^{\otimes \vwq(\hv)-\bbone_\cbr{\hv_1}(\hv)-\pi(\hv)}}
  \otimes {\diffker^n_{j_0}}^{\otimes q-1-\sum_{\hv\in\whv_1}\pi(\hv)}}.
\end{align*}
Thus the functional $\babr{D\cali_n, u_n^{(q)}}$ is decomposed as 
\begin{align}
  \babr{D\cali_n, u_n^{(q)}}&=
  \cali_n^{(A)}+
  \sum_{C_1\in\Comp_+(G)}
  \sum_{\substack{\hv_1\in\whv_1\\\vwq(\hv_1)>0}}
  \sum_{\pi\in\prodFormulaIndex}
  \cali_n^{(C_1,\hv_1,\pi)}.
  \label{eq:240614.1408}
\end{align}

We will write the above functionals in terms of weighted graphs.
First we consider $\cali_n^{(C_1,\hv_1,\pi)}$.
Define the weighted graph 
$\newCompProd = 
(\widetilde V_1, \ewt_{\hv_1,\pi}, \vwq_{\hv_1,\pi})$ by
\begin{itemize} 
  \item $\widetilde V_1 = V_1\sqcup \cbr{v_0}$,
  % \item $\ewt_{\hv_1,\pi} = 
  % \ewt\bbone_{\tp(V_1)} + 
  % \sum_{\hv\in\whv_1} \bbone_{[\hv,(v_0,2)]}\times
  % (\pi(\hv)+\bbone_\cbr{\hv_1}(\hv))$.

  \item $\ewt_{\hv_1,\pi} = 
  \ewt\bbone_{\tp(V_1)} 
  +\bbone_{[\hv_1,(v_0,2)]}
  +\sum_{\hv\in\whv_1} \bbone_{[\hv,(v_0,2)]} \times\pi(\hv)$,

  \item $\vwq_{\hv_1,\pi}=
  (\vwq-\pi)\bbone_{\wtv_1}-\bbone_\cbr{\hv_1}
  +\bbone_\cbr{(v_0,2)}\times(q-1-\sum_{\hv\in\whv_1}\pi(\hv))$.
\end{itemize}
Also define
$\newGraph^{(C_1,\hv_1,\pi)}$ and  $\newWeight^{(C_1,\hv_1,\pi)}$ by
\begin{align*}
  \newGraph^{(C_1,\hv_1,\pi)}&=
  \newCompProd\vee\Brbr{
  \mathop\vee_{\substack{C\in\Comp(G), C\neq C_1}} C},
  \\% \tand% \\
  \newWeight_{n,\extIndex}^{(C_1,\hv_1,\pi)}&=
  \vwq(\hv_1)\times c(\pi)\times
  A_{n,j} A'_{n,j_0}
  \tfor \extIndex\in[n-1]^{\extVertex}.
\end{align*}
% {\rb 英語的にチェックせよ．$\newWeight^{(C_1,\hv_1,\pi)}$を雑に定義しているのは構わない．}
Setting 
$\alpha^{(C_1,\hv_1,\pi)}=q H-\half$,
we can write 
\begin{align}
  \cali_n^{(C_1,\hv_1,\pi)}
  &=
  n^{\alpha^{(C_1,\hv_1,\pi)}}\times
  % n^{q H-\half}\times
  \cali_n\rbr{\newGraph^{(C_1,\hv_1,\pi)}, \newWeight^{(C_1,\hv_1,\pi)},\extKer}.
  \label{eq:240617.1821}
\end{align}

As for $\cali_n^{(A)}$,
we define $\newGraph^{(A)}$ by 
$\newGraph^{(A)} = \newComp_0\vee G$ with 
\begin{align*}
  \newComp_0=(\cbr{v_0},\emptyset,
  ((v_0,2)\mapsto q-1, (v_0,1)\mapsto 0))
\end{align*}
and the weight functional $\newWeight^{(A)}$ by
\begin{align}
  \newWeight_{n,\extIndex}^{(A)} = 
  n^{2H}\abr{DA_{n,j}, \diffker^n_{j_0}} A'_{n,j_0}
  \label{eq:240617.1822}
\end{align}
for $\extIndex\in[n-1]^{\extVertex}$.
Setting 
$\alpha^{(A)}=(q -2)H-\half$,
we have 
\begin{align*}
  \cali_n^{(A)}&=
  n^{\alpha^{(A)}}\times
  % n^{(q -2)H-\half}\times
  \cali_n\rbr{\newGraph^{(A)}, \newWeight^{(A)},\extKer}.
\end{align*}

Next we need to calculate the exponent of the weighted graphs 
$\newGraph^{(C_1,\hv_1,\pi)}$ and $\newGraph^{(A)}$.
First we deal with $\newGraph^{(C_1,\hv_1,\pi)}$.
The following elementary relations for $\newCompProd$ hold:
\begin{itemize}
  \item $\abs{V(\newCompProd)}=\abs{V_1}+1$.

  \item $\bartheta(\newCompProd)=
  \bartheta(C_1)+\sum_{\hv\in\whv_1}\pi(\hv)+1$.

  \item $\barq(\newCompProd)=
    \barq(C_1)
    - 2\sum_{\hv\in\whv_1}\pi(\hv) -2 + q$.    
\end{itemize}
We can also show that 
$\ell_2(\newCompProd)\geq\ell_2(C_1)+1$ since 
$[v_1,v_0]\in E_2(\newCompProd)$,
where $v_1$ is specified by 
$\hv_1=(v_1,\kap_1)$.
Hence we have 
\begin{align*}
  \et(\newCompProd) &=
  1 - 2H\bartheta(\newCompProd) 
  + (2H-1) \brbr{\babs{V(\newCompProd)}-1-\ell_2(\newCompProd)},
  % \\&=
  % 1-2H(\bartheta(C_1)+\sum_{\hv\in\whv_1}\pi(\hv)+1)
  % + (2H-1) \brbr{(\abs{V_1}+1)-1-\ell_2(\newCompProd)}
  \\&\leq
  1-2H(\bartheta(C_1)+\sum_{\hv\in\whv_1}\pi(\hv)+1)
  + (2H-1) \brbr{(\abs{V_1}+1)-1-(\ell_2(C_1)+1)}
  % \\&=
  % 1-2H\bartheta(C_1)
  % -2H(\sum_{\hv\in\whv_1}\pi(\hv)+1)
  % + (2H-1) \brbr{\abs{V_1}-1-\ell_2(C_1)}
  \\&=
  \et(C_1)
  -2H(\sum_{\hv\in\whv_1}\pi(\hv)+1)
  \\&=
  \et(C_1)
  -H(\barq(C_1)+ q-\barq(\newCompProd))
\end{align*}

\paragraph{(1-1)}
Consider the case where $C_1\in\Comp_2(G)$ and $\barq(\newCompProd)>0$.
Obviously, we have $\barq_1(\newCompProd)=0$, and hence 
$\eq(\newCompProd)=-\half-H\barq(\newCompProd)$.
The exponent $e(\newCompProd)$ is estimated as 
\begin{align*}
  e(\newCompProd)&\leq
  \et(C_1)
  -H(\barq(C_1)+ q-\barq(\newCompProd))
  -\half-H\barq(\newCompProd)
  % \\&=
  % \et(C_1)
  % -H(\barq(C_1)+ q)
  % -\half
  % =
  % \et(C_1)+\eq(C_1) -H q
  =e(C_1) -qH,
\end{align*}
and hence
\begin{align}
  \alpha^{(C_1,\hv_1,\pi)} + e(\newGraph^{(C_1,\hv_1,\pi)})=
  (q H-\half) + e(\newCompProd)
  +\sum_{\substack{C\in\Comp(G), C\neq C_1}} e(C)\leq
  e(G)-\half.
  \label{eq:240617.1831}
\end{align}

\paragraph{(1-2)}
The case where $C_1\in\Comp_2(G)$ and $\barq(\newCompProd)=0$.
Since $\eq(\newCompProd)=0$ by definition, 
we have 
\begin{align*}
  e(\newCompProd)&=\et(\newCompProd)
  \leq% \\&\leq
  % \et(C_1)-H(\barq(C_1)+ q-\barq(\newCompProd))
  % \\&=
  \et(C_1)-H(\barq(C_1)+ q)
  % =
  % \et(C_1)+\eq(C_1)+\half - Hq 
  =
  e(C_1)+\half - Hq %\eq(C)
\end{align*}
and 
\begin{align}
  \alpha^{(C_1,\hv_1,\pi)} + e(\newGraph^{(C_1,\hv_1,\pi)})
  % =
  % (q H-\half) + e(\newCompProd)
  % +\sum_{\substack{C\in\Comp(G), C\neq C_1}} e(C)
  \leq 
  e(G).
  \label{eq:240617.1832}
\end{align}

\paragraph{(1-3)}
The case where $C_1\in\Comp_1(G)$ and $\barq(\newCompProd)>0$.
% By \eqref{eq:230818.1803}
From the definition of $\eq$, we have 
$\eq(\newCompProd)\leq-1-H(\barq(\newCompProd)-1)$,
and
\begin{align*}
  e(\newCompProd)&\leq
  \et(C_1)
  -H(\barq(C_1)+ q-\barq(\newCompProd))
  -1-H(\barq(\newCompProd)-1)
  % =
  % \et(C_1)+\eq(C_1)-Hq
  =
  e(C_1)-Hq
\end{align*}
since 
$\eq(C_1) = -1-H(\barq(C_1)-1)$.
Hence we obtain
\begin{align}
  \alpha^{(C_1,\hv_1,\pi)} + e(\newGraph^{(C_1,\hv_1,\pi)})
  % =
  % (q H-\half) + e(\newCompProd)
  % +\sum_{\substack{C\in\Comp(G), C\neq C_1}} e(C)
  \leq 
  e(G)-\half.
  \label{eq:240617.1833}
\end{align}

\paragraph*{(1-4)}
The case where $C_1\in\Comp_1(G)$ and $\barq(\newCompProd)=0$.
Again by definition we have $\eq(\newCompProd)=0$.
Hence it holds that
\begin{align*}
  e(\newCompProd)&\leq
  \et(C_1) -H(\barq(C_1)+ q)
  % =
  % \et(C_1)+\eq(C_1)+1-H -Hq
  =
  e(C_1)+1-H -Hq
\end{align*}
and 
\begin{align}
  \alpha^{(C_1,\hv_1,\pi)} + e(\newGraph^{(C_1,\hv_1,\pi)})
  % =
  % (q H-\half) + e(\newCompProd)
  % +\sum_{\substack{C\in\Comp(G), C\neq C_1}} e(C)
  \leq 
  e(G)+\half-H .
  \label{eq:240617.1834}
\end{align}

\paragraph*{(2)}
Next we consider the exponent of $\newGraph^{(A)}$.
By the definition of $\newComp_0$ and the assumption $q\geq2$,
we have 
$\et(\newComp_0)=1$ and $\eq(\newComp_0)=-\half-H(q-1)$.
% $\barq_2(\newComp_0)>0$ and $\barq_1(\newComp_0)=0$, 
Hence it holds that
\begin{align}
  \alpha^{(A)} + e(\newGraph^{(A)})=
  ((q -2)H-\half) + e(\newComp_0) +e(G)=
  e(G)-H.
  \label{eq:240617.1835}
\end{align}

\vspsm
Recalling the decomposition \eqref{eq:240614.1408}
of $\babr{D\cali_n, u_n^{(q)}}$, 
the representations \eqref{eq:240617.1821} and \eqref{eq:240617.1822}
of functionals $\cali_n^{(C_1,\hv_1,\pi)}$ and $\cali_n^{(A)}$, 
and the estimates 
\eqref{eq:240617.1831},
\eqref{eq:240617.1832},
\eqref{eq:240617.1833},
\eqref{eq:240617.1834} and
\eqref{eq:240617.1835},
we obtain the first assertion of the proposition.

For the case (i) $\Comp_+(G)=\emptyset$,
the decomposition \eqref{eq:240614.1408} reduces to 
\begin{align*}
  \babr{D\cali_n, u_n^{(q)}}&=
  \cali_n^{(A)}
\end{align*}
and the estimate of {\bf(2)} gives the proof.

Consider the case (ii) $\barq(C)\neq q$ for all $C\in\Comp_+(G)$.
By the definition of $\prodFormulaIndex$,
it holds that 
\begin{align*}
  \sum_{\hv\in\whv_1}\pi(\hv)&\leq 
  (\barq(C_1)\wedge q)-1.
\end{align*}
in general, which implies 
\begin{align*}
  \barq(\newCompProd)&=\barq(C_1)- 2\sum_{\hv\in\whv_1}\pi(\hv) -2 + q
  \geq% \nn\\&\geq
  % \barq(C_1) -2\rbr{(\barq(C_1)\wedge q)-1} -2 + q=
  \barq(C_1)+ q -2(\barq(C_1)\wedge q).
  %\label{eq:240614.1654}
\end{align*}
Thus if $\barq(C)\neq q$ holds for all $C\in\Comp_+(G)$, then
we have $\barq(\newCompProd)>0$ for any $(C_1,\hv_1,\pi)$ 
in \eqref{eq:240614.1408}.
This means that only the cases {\bf(1-1), (1-3)} and {\bf(2)}
(i.e. \eqref{eq:240617.1831}, 
\eqref{eq:240617.1833} and 
\eqref{eq:240617.1835})
appear in the decomposition \eqref{eq:240614.1408}.
\end{proof}

% \newpage

\section{Asymptotic expansion}\label{sec:240617.1526}
\subsection{Stochastic expansion of $Z_n$}\label{sec:240626.1840}
In order to apply the general theory from \cite{nualart2019asymptotic},
we shall write the functional $Z_n$
as a perturbation of a Skorohod integral;
recall that $Z_n=n^\half(S_n-S_\infty)$ is
the rescaled error of the convergence of 
the weighted power variation $S_n$ to its limit $S_\infty$, 
that is defined
at \eqref{eq:240628.1935} and \eqref{eq:240628.1936}, respectively,
and the variation $S_n$ is based on 
the second order difference $\secDiff{n}{j}X$ of the process $X_t$:
\begin{align*}
  \secDiff{n}{j}X=
  % \diff{n}{j+1}X-\diff{n}{j}X=
  X_\tnjp-2X_\tnj+X_\tnjm.
\end{align*}
The following lemma gives the decomposition of $\secDiff{n}{j}X$.
\begin{lemma}\label{lemma:240523.1438}
  For $n\geq2$ and $j\in[n-1]$, $\secDiff{n}{j}X$
  decomposes as follows:
  \begin{align*}
    \secDiff{n}{j}X &= 
    T^{(0)}_{n,j} + T^{(1)}_{n,j} + T^{(2)}_{n,j}
  \end{align*}
  with 
  \begin{align*}
    T^{(0)}_{n,j}&=V^{[1]}_\tnj I_1(\diffker^{n}_{j})
    \\
    % T^{(1)}_{n,j}&=
    % V^{[1,1;1]}_\tnj n^{-2H}+
    % \half V^{[1,1;1]}_\tnj I_2((\bbone^n_{j+1})^{\otimes2}) +
    % \half V^{[1,1;1]}_\tnj I_2((\bbone^n_j)^{\otimes2})
    % \\&\quad+
    % V^{[1,1;2]}_\tnj I_1(\bbone^+_{n,j+1}) \times n^{-1} + 
    % V^{[1,1;2]}_\tnj I_1(\bbone^-_{n,j}) \times n^{-1}
    % \\&\quad+ 
    % V^{[2,1;1]}_\tnj I_1(\bbone^-_{n,j+1}) \times n^{-1} + 
    % V^{[2,1;1]}_\tnj I_1(\bbone^+_{n,j}) \times n^{-1}
    T^{(1)}_{n,j}&=
    \fvoo_\tnj n^{-2H}+
    \half \fvoo_\tnj 
    \cbr{I_2((\bbone^n_{j+1})^{\otimes2}) + I_2((\bbone^n_j)^{\otimes2})}
    \\&\quad+
    \fvot_\tnj 
    \cbr{I_1(\bbone^+_{n,j+1}) + I_1(\bbone^-_{n,j})} \times n^{-1} + 
    \\&\quad+ 
    \fvto_\tnj 
    \cbr{I_1(\bbone^-_{n,j+1}) + I_1(\bbone^+_{n,j})} \times n^{-1}
    \\
    T^{(2)}_{n,j}&=\hat O(n^{-2}\vee n^{-3H}),
  \end{align*}
  where we define 
  \begin{align}
    &\bbone^+_{n,j}(t)=\bbone^n_{j}(t) (t-\tnjm) n,\quad
    \bbone^-_{n,j}(t)=\bbone^n_{j}(t) (\tnj-t) n,
    \label{eq:240606.1219}
    \\
    &\fvoo_t=\dfvo_t \fvo_t,\quad
    \fvot_t=\dfvo_t \fvt_t,\quad
    \fvto_t=\dfvt_t \fvo_t.
    \nn% \label{eq:240606.1731}
  \end{align}
\end{lemma}
% \noindent
Note that 
$\bbone^n_j=\bbone_\sbr{\tnjm,\tnj}$ and
$\diffker^n_j=\bbone^n_{j+1}-\bbone^n_j$
correspond to 
the difference and the second-order difference of the fBm $B$, respectively.
The proof of this lemma will be given in Section \ref{sec:240606.1759}.

\vspsm
Since it holds that 
$T^{(1)}_{n,j}+T^{(2)}_{n,j}=O_M(n^{-2H})$,
the $2k$-th power of the second order difference is written as follows:
\begin{align*}
  (\secDiff{n}{j}X)^{2k} &= 
  (T^{(0)}_{n,j} + T^{(1)}_{n,j} + T^{(2)}_{n,j})^{2k}
  \nn\\&=
  (T^{(0)}_{n,j})^{2k}+
  2k\times (T^{(0)}_{n,j})^{2k-1}\times T^{(1)}_{n,j}+
  2k\times (T^{(0)}_{n,j})^{2k-1}\times T^{(2)}_{n,j}
  \nn\\&\quad+
  (\text{const.})\times 
  (T^{(0)}_{n,j})^{2k-2}\times O(n^{-2H})^2
  \nn\\&=
  (T^{(0)}_{n,j})^{2k}+
  2k\times (T^{(0)}_{n,j})^{2k-1}\times T^{(1)}_{n,j}
  +\hat O(n^{-2kH+H-2}\vee n^{-2kH-2H})
  % \label{eq:240514.1542}
\end{align*}
Then the weighted power variation $S_n$ is decomposed as follows:
\begin{align}
  S_n&=
  n^{2kH-1}\sum_{j=1}^{n-1}
  f(X_\tnj) (\secDiff{n}{j}X)^{2k}
  =% \\&=
  S_n^{(0)} + \sum_{\ell=1}^k S_n^{(1,\ell)}
  +S_n^{(2)}+S_n^{(3)},
  \label{eq:240606.0959}
\end{align}
where we define
\begin{align}
  S_n^{(0)}&=
  n^{2kH-1}\sum_{j=1}^{n-1} f(X_\tnj) 
  (V^{[1]}_\tnj )^{2k}
  \times% \sum_{\ell=0}^{k} 
  \cnstVar{2k}{0} n^{-2Hk}
  \nn\\
  S_n^{(1,\ell)}&=
  n^{2kH-1}\sum_{j=1}^{n-1} f(X_\tnj) 
  (V^{[1]}_\tnj )^{2k}
  \times% \sum_{\ell=0}^{k} 
  \cnstVar{2k}{2\ell} n^{-2H(k-\ell)}
  I_{2\ell}((\diffker^{n}_{j})^{\otimes2\ell})
  \nn\\% \label{eq:240605.1725}\\
  S_n^{(2)}&=
  2k\times n^{2kH-1}\sum_{j=1}^{n-1} f(X_\tnj) \times
  (T^{(0)}_{n,j})^{2k-1}\times T^{(1)}_{n,j}
  \nn\\
  S_n^{(3)}&=\hat O_M(n^{(H-2)\vee(-2H)}).
  % \quad \text{(negligible.)}
  \label{eq:240403.1217}
\end{align}
Here we have used the product formula to expand the factor 
$I_1(\diffker^{n}_{j})^{2k}$ in $(T^{(0)}_{n,j})^{2k}$ and 
the constants 
$\cnstVar{2k}{2\ell}$ ($\ell=0,...,k$) is written as 
\begin{align}\label{eq:240605.1731}
  \cnstVar{2k}{2\ell} = 
  \binom{2k}{2\ell} (2k-2\ell-1)!!
  \times (\innerProdDiffker)^{k-\ell},
\end{align}
where $\innerProdDiffker$ is defined at \eqref{230925.1607}.

The functional $S_n^{(1,\ell)}$ is written as 
\begin{align}
  S_n^{(1,\ell)}&=
  \cnstVar{2k}{2\ell}\times
  n^{2H\ell-1}
  \sum_{j=1}^{n-1} 
  a(X_\tnj) 
  % f(X_\tnj) (V^{[1]}_\tnj )^{2k}\times
  I_{2\ell}((\diffker^{n}_{j})^{\otimes2\ell}),
  \label{eq:240607.1035}
\end{align}
where we introduce the function 
$a(x)=f(x) (V^{[1]}(x))^{2k}$ 
for brevity.
The following lemma shows that
a Skorohod integral works as the principal term of this functional.
\begin{lemma}\label{lemma:240523.1805}
  For $\ell=1,...,k$,
  the functional $S_n^{(1,\ell)}$ is expanded as 
  % has the following decomposition 
  % involving the Skorohod integral of $u_n^{\ell}$:
  \begin{align*}
    S_n^{(1,\ell)}&=
    n^{-\half}\rbr{\delta(u_n^{(\ell)}) + \rsdOne{\ell}},
    % \label{eq:240403.1204}
  \end{align*}
  where we define
  \begin{align}  
    u_n^{(\ell)}&=
    \cnstVar{2k}{2\ell}\times
    n^{2H\ell-\half}
    \sum_{j=1}^{n-1} 
    a(X_\tnj)
    % f(X_\tnj) (V^{[1]}_\tnj )^{2k} \times
    I_{2\ell-1}({\diffker^{n}_{j}}^{\otimes2\ell-1})
    \diffker^{n}_{j}.
    \label{eq:240403.1206}
    \\
    \rsdOne{\ell}&=
    \cnstVar{2k}{2\ell}\times
    n^{2H\ell-\half}
    \sum_{j=1}^{n-1} 
    \abr{D\rbr{a(X_\tnj)},\diffker^{n}_{j}}
    % \abr{D\rbr{f(X_\tnj) (V^{[1]}_\tnj )^{2k}},\diffker^{n}_{j}}
    I_{2\ell-1}({\diffker^{n}_{j}}^{\otimes2\ell-1})
    \nn
    % \\&=
    % \cnstVar{2k}{2\ell}\times
    % n^{2H(\ell-1)-\half}
    % \sum_{j=1}^{n-1} 
    % n^{2H}
    % \abr{D\rbr{f(X_\tnj) (V^{[1]}_\tnj )^{2k}},\diffker^{n}_{j}}
    % I_{2\ell-1}((\diffker^{n}_{j})^{\otimes2\ell-1})
  \end{align}
  and 
  $\rsdOne{\ell}=O_M(n^{-H})$.
\end{lemma}
\begin{proof}
  This decomposition follows from %the definition of $u_n^{(\ell)}$ and
  the basic property of the divergence operator $\delta$ 
  (cf. Proposition 1.3.3 of \cite{nualart2006malliavin}).
  % {\rb ここでpreliminaryを引用する．}

  The residual functional $\rsdOne{\ell}$ is written as 
  \begin{align*}
    \rsdOne{\ell}&=
    % \cnstVar{2k}{2\ell}\times
    % n^{2H\ell-\half}
    % \sum_{j=1}^{n-1} 
    % \abr{D\rbr{a(X_\tnj)},\diffker^{n}_{j}}
    % % \abr{D\rbr{f(X_\tnj) (V^{[1]}_\tnj )^{2k}},\diffker^{n}_{j}}
    % I_{2\ell-1}((\diffker^{n}_{j})^{\otimes2\ell-1})
    % \nn
    % \\&=
    \cnstVar{2k}{2\ell}\times
    n^{2H(\ell-1)-\half}
    \sum_{j=1}^{n-1} 
    n^{2H}
    \abr{D\rbr{a(X_\tnj)},\diffker^{n}_{j}}
    % \abr{D\rbr{f(X_\tnj) (V^{[1]}_\tnj )^{2k}},\diffker^{n}_{j}}
    I_{2\ell-1}({\diffker^{n}_{j}}^{\otimes2\ell-1}).
  \end{align*}
  % The factor 
  % $n^{2H}
  % \abr{D\rbr{f(X_\tnj) (V^{[1]}_\tnj )^{2k}},\diffker^{n}_{j}}$
  % satisfies the condition 
  We can apply the theory of order estimate by the exponent
  to $n^{-(2H(\ell-1)-\half)}\rsdOne{\ell}$.
  Since the corresponding exponent to this factor is calculated as 
  \begin{align*}
    1+(-\half-H(2\ell-1))=\half-H(2\ell-1),
  \end{align*}
  we have 
  \begin{align*}
    \rsdOne{\ell} = n^{2H(\ell-1)-\half}\times O_M(n^{\half-H(2\ell-1)})
    =O_M(n^{-H}).
  \end{align*}  
\end{proof}

For the functional $S_n^{(2)}$,
we have the following lemma:
\begin{lemma}\label{lemma:240523.1809}
  The functional $S_n^{(2)}$ decomposes as 
% \begin{align*}
%   S_n^{(2)}&=
%   \sum_{\ell=1}^k \cbr{S_n^{(2,2;\ell,1)}+S_n^{(2,2;\ell,2)}}+
%   O_M(n^{(-H-\half)\vee(H-2)}),
% \end{align*}
\begin{align*}
  S_n^{(2)}&=
  % \sum_{\ell=1}^k 
  % \cbr{S_n^{(2,2;\ell,1)}+S_n^{(2,2;\ell,2)}}
  % +O_M(n^{-\half-H})
  % +O_M(n^{-H-\half})
  % +O_M(n^{H-2})
  % \\&=
  \sum_{\ell=1}^k 
  \cbr{S_n^{(2,2;\ell,1)}+S_n^{(2,2;\ell,2)}}
  +O_M(n^{(-H-\half)\vee(H-2)})
  % \label{eq:240403.1212}
\end{align*}
where the functionals 
$S_n^{(2,2;\ell,1)}$ and $S_n^{(2,2;\ell,2)}$ are defined at 
\eqref{eq:240403.1215} and \eqref{eq:240403.1216}, respectively,
and they are estimated as 
$S_n^{(2,2;\ell,i)}=O_M(n^{-1})$ $(i=1,2)$.
%
% \begin{align}
%   S_n^{(2,2;\ell,1)}&=
%   2k\times n^{2H\ell-1} 
%   \sum_{j=1}^{n-1} 
%   {\rb f_\tnj A^{(1,2)}_{n,j,\ell}} \times 
%   I_{2\ell+1}((\diffker^{n}_{j})^{\otimes 2\ell-1}
%   \otimes(\bbone^n_j)^{\otimes2})
%   \label{eq:240403.1215}
%   \\
%   S_n^{(2,2;\ell,2)}&=
%   2k\times n^{2H\ell-1} 
%   \sum_{j=1}^{n-1} 
%   {\rb f_\tnj A^{(1,2)}_{n,j,\ell}} \times 
%   I_{2\ell+1}((\diffker^{n}_{j})^{\otimes 2\ell-1}
%   \otimes(\bbone^n_{j+1})^{\otimes2}),
%   \label{eq:240403.1216}
% \end{align}
\end{lemma}

\begin{proof}
We write 
% \begin{align*}
  $T^{(1)}_{n,j}=
  \sum_{i=1}^4 T^{(1,i)}_{n,j}$
% \end{align*}
with
\begin{align*}
  T^{(1,1)}_{n,j}&=
  V^{[1,1;1]}_\tnj n^{-2H},\\
  T^{(1,2)}_{n,j}&=
  \half V^{[1,1;1]}_\tnj 
  \cbr{I_2((\bbone^n_{j+1})^{\otimes2}) + I_2((\bbone^n_j)^{\otimes2})},
  \\
  T^{(1,3)}_{n,j}&=
  V^{[1,1;2]}_\tnj 
  \cbr{I_1(\bbone^+_{n,j+1}+ \bbone^-_{n,j})} \times n^{-1},\\
  % \cbr{I_1(\bbone^+_{n,j+1}) + I_1(\bbone^-_{n,j})} \times n^{-1},\\
  T^{(1,4)}_{n,j}&=
  V^{[2,1;1]}_\tnj 
  \cbr{I_1(\bbone^-_{n,j+1}+ \bbone^+_{n,j})} \times n^{-1}
  % \cbr{I_1(\bbone^-_{n,j+1}) + I_1(\bbone^+_{n,j})} \times n^{-1}
\end{align*}
and 
% For $i=1,..,4$, 
let us define
\begin{align*}
  S_n^{(2,i)}&=
  2k\times n^{2kH-1}\sum_{j=1}^{n-1} f(X_\tnj) \times
  (T^{(0)}_{n,j})^{2k-1}\times T^{(1,i)}_{n,j}.
\end{align*}
By the definition of $S_n^{(2)}$, we have
$% \begin{align*}
  S_n^{(2)}=\sum_{i=1}^4S_n^{(2,i)}.
$ %\end{align*}
By the product formula, we can write
\begin{align}
  I_1(\diffker^{n}_{j})^{2k-1}&=
  % \sum_{\ell=1}^k 
  % a^{2k-1}_{2\ell-1} 
  % I_{2\ell-1}((\diffker^{n}_{j})^{\otimes 2\ell-1})
  % \norm{\diffker^{n}_{j}}^{2(k-\ell)}
  % \\&=
  \sum_{\ell=1}^k 
  \cnstVar{2k-1}{2\ell-1}
  I_{2\ell-1}((\diffker^{n}_{j})^{\otimes 2\ell-1})
  n^{-2H(k-\ell)}
  \label{eq:240605.2109}
\end{align}
with some positive integers $\cnstVar{2k-1}{2\ell-1}$.
% {\blumy with 
% $\cnstVar{2k-1}{2\ell-1}=
% \binom{2k-1}{2\ell-1} (2k-2\ell-1)!!
% \times (\innerProdDiffker)^{k-\ell}$.}
% \begin{align}
%   \cnstVar{2k-1}{2\ell-1} = a^{2k-1}_{2\ell-1} c_0^{k-\ell}.
%   \label{eq:240402.1651}
% \end{align}

%%%%%%%%%%%%%%%%%%%%%%%%%
%%%%% $S_n^{(2,1)}$ %%%%%
%%%%%%%%%%%%%%%%%%%%%%%%%
For $S_n^{(2,1)}$, 
using the above expansion,
we have 
$S_n^{(2,1)}=\sum_{\ell=1}^k S_n^{(2,1;\ell)}$,
% \begin{align*}
%   S_n^{(2,1)}&=
%   % 2k\times n^{2kH-1}\sum_{j=1}^{n-1} f(X_\tnj) \times
%   % (T^{(0)}_{n,j})^{2k-1}\times T^{(1,1)}_{n,j}
%   % \\&=
%   2k\times n^{2kH-1}\sum_{j=1}^{n-1} f(X_\tnj) \times
%   \sum_{\ell=1}^k 
%   \cnstVar{2k}{2\ell-1}
%   (V^{[1]}_\tnj)^{2k-1} V^{[1,1;1]}_\tnj \times 
%   n^{-2H-2H(k-\ell)} \times 
%   I_{2\ell-1}((\diffker^{n}_{j})^{\otimes 2\ell-1})
%   \\&=
%   \sum_{\ell=1}^k S_n^{(2,1;\ell)},
% \end{align*}
where
\begin{align*}
  S_n^{(2,1;\ell)}&=
  2k\times \cnstVar{2k-1}{2\ell-1}
  \times n^{2H(\ell-1)-1} %\times 
  \sum_{j=1}^{n-1} 
  f(X_\tnj) (V^{[1]}_\tnj)^{2k-1} V^{[1,1;1]}_\tnj 
  % f_\tnj (V^{[1]}_\tnj)^{2k-1} V^{[1,1;1]}_\tnj 
  \times 
  I_{2\ell-1}((\diffker^{n}_{j})^{\otimes 2\ell-1}).
\end{align*}
The exponent corresponding to the factor 
$n^{-(2H(\ell-1)-1)}S_n^{(2,1;\ell)}$ is 
$1+(-\half-H(2\ell-1))=\half-H(2\ell-1)$, 
and hence we have 
\begin{align*}
  S_n^{(2,1;\ell)}=
  n^{2H(\ell-1)-1} \times O_M(n^{\half-H(2\ell-1)})=
  O_M(n^{-H-\half}).
\end{align*}

%%%%%%%%%%%%%%%%%%%%%%%%%
%%%%% $S_n^{(2,3)}$ %%%%%
%%%%%%%%%%%%%%%%%%%%%%%%%
For $S_n^{(2,3)}$, similarly, we have 
$S_n^{(2,3)}=
\sum_{\ell=1}^k S_n^{(2,3;\ell)}$ with 
\begin{align*}
  S_n^{(2,3;\ell)}&=
  2k\times \cnstVar{2k-1}{2\ell-1}\times
  n^{2H\ell-2}
  \sum_{j=1}^{n-1} f(X_\tnj)
  (V^{[1]}_\tnj)^{2k-1}
  V^{[1,1;2]}_\tnj 
  \times
  I_{2\ell-1}((\diffker^{n}_{j})^{\otimes 2\ell-1})
  I_1(\bbone^+_{n,j+1}+\bbone^-_{n,j}).
  % \cbr{I_1(\bbone^+_{n,j+1}) + I_1(\bbone^-_{n,j})}.
\end{align*}
By another use of the product formula, we have
\begin{align*}
  &I_{2\ell-1}((\diffker^{n}_{j})^{\otimes 2\ell-1})
  I_1(\bbone^+_{n,j+1}+\bbone^-_{n,j})
  \\*&=
  I_{2\ell}((\diffker^{n}_{j})^{\otimes 2\ell-1}
  \otimes(\bbone^+_{n,j+1}+\bbone^-_{n,j}))+
  (2\ell-1)
  I_{2\ell-2}((\diffker^{n}_{j})^{\otimes 2\ell-2})
  \nabr{\diffker^{n}_{j},\bbone^+_{n,j+1}+\bbone^-_{n,j}}
  \\*&=
  I_{2\ell}((\diffker^{n}_{j})^{\otimes 2\ell-1}
  \otimes(\bbone^+_{n,j+1}+\bbone^-_{n,j})),
\end{align*}
where the last equality follows from an elementary observation
$\nabr{\diffker^{n}_{j},\bbone^+_{n,j+1}}=
-\nabr{\diffker^{n}_{j},\bbone^-_{n,j}}$.
Hence we can write 
$S_n^{(2,3;\ell)}=S_n^{(2,3;\ell,1)}+S_n^{(2,3;\ell,2)}$
with
\begin{align*}
  % S_n^{(2,3;\ell)}&=
  % 2k\times \cnstVar{2k-1}{2\ell-1}\times
  % n^{2H\ell-2}
  % \sum_{j=1}^{n-1} f(X_\tnj)
  % (V^{[1]}_\tnj)^{2k-1}
  % V^{[1,1;2]}_\tnj 
  % \times
  % I_{2\ell}((\diffker^{n}_{j})^{\otimes 2\ell-1}
  % \otimes(\bbone^+_{n,j+1}+\bbone^-_{n,j})),
  % \\
  S_n^{(2,3;\ell,1)}&=
  2k\times \cnstVar{2k-1}{2\ell-1}\times
  n^{2H\ell-2}
  \sum_{j=1}^{n-1} f(X_\tnj)
  (V^{[1]}_\tnj)^{2k-1}
  V^{[1,1;2]}_\tnj 
  \times
  I_{2\ell}((\diffker^{n}_{j})^{\otimes 2\ell-1}
  \otimes\bbone^+_{n,j+1})
  \\
  S_n^{(2,3;\ell,2)}&=
  2k\times \cnstVar{2k-1}{2\ell-1}\times
  n^{2H\ell-2}
  \sum_{j=1}^{n-1} f(X_\tnj)
  (V^{[1]}_\tnj)^{2k-1}
  V^{[1,1;2]}_\tnj 
  \times
  I_{2\ell}((\diffker^{n}_{j})^{\otimes 2\ell-1}
  \otimes\bbone^-_{n,j}).
\end{align*}
The exponent for the factor 
$n^{-(2H\ell-2)}S_n^{(2,3;\ell,1)}$ is calculated as
$1+(-1-H(2\ell-1))=-H(2\ell-1)$, and hence 
we have 
\begin{align*}
  S_n^{(2,3;\ell,1)}&=
  n^{2H\ell-2}\times O_M(n^{-H(2\ell-1)})=
  O_M(n^{H-2}).
\end{align*}
Similarly, we have 
$S_n^{(2,3;\ell,2)}=O_M(n^{H-2})$.
Thus it holds that 
$S_n^{(2,3)}=O_M(n^{H-2})$.

%%%%%%%%%%%%%%%%%%%%%%%%%
%%%%% $S_n^{(2,4)}$ %%%%%
%%%%%%%%%%%%%%%%%%%%%%%%%
For $S_n^{(2,4)}$, by a similar argument, it holds that
\begin{align*}
  S_n^{(2,4)}&=O_M(n^{H-2}).
\end{align*}

%%%%%%%%%%%%%%%%%%%%%%%%%
%%%%% $S_n^{(2,2)}$ %%%%%
%%%%%%%%%%%%%%%%%%%%%%%%%

Next we handle $S_n^{(2,2)}$.
Using again \eqref{eq:240605.2109},
we have
$S_n^{(2,2)}=\sum_{\ell=1}^k S_n^{(2,2;\ell)}$
with
\begin{align*}
  S_n^{(2,2;\ell)}&=
  % 2k\times n^{2H\ell-1} 
  % \sum_{j=1}^{n-1} 
  % f_\tnj A^{(1,2)}_{n,j,\ell} \times
  % I_{2\ell-1}((\diffker^{n}_{j})^{\otimes 2\ell-1}) \times
  % \cbr{I_2((\bbone^n_{j+1})^{\otimes2}) + I_2((\bbone^n_j)^{\otimes2})}
  %
  % 2k\times\cnstVar{2k-1}{2\ell-1}\times\half
  % \\&\quad
  % \times n^{2H\ell-1} 
  % \sum_{j=1}^{n-1} 
  % f(X_\tnj) (V^{[1]}_\tnj)^{2k-1} V^{[1,1;1]}_\tnj 
  % % f_\tnj A^{(1,2)}_{n,j,\ell} \times
  % I_{2\ell-1}((\diffker^{n}_{j})^{\otimes 2\ell-1}) \times
  % \cbr{I_2((\bbone^n_{j+1})^{\otimes2}) + I_2((\bbone^n_j)^{\otimes2})},
  % \\&=
  n^{2H\ell-1} 
  \sum_{j=1}^{n-1} 
  A^{(2,2;\ell)}_{n,j}\times
  I_{2\ell-1}((\diffker^{n}_{j})^{\otimes 2\ell-1}) \times
  \cbr{I_2((\bbone^n_{j+1})^{\otimes2}) + I_2((\bbone^n_j)^{\otimes2})},
\end{align*}
where we define 
$A^{(2,2;\ell)}_{n,j}= 
k\times\cnstVar{2k-1}{2\ell-1}\times
f(X_\tnj) (V^{[1]}_\tnj)^{2k-1} V^{[1,1;1]}_\tnj$.
If $\ell\geq2$, by the product formula, we have
\begin{align*}
  I_{2\ell-1}((\diffker^{n}_{j})^{\otimes 2\ell-1})
  I_2((\bbone^n_j)^{\otimes2})
  &=% \\&=
  % (\text{const.})
  C^{(2)}_{2\ell-1,2}\times
  I_{2\ell-3}((\diffker^{n}_{j})^{\otimes 2\ell-3})
  \abr{\diffker^{n}_{j}, \bbone^n_j}^2
  \\&\quad+
  % (\text{const.})
  C^{(1)}_{2\ell-1,2}\times
  I_{2\ell-1}((\diffker^{n}_{j})^{\otimes 2\ell-2}
  \otimes \bbone^n_j)
  \abr{\diffker^{n}_{j}, \bbone^n_j}
  \\&\quad+
  I_{2\ell+1}((\diffker^{n}_{j})^{\otimes 2\ell-1}
  \otimes(\bbone^n_j)^{\otimes2})
  \\
  I_{2\ell-1}((\diffker^{n}_{j})^{\otimes 2\ell-1})
  I_2((\bbone^n_{j+1})^{\otimes2})
  &=% \\&=
  % (\text{const.})
  C^{(2)}_{2\ell-1,2}\times
  I_{2\ell-3}((\diffker^{n}_{j})^{\otimes 2\ell-3})
  \abr{\diffker^{n}_{j}, \bbone^n_{j+1}}^2
  \\&\quad+
  % (\text{const.})
  C^{(1)}_{2\ell-1,2}\times
  I_{2\ell-1}((\diffker^{n}_{j})^{\otimes 2\ell-2}
  \otimes \bbone^n_{j+1})
  \abr{\diffker^{n}_{j}, \bbone^n_{j+1}}
  \\&\quad+
  I_{2\ell+1}((\diffker^{n}_{j})^{\otimes 2\ell-1}
  \otimes(\bbone^n_{j+1})^{\otimes2}),
\end{align*}
and if $\ell=1$,
\begin{align*}
  I_{1}(\diffker^{n}_{j})
  I_2((\bbone^n_j)^{\otimes2})
  &=% \\&=
  2I_{1}(\bbone^n_j)
  \abr{\diffker^{n}_{j}, \bbone^n_j}
  +% \\&\quad+
  I_{3}(\diffker^{n}_{j}\otimes(\bbone^n_j)^{\otimes2})
  \\
  I_{1}(\diffker^{n}_{j})
  I_2((\bbone^n_{j+1})^{\otimes2})
  &=% \\&=
  2I_{1}(\bbone^n_{j+1})
  \abr{\diffker^{n}_{j}, \bbone^n_{j+1}}
  +% \\&\quad+
  I_{3}(\diffker^{n}_{j}\otimes(\bbone^n_{j+1})^{\otimes2}).
\end{align*}
Since we have 
$-\nabr{\diffker^{n}_{j},\bbone^n_j}=
\nabr{\diffker^{n}_{j},\bbone^n_{j+1}}=
2^{-1}\nabr{\diffker^{n}_{j},\diffker^{n}_{j}}$,
% =2^{-1} c_0 n^{-2H}$,
we obtain the following formulas:
% \begin{itembox}[l]{Memo}
%   \begin{itemize}
%     \item 
%     $\nnorm{\diffker^{n}_{j}}^2=\nabr{\diffker^{n}_{j},\diffker^{n}_{j}}=
%     c_0 n^{-2H}$
%     \item $\nabr{\diffker^{n}_{j},\diffker^{n}_{j}}=
%     \nabr{\diffker^{n}_{j},\bbone^n_{j+1}-\bbone^n_{j}}=
%     \nabr{\diffker^{n}_{j},\bbone^n_{j+1}}
%     -\nabr{\diffker^{n}_{j},\bbone^n_j}=
%     2\nabr{\diffker^{n}_{j},\bbone^n_{j+1}}$
%   \end{itemize}  
% \end{itembox}

\noindent
If $\ell\geq2$,
\begin{align*}
  &I_{2\ell-1}((\diffker^{n}_{j})^{\otimes 2\ell-1})
  \cbr{I_2((\bbone^n_{j+1})^{\otimes2}) + I_2((\bbone^n_j)^{\otimes2})}
  \\&=% \\&=
  2\times C^{(2)}_{2\ell-1,2}\times
  I_{2\ell-3}((\diffker^{n}_{j})^{\otimes 2\ell-3})
  \nabr{\diffker^{n}_{j},\bbone^n_{j+1}}^2
  % \times\cbr{2^{-1} c_0 n^{-2H}}^2
  \\&\quad+
  C^{(1)}_{2\ell-1,2}\times
  I_{2\ell-1}((\diffker^{n}_{j})^{\otimes 2\ell-1})
  \nabr{\diffker^{n}_{j},\bbone^n_{j+1}}
  % \times2^{-1} c_0 n^{-2H}
  \\&\quad+
  I_{2\ell+1}((\diffker^{n}_{j})^{\otimes 2\ell-1}
  \otimes(\bbone^n_j)^{\otimes2})
  +% \\&\quad+
  I_{2\ell+1}((\diffker^{n}_{j})^{\otimes 2\ell-1}
  \otimes(\bbone^n_{j+1})^{\otimes2}).
  % \\&{\rb=:
  % F_{n,j}^{(\ell,4)}+F_{n,j}^{(\ell,3)}+F_{n,j}^{(\ell,1)}+F_{n,j}^{(\ell,2)}}
\end{align*}
If $\ell=1$,
\begin{align*}
  &I_{1}(\diffker^{n}_{j})
  \cbr{I_2((\bbone^n_j)^{\otimes2})+I_2((\bbone^n_{j+1})^{\otimes2})}
  % \\*&=% \\&=
  % 2I_{1}(\bbone^n_j)
  % \abr{\diffker^{n}_{j}, \bbone^n_j} +
  % 2I_{1}(\bbone^n_{j+1})
  % \abr{\diffker^{n}_{j}, \bbone^n_{j+1}} +
  % I_{3}(\diffker^{n}_{j}\otimes(\bbone^n_j)^{\otimes2}) +
  % I_{3}(\diffker^{n}_{j}\otimes(\bbone^n_{j+1})^{\otimes2})
  \\*&=% \\&=
  2I_{1}(\diffker^{n}_{j})
  \abr{\diffker^{n}_{j},\bbone^n_{j+1}} +
  % \times\cbr{2^{-1} c_0 n^{-2H}}+
  I_{3}(\diffker^{n}_{j}\otimes(\bbone^n_j)^{\otimes2}) +
  I_{3}(\diffker^{n}_{j}\otimes(\bbone^n_{j+1})^{\otimes2}).
  % \\&{\rb=:
  % F_{n,j}^{(\ell,3)}+F_{n,j}^{(\ell,1)}+F_{n,j}^{(\ell,2)}}
\end{align*}

Hence, we define 
$S_n^{(2,2;\ell,i)}$ ($i=1,2,3$) for $\ell\geq1$ and 
$S_n^{(2,2;\ell,4)}$ for $\ell\geq2$ by
\begin{align}
  S_n^{(2,2;\ell,1)}&=
  n^{2H\ell-1} 
  \sum_{j=1}^{n-1} 
  A^{(2,2;\ell)}_{n,j}\times
  I_{2\ell+1}((\diffker^{n}_{j})^{\otimes 2\ell-1}
  \otimes(\bbone^n_j)^{\otimes2})
  % I_{2\ell-1}((\diffker^{n}_{j})^{\otimes 2\ell-1}) \times
  % \cbr{I_2((\bbone^n_{j+1})^{\otimes2}) + I_2((\bbone^n_j)^{\otimes2})}
  \label{eq:240403.1215}\\
  S_n^{(2,2;\ell,2)}&=
  n^{2H\ell-1} 
  \sum_{j=1}^{n-1} 
  A^{(2,2;\ell)}_{n,j}\times
  I_{2\ell+1}((\diffker^{n}_{j})^{\otimes 2\ell-1}
  \otimes(\bbone^n_{j+1})^{\otimes2})
  \label{eq:240403.1216}\\
  S_n^{(2,2;\ell,3)}&=
  n^{2H\ell-1} 
  \sum_{j=1}^{n-1} 
  A^{(2,2;\ell)}_{n,j}\times
  C^{(1)}_{2\ell-1,2}\times
  I_{2\ell-1}((\diffker^{n}_{j})^{\otimes 2\ell-1})
  \abr{\diffker^{n}_{j},\bbone^n_{j+1}}
  \nn\\
  S_n^{(2,2;\ell,4)}&=
  n^{2H\ell-1} 
  \sum_{j=1}^{n-1} 
  A^{(2,2;\ell)}_{n,j}\times
  2\times C^{(2)}_{2\ell-1,2}\times
  I_{2\ell-3}((\diffker^{n}_{j})^{\otimes 2\ell-3})
  \abr{\diffker^{n}_{j},\bbone^n_{j+1}}^2
  \nn
\end{align}
and we can write 
\begin{align*}
  S_n^{(2,2;\ell)}&=
  \sum_{i=1}^4 S_n^{(2,2;\ell,i)}
  \tfor \ell\geq2,&
  % \\
  S_n^{(2,2;\ell)}&=
  \sum_{i=1}^3 S_n^{(2,2;\ell,i)}
  \tfor \ell=1.
\end{align*}

By the relation 
$%-\nabr{\diffker^{n}_{j},\bbone^n_j}=
\nabr{\diffker^{n}_{j},\bbone^n_{j+1}}=
2^{-1}\nabr{\diffker^{n}_{j},\diffker^{n}_{j}}=
2^{-1} \innerProdDiffker\times n^{-2H}$,
the functionals $S_n^{(2,2;\ell,i)}$ ($i=1,...,4$) are estimated as follows:
\begin{align*}
  S_n^{(2,2;\ell,1)},S_n^{(2,2;\ell,2)}&=
  n^{2H\ell-1}\times O_M(n^{1+(-1-H\times2\ell)})=
  O_M(n^{-1})
  \\
  S_n^{(2,2;\ell,3)}&=
  n^{2H\ell-1}\times n^{-2H}\times O_M(n^{1+(-\half-H(2\ell-1))})=
  O_M(n^{-\half-H})
  \\
  S_n^{(2,2;\ell,4)}&=
  n^{2H\ell-1}\times n^{-4H}\times O_M(n^{1+(-\half-H(2\ell-3))})=
  O_M(n^{-\half-H})
  \quad\tforsm \ell\geq2.
\end{align*}
Therefore, we obtain 
\begin{align*}
  S_n^{(2,2)}=\sum_{\ell=1}^k S_n^{(2,2;\ell)}
  =\sum_{\ell=1}^k 
  \cbr{S_n^{(2,2;\ell,1)}+S_n^{(2,2;\ell,2)}}
  +O_M(n^{-\half-H})
\end{align*}

\end{proof}

The functional $S_n^{(0)}$ is the principal term of $S_n$, and
written as the following Riemann sum:
\begin{align*}
  S_n^{(0)}&=
  \cnstVar{2k}{0}\times
  n^{-1}
  \sum_{j=1}^{n-1} 
  a(X_\tnj).
\end{align*}
The convergence error of $S_n^{(0)}$ to the corresponding integral,
that is $S_\infty$,
is expanded by the next lemma:
\begin{lemma}\label{lemma:240523.1807}
  It holds that
  \begin{align*}
    S_n^{(0)}-S_\infty&=
    \cnstVar{2k}{0}\times
    (S_n^{(0;1;1)\dagger}
    -S_n^{(\infty;1)\dagger})+
    \hat O_M(n^{-2H})
  \end{align*}
  where
  \begin{align}
    S_n^{(0;1;1)\dagger}&=
    2^{-1} n^{-1}
    \sum_{j=1}^{n-1}
    a'(X_\tnj) V^{[1]}_\tnj \times 
    \cbr{-B\brbr{\bbone^-_{2n,2j+1}} + B\brbr{\bbone^+_{2n,2j}}},
    \label{eq:240606.1236}
    % \label{eq:240523.1821}: S_n^{(0,1)}についてのlabel
    \\
    S_n^{(\infty;1)\dagger}&=
    \frac{1}{2n}%\times
    \brbr{a(X_0)+a(X_1)}
    \label{eq:240606.1237}.
    % \label{eq:240403.1200}: S_n^{(\infty,1)}についてのlabel
  \end{align}
  The two above functionals have the order of $O_M(n^{-1})$.
\end{lemma}
The proof %of this lemma 
will be given in Section \ref{sec:240606.1341}.

\paragraph{Stochastic expansion of $Z_n$.}
%%%%%%%%%%%%%%%%%%%%%%%%%
%%%%% $Z_n$の確率展開 %%%%%
%%%%%%%%%%%%%%%%%%%%%%%%%
Summing up the above three lemmas with \eqref{eq:240606.0959}
and \eqref{eq:240403.1217},
the functional
$Z_n=n^\half (S_n-S_\infty)$ is decomposed as follows:
\begin{align*}
  Z_n&=
  % n^\half(S_n-S_\infty)
  % \\&=
  % n^\half(S_n^{(0)} + \sum_{\ell=1}^k S_n^{(1,\ell)}
  % +S_n^{(2)}+S_n^{(3)}
  % -S_\infty)
  % \\&=
  n^\half(S_n^{(0)}-S_\infty) 
  + \sum_{\ell=1}^k n^\half S_n^{(1,\ell)}
  +n^\half S_n^{(2)}
  +n^\half S_n^{(3)}
  \\&=
  n^\half
  \cnstVar{2k}{0}\times
  (S_n^{(0;1;1)\dagger}
  -S_n^{(\infty;1)\dagger})+
  \hat O_M(n^{\half-2H})
  \\&\quad
  + \sum_{\ell=1}^k \delta(u_n^{(\ell)}) + O_M(n^{-H})
  \\&\quad
  +n^\half
  \sum_{\ell=1}^k 
  \cbr{S_n^{(2,2;\ell,1)}+S_n^{(2,2;\ell,2)}}
  +O_M(n^{(-H)\vee(H-\frac32)})
  \\&\quad
  +\hat O_M(n^{(H-\frac32)\vee(\half-2H)}).
  % +n^\half \hat O_M(n^{(H-2)\vee(-2H)})
\end{align*}
We set $r_n=n^{-\half}$, and we define $u_n$ and $N_n$ by 
\begin{align}
  u_n&=
  \sum_{\ell=1}^k u_n^{(\ell)}
  \label{eq:240606.1421}\\
  N_n&=
  n\times\cnstVar{2k}{0}
  (S_n^{(0;1;1)\dagger}
  -S_n^{(\infty;1)\dagger})
  +n\sum_{\ell=1}^k 
  \cbr{S_n^{(2,2;\ell,1)}+S_n^{(2,2;\ell,2)}}
  +N_n'
  \label{eq:240606.1422}
\end{align}
with some functional 
$N_n'=
% \hat O_M(n^{1-2H})
% +O_M(n^{(\half-H)\vee(H-1)})
% +\hat O_M(n^{(H-1)\vee(1-2H)})=
\hat O_M(n^{(H-1)\vee(\half-H)})$.
Notice that 
$(H-1)\vee(\half-H)<0$ for $H\in(\half,1)$.
Following the framework of the general theory from \cite{nualart2019asymptotic},
we write 
$M_n=\delta(u_n)$
% \begin{align*}
%   M_n=\delta(u_n)
% \end{align*}
for the principal term of $Z_n$.

Therefore,
the functional $Z_n$ decomposes into 
a sum of a Skorohod integral and 
a perturbation functional of a minor order.
%%%%% $Z_n$の確率展開 %%%%%
\begin{proposition}\label{prop:240522.2303}
  The functional $Z_n$ has the following stochastic expansion
  \begin{align*}
    Z_n&=
    M_n + r_nN_n,
    % \delta(u_n) + r_nN_n,
  \end{align*}
  where the functional $N_n$ is estimated as $O_M(1)$.
\end{proposition}
\begin{proof}
  By Lemmas \ref{lemma:240523.1809} and \ref{lemma:240523.1807},
  the functionals
  $S_n^{(0;1;1)\dagger}$, $S_n^{(\infty;1)\dagger}$ and
  $S_n^{(2,2;\ell,i)}$ $(i=1,2)$ are estimated as
  $O_M(n^{-1})$.
\end{proof}

% \newpage
\subsubsection{Order estimate of $u_n$}
The principal term of the functional $Z_n$ is written as the Skorhod integral of 
$u_n$, which is defined at \eqref{eq:240606.1421}.
Since it is essential to hold %the condition 
that the Sobolev norm $\norm{u_n}_{k,p}$ 
(in the sense of the Malliavin caluculus) is bounded,
we will check it in the following proposition.
\begin{proposition}\label{prop:240607.1051}
  The functional $u_n$ is of $O_M(1)$, that is 
  \begin{align*}
    \norm{u_n}_{k,p}=O(1)
  \end{align*}
  for any $k\geq0$ and $p\geq1$.
\end{proposition}
Before we proceed to the proof of this proposition,
for notational convenience we introduce 
\begin{align}  
  u_n^{(\ell)\dagger}&=
  % \times
  n^{2H\ell-\half}
  \sum_{j=1}^{n-1} 
  a(X_\tnj)
  % f(X_\tnj) (V^{[1]}_\tnj )^{2k} \times
  I_{2\ell-1}({\diffker^{n}_{j}}^{\otimes2\ell-1})
  \diffker^{n}_{j}
  \label{eq:240606.1927}
  \\
  \unlconst{\ell}&=\cnstVar{2k}{2\ell}.
  \label{eq:240606.1928}
\end{align}
Then we have 
$u_n^{(\ell)}=\unlconst{\ell} u_n^{(\ell)\dagger}$ and
\begin{align*}
  u_n&=
  \sum_{\ell=1}^k \unlconst{\ell} u_n^{(\ell)\dagger}
\end{align*}
\begin{proof}
  It is enough to show that 
  \begin{align*}
    \norm{\bnorm{D^i u_n^{(\ell)\dagger}}_{\calh^{\otimes i+1}}}_p
    =O(1)
  \end{align*}
  for any $p\geq2$ and $i\geq0$.

  The $\calh^{\otimes i+1}$-valued random variable $D^i u_n^{(\ell)\dagger}$
  is decomposed as follows:
  \begin{align*}
    D^i u_n^{(\ell)\dagger}&=
    n^{2\ell H-\half}\sum_{j=1}^{n-1} 
    D^i \cbr{a(X_\tnj) I_{2\ell-1}({\diffker^n_{j}}^{\otimes 2\ell-1})}
    \otimes \diffker^n_{j}
    \\&=
    n^{2\ell H-\half}\sum_{j=1}^{n-1} \cbr{
    \sum_{q=0}^i \binom{i}{q} 
    D^{i-q}a(X_\tnj) \widetilde\otimes 
    D^q I_{2\ell-1}({\diffker^n_{j}}^{\otimes 2\ell-1})}
    \otimes \diffker^n_{j}
    \\&=
    \sum_{q=0}^i \binom{i}{q}
    n^{2\ell H-\half}\sum_{j=1}^{n-1}\cbr{
    D^{i-q}a(X_\tnj) \widetilde\otimes 
    D^q I_{2\ell-1}({\diffker^n_{j}}^{\otimes 2\ell-1})}
    \otimes \diffker^n_{j}
    \\&=:
    \sum_{q=0}^i \binom{i}{q}
    \widetilde\cali_n^{(\ell,i,q)}
  \end{align*}
  We set
  \begin{align*}
    \cali_n^{(\ell,i,q)}&=
    n^{2\ell H-\half}\sum_{j=1}^{n-1}
    D^{i-q}a(X_\tnj) \otimes 
    D^q I_{2\ell-1}({\diffker^n_{j}}^{\otimes 2\ell-1})
    \otimes \diffker^n_{j}.
  \end{align*}
  Since we have 
  % \begin{align*}
  $\bnorm{\widetilde\cali_n^{(\ell,i,q)}}_{\calh^{\otimes i+1}}\leq
  \bnorm{\cali_n^{(\ell,i,q)}}_{\calh^{\otimes i+1}}$,
  % \end{align*}
  it suffices to show that
  \begin{align*}
    \norm{\bnorm{\cali_n^{(\ell,i,q)}}_{\calh^{\otimes i+1}}}_p
    =O(1)
    % \label{eq:240606.2037}
  \end{align*}
  for $q=0,...,i$.
  Note that it holds that
  \begin{align*}
    \norm{\bnorm{\cali_n^{(\ell,i,q)}}_{\calh^{\otimes i+1}}}_p
    =
    \norm{\babr{\cali_n^{(\ell,i,q)},\cali_n^{(\ell,i,q)}}
    _{\calh^{\otimes i+1}}}_{\frac{p}2}^{\half}
    % \label{eq:240606.2038}
  \end{align*}
  for $p\geq2$.

  We can write
  \begin{align*}
    &\babr{\cali_n^{(\ell,i,q)},\cali_n^{(\ell,i,q)}}
    _{\calh^{\otimes i+1}}
    \\&=
    n^{4\ell H-1}
    \sum_{j_1,j_2=1}^{n-1}
    \abr{D^{i-q}a(X_\tnjo), D^{i-q}a(X_\tnjt)}
    \abr{D^q I_{2\ell-1}({\diffker^n_{j_1}}^{\otimes 2\ell-1}),
    D^q I_{2\ell-1}({\diffker^n_{j_2}}^{\otimes 2\ell-1})}
    \abr{\diffker^n_{j_1},\diffker^n_{j_2}}
  \end{align*}
  
  \item [(i)] If $q>2\ell-1$, then % (i.e. $q\geq2\ell$), 
  we have
  $D^q I_{2\ell-1}({\diffker^n_{j}}^{\otimes 2\ell-1})=0$, and hence
  $\babr{\cali_n^{(\ell,i,q)},\cali_n^{(\ell,i,q)}}
  _{\calh^{\otimes i+1}}=0$.

  \item [(ii)] If $q\leq2\ell-1$, then
  \begin{align*}
    D^q I_{2\ell-1}({\diffker^n_{j_1}}^{\otimes 2\ell-1})&=
    C_{\ell,q}\times
    I_{2\ell-1-q}({\diffker^n_{j_1}}^{\otimes 2\ell-1-q})
    {\diffker^n_{j_1}}^{\otimes q}
  \end{align*}
  with some positive integer $C_{\ell,q}$.
  By the product formula, we have
  \begin{align*}
    &I_{2\ell-1-q}({\diffker^n_{j_1}}^{\otimes 2\ell-1-q})
    I_{2\ell-1-q}({\diffker^n_{j_2}}^{\otimes 2\ell-1-q})
    \\&=
    \sum_{r=0}^{2\ell-1-q} C_{\ell,q,r}\times
    I_{2(2\ell-1-q-r)}({\diffker^n_{j_1}}^{\otimes 2\ell-1-q-r}\otimes
    {\diffker^n_{j_2}}^{\otimes 2\ell-1-q-r})
    \abr{\diffker^n_{j_1},\diffker^n_{j_2}}^r
  \end{align*}
  Thus, 
  \begin{align*}
    &\babr{\cali_n^{(\ell,i,q)},\cali_n^{(\ell,i,q)}}
    _{\calh^{\otimes i+1}}
    \\&=
    (C_{\ell,q})^2\times
    \sum_{r=0}^{2\ell-1-q} C_{\ell,q,r}\times
    % \\&\qquad
    n^{4\ell H-1}
    \sum_{j_1,j_2=1}^{n-1}
    % A_{n,j}^{(i-q)}
    \babr{D^{i-q}a(X_\tnjo), D^{i-q}a(X_\tnjt)}_{\calh^{\otimes i-q}}
    \\&\hspace{50pt}\times
    I_{2(2\ell-1-q-r)}({\diffker^n_{j_1}}^{\otimes 2\ell-1-q-r}\otimes
    {\diffker^n_{j_2}}^{\otimes 2\ell-1-q-r})\times
    \abr{\diffker^n_{j_1},\diffker^n_{j_2}}^{r+q+1}
    \\&=:
    (C_{\ell,q})^2\times
    \sum_{r=0}^{2\ell-1-q} C_{\ell,q,r}\times
    \cali_n^{(\ell,i,q,r)}.
  \end{align*}
  % where we define 
  % $A_{n,j}^{(i-q)}:=
  % \babr{D^{i-q}a(X_\tnjo), D^{i-q}a(X_\tnjt)}_{\calh^{\otimes i-q}}$.

  We can estimate the $L^p$-norm of the functional $\cali_n^{(\ell,i,q,r)}$
  using the exponent introduced in Section \ref{sec:240626.1834}.
  \begin{itemize}
    \item If $2(2\ell-1-q-r)>0$,
    the exponent corresponding to the factor 
    $n^{-(4\ell H-1)}\cali_n^{(\ell,i,q,r)}$ is 
    $(1-2H(r+q+1))+(-\half-H(2(2\ell-1-q-r)))=
    \half-4\ell H$.
    Hence we have $\cali_n^{(\ell,i,q,r)}=O_M(n^{-\half})$.
  
    \item If $2(2\ell-1-q-r)=0$, 
    i.e. $2\ell=1+q+r$,
    the exponent for the factor 
    $n^{-(4\ell H-1)}\cali_n^{(\ell,i,q,r)}$ is 
    $1-2H(r+q+1)=1-4\ell H$, and it holds
    $\cali_n^{(\ell,i,q,r)}=O_M(n^{0})$.
  \end{itemize}
  Hence, we obtain  
  \begin{align*}
    &\babr{\cali_n^{(\ell,i,q)},\cali_n^{(\ell,i,q)}}_{\calh^{\otimes i+1}}
    =O_M(1),
  \end{align*}
  for $q\leq2\ell-1$.
\end{proof}

% \newpage

% \newpage
\subsection{Decomposition of functionals}
\label{sec:240628.1844}
We give the decompositions of the functionals 
$\abr{DM_n, u_n}$, % and 
$(D_{u_n})^2M_n$ and 
$\abr{DN_n, u_n}$, 
which will be repeatedly used in the sequel.
\subsubsection{Decomposition of the functional $\abr{DM_n, u_n}$}
Recall that we have defined 
$M_n=\delta(u_n)$ with 
$u_n=\sum_{\ell=1}^k u_n^{(\ell)}$
and $u_n^{(\ell)}$ is defined at \eqref{eq:240403.1206}.
For brevity, we write 
$M_n^{(\ell)}=\delta(u_n^{(\ell)})$
for $\ell=1,...,k$.
In the light of Lemma \ref{lemma:240523.1805}, 
we also denote 
% we also denote the principal part of $\delta(u_n^{(\ell)})=M_n^{(\ell)}$,
% that is $n^\half S_n^{(1,\ell)}$,
% $M_n^{\prime(\ell)}=n^\half S_n^{(1,\ell)}$
by $M_n^{\prime(\ell)}$
the principal part of $M_n^{(\ell)}$,
% the principal part of $M_n^{(\ell)}:=\delta(u_n^{(\ell)})$,
namely
% which is explicitly written as
\begin{align}
  M_n^{\prime(\ell)}
  (:=n^\half S_n^{(1,\ell)})
  &=
  \unlconst{\ell}\times
  % \cnstVar{2k}{2\ell}\times
  n^{2H\ell-\half}
  \sum_{j=1}^{n-1} 
  a(X_\tnj) 
  I_{2\ell}({\diffker^{n}_{j}}^{\otimes2\ell}),
  \label{eq:240612.0901}
\end{align}
where the functional $S_n^{(1,\ell)}$ and 
the constant $\unlconst{\ell}$ are defined at
\eqref{eq:240607.1035}
% \eqref{eq:240605.1725} (or \eqref{eq:240607.1035})
and \eqref{eq:240606.1928}, respectively.
By Lemma \ref{lemma:240523.1805}, 
we have 
\begin{align*}
  M_n^{(\ell)}=
  % \delta(u_n^{(\ell)})(=M_n^{(\ell)})=
  M_n^{\prime(\ell)}+O_M(n^{-H}).
\end{align*}

With the notations defined above,
the functional $\abr{DM_n,u_n}$ is expanded as follows:
\begin{align}
  \abr{DM_n,u_n}&=
  % \sum_{\ell_1,\ell_2=1}^k
  % \babr{DM_n^{(\ell_1)}, u_n^{(\ell_2)}}
  % \nn\\&=
  \sum_{\ell_1,\ell_2=1}^k
  \babr{DM_n^{\prime(\ell_1)}, u_n^{(\ell_2)}}
  +O_M(n^{-H}),
  \label{eq:240607.1234}
\end{align}
where we used Proposition \ref{prop:240607.1051}.
% {\rb 冗長か？[240607.1059]}
%
We write 
$G_n^{(\ell_1,\ell_2)}=
\babr{DM_n^{\prime(\ell_1)}, u_n^{(\ell_2)}}$.
% 
% %
% \begin{itembox}[l]{Memo}
%   \begin{align*}
%     M_n^{\prime(\ell)}
%     (:=n^\half S_n^{(1,\ell)})
%     &=
%     \unlconst{\ell}\times
%     % \cnstVar{2k}{2\ell}\times
%     n^{2H\ell-\half}
%     \sum_{j=1}^{n-1} 
%     a(X_\tnj) 
%     I_{2\ell}({\diffker^{n}_{j}}^{\otimes2\ell}),
%   \end{align*}
%   \begin{align*}  
%     % n^\half
%     % S_n^{(1,\ell)}&=
%     % n^\half
%     % n^{2kH-1}\sum_{j=1}^{n-1} f(X_\tnj) 
%     % (V^{[1]}_\tnj )^{2k}
%     % \times% \sum_{\ell=0}^{k} 
%     % \cnstVar{2k}{2\ell} n^{-2H(k-\ell)}
%     % I_{2\ell}((\diffker^{n}_{j})^{\otimes2\ell})
%     % % \label{eq:240605.1725}
%     % \\&=
%     % \cnstVar{2k}{2\ell} 
%     % n^{2H\ell-\half}
%     % % n^{2kH-\half}n^{-2H(k-\ell)}
%     % \sum_{j=1}^{n-1} a(X_\tnj) 
%     % I_{2\ell}({\diffker^{n}_{j}}^{\otimes2\ell})
%     % \\
%     u_n^{(\ell)\dagger}&=
%     n^{2H\ell-\half}
%     \sum_{j=1}^{n-1} 
%     a(X_\tnj)
%     % f(X_\tnj) (V^{[1]}_\tnj )^{2k} \times
%     I_{2\ell-1}({\diffker^{n}_{j}}^{\otimes2\ell-1})
%     \diffker^{n}_{j}
%     % \label{eq:240606.1927}
%     \\
%     \unlconst{\ell}&=\cnstVar{2k}{2\ell}.
%     % \label{eq:240606.1928}
%   \end{align*}
%   Then we have 
%   $u_n^{(\ell)}=\unlconst{\ell} u_n^{(\ell)\dagger}$ and
% \end{itembox}
% 
Furthermore, we have
\begin{align}
  G_n^{(\ell_1,\ell_2)}&=
  % \babr{DM_n^{\prime(\ell_1)}, u_n^{(\ell_2)}}
  % \nn\\&=
  \Bigg\langle
    D\rbr{\unlconst{\ell_1}\times
    % D\rbr{\cnstVar{2k}{2\ell_1}\times
    n^{2H\ell_1-\half}
    \sum_{j=1}^{n-1} 
    a(X_\tnj) 
    I_{2\ell_1}({\diffker^{n}_{j}}^{\otimes2\ell_1})}, 
  \nn\\&\hspace{50pt}
  \unlconst{\ell_2}\times
  % \cnstVar{2k}{2\ell_2}\times
  n^{2H\ell_2-\half}
  \sum_{j=1}^{n-1} 
  a(X_\tnj) 
  I_{2\ell_2-1}({\diffker^{n}_{j}}^{\otimes2\ell_2-1})
  \diffker^{n}_{j}\Bigg\rangle
  % % 
  % \nn\\&=
  % \unlconst{\ell_1} \unlconst{\ell_2}
  % % \cnstVar{2k}{2\ell_1} \cnstVar{2k}{2\ell_2}
  % n^{2H(\ell_1+\ell_2)-1} 
  % \sum_{j_1,j_2=1}^{n-1}\abr{
  % D\rbr{a(X_\tnjo) %\wtFunc{n}{j_1} 
  % I_{2\ell_1}({\diffker^{n}_{j_1}}^{\otimes2\ell_1})}, 
  % a(X_\tnjt) %\wtFunc{n}{j_2} 
  % I_{2\ell_2-1}({\diffker^{n}_{j_2}}^{\otimes2\ell_2-1})
  % \diffker^{n}_{j_2}}
  % %
  % \nn\\&=
  % \unlconst{\ell_1} \unlconst{\ell_2}
  % n^{2H(\ell_1+\ell_2)-1} 
  % \sum_{j_1,j_2=1}^{n-1}\abr{
  % D{a(X_\tnjo)} I_{2\ell_1}({\diffker^{n}_{j_1}}^{\otimes2\ell_1}), 
  % a(X_\tnjt) 
  % I_{2\ell_2-1}({\diffker^{n}_{j_2}}^{\otimes2\ell_2-1})
  % \diffker^{n}_{j_2}}
  % \nn\\&\quad+
  % \unlconst{\ell_1} \unlconst{\ell_2}
  % n^{2H(\ell_1+\ell_2)-1} 
  % \sum_{j_1,j_2=1}^{n-1}\abr{
  % a(X_\tnjo) 
  % D{I_{2\ell_1}({\diffker^{n}_{j_1}}^{\otimes2\ell_1})}, 
  % a(X_\tnjt) 
  % I_{2\ell_2-1}({\diffker^{n}_{j_2}}^{\otimes2\ell_2-1})
  % \diffker^{n}_{j_2}}
  %
  \nn\\&=
  \unlconst{\ell_1} \unlconst{\ell_2}
  n^{2H(\ell_1+\ell_2)-1} 
  \sum_{j_1,j_2=1}^{n-1}
  \abr{D{a(X_\tnjo)}, \diffker^{n}_{j_2}} a(X_\tnjt) 
  I_{2\ell_1}({\diffker^{n}_{j_1}}^{\otimes2\ell_1})
  I_{2\ell_2-1}({\diffker^{n}_{j_2}}^{\otimes2\ell_2-1})
  \nn\\&\quad+
  \unlconst{\ell_1} \unlconst{\ell_2}
  n^{2H(\ell_1+\ell_2)-1} 
  \sum_{j_1,j_2=1}^{n-1}
  a(X_\tnjo) a(X_\tnjt) 
  \abr{D{I_{2\ell_1}({\diffker^{n}_{j_1}}^{\otimes2\ell_1})}, 
  \diffker^{n}_{j_2}}
  I_{2\ell_2-1}({\diffker^{n}_{j_2}}^{\otimes2\ell_2-1})
  %
  % \nn\\&{\redmy=
  % \unlconst{\ell_1} \unlconst{\ell_2}
  % n^{2H(\ell_1+\ell_2)-1} 
  % \sum_{j_1,j_2=1}^{n-1}
  % \abr{D{a(X_\tnjo)}, \diffker^{n}_{j_2}} a(X_\tnjt) 
  % I_{2\ell_1}({\diffker^{n}_{j_1}}^{\otimes2\ell_1})
  % I_{2\ell_2-1}({\diffker^{n}_{j_2}}^{\otimes2\ell_2-1})}
  % \nn\\&{\redmy\quad+
  % \unlconst{\ell_1} \unlconst{\ell_2}
  % n^{2H(\ell_1+\ell_2)-1} 
  % \sum_{j_1,j_2=1}^{n-1}
  % a(X_\tnjo) a(X_\tnjt) }
  % \nn\\&{\redmy\hspace{50pt}\times
  % 2\ell_1
  % I_{2\ell_1-1}({\diffker^{n}_{j_1}}^{\otimes2\ell_1-1})
  % \abr{\diffker^{n}_{j_1}, 
  % \diffker^{n}_{j_2}}
  % I_{2\ell_2-1}({\diffker^{n}_{j_2}}^{\otimes2\ell_2-1})}
  % 
  \nn\\&=:
  G_n^{(\ell_1,\ell_2;2)} + G_n^{(\ell_1,\ell_2;1)}
  \label{eq:240403.1754}
\end{align}

For $G_n^{(\ell_1,\ell_2;2)}$, we can write
\begin{align*}
  G_n^{(\ell_1,\ell_2;2)}&=
  \unlconst{\ell_1} \unlconst{\ell_2}
  n^{2H(\ell_1+\ell_2)-1-2H} 
  \\&\qquad\times
  \sum_{j_1,j_2=1}^{n-1}
  n^{2H}\abr{D{a(X_\tnjo)}, \diffker^{n}_{j_2}} a(X_\tnjt) 
  I_{2\ell_1}({\diffker^{n}_{j_1}}^{\otimes2\ell_1})
  I_{2\ell_2-1}({\diffker^{n}_{j_2}}^{\otimes2\ell_2-1})
\end{align*}
The exponent for the factor 
$n^{-(2H(\ell_1+\ell_2)-1-2H)} G_n^{(\ell_1,\ell_2;2)}$ is 
$1+(-\half-2\ell_1H) + 1+(-\half-(2\ell_2-1)H)=
1 -(2\ell_1+2\ell_2-1)H$,
and hence by Theorem \ref{thm:240628.1910} we have 
\begin{align}
  G_n^{(\ell_1,\ell_2;2)}=
  n^{2H(\ell_1+\ell_2)-1-2H} \times 
  O_M(n^{1 -(2\ell_1+2\ell_2-1)H})
  =O_M(n^{-H}).
  % O_M(n^{2H(\ell_1+\ell_2)-1-2H + 1 -(2\ell_1+2\ell_2-1)H})
  \label{eq:240607.1235}
\end{align}

The functional $G_n^{(\ell_1,\ell_2;1)}$ can be written as 
\begin{align*}
  G_n^{(\ell_1,\ell_2;1)}&=
  2\ell_1\times
  \unlconst{\ell_1} \unlconst{\ell_2} \times 
  n^{2H(\ell_1+\ell_2)-1} 
  \sum_{j_1,j_2=1}^{n-1}
  a(X_\tnjo) a(X_\tnjt)
  \nn\\&\hspace{50pt}\times
  I_{2\ell_1-1}({\diffker^{n}_{j_1}}^{\otimes2\ell_1-1})
  I_{2\ell_2-1}({\diffker^{n}_{j_2}}^{\otimes2\ell_2-1})
  \abr{\diffker^{n}_{j_1}, \diffker^{n}_{j_2}}
\end{align*}
and we set 
\begin{align*}
  G_n^{\prime(\ell_1,\ell_2)\dagger}&=
  n^{2H(\ell_1+\ell_2)-1} 
  \sum_{j_1,j_2=1}^{n-1}
  a(X_\tnjo) a(X_\tnjt)
  % \nn\\&\hspace{50pt}
  I_{2\ell_1-1}({\diffker^{n}_{j_1}}^{\otimes2\ell_1-1})
  I_{2\ell_2-1}({\diffker^{n}_{j_2}}^{\otimes2\ell_2-1})
  \abr{\diffker^{n}_{j_1}, \diffker^{n}_{j_2}}
\end{align*}
so that we have 
\begin{align}
  G_n^{(\ell_1,\ell_2;1)}&=
  2\ell_1\times
  \unlconst{\ell_1} \unlconst{\ell_2} \times
  G_n^{\prime(\ell_1,\ell_2)\dagger}.
  \label{eq:240607.1237}
\end{align}
% for notation convenience.
%
By the product formula %, we have 
\begin{align}
  &I_{2\ell_1-1}({\diffker^{n}_{j_1}}^{\otimes2\ell_1-1})
  I_{2\ell_2-1}({\diffker^{n}_{j_2}}^{\otimes2\ell_2-1})
  \nn\\&=
  \sum_{m=0}^{(2\ell_1-1)\wedge(2\ell_2-1)}
  \pfqtanconst{\ell_1,\ell_2}{m}\times
  % \cnstPF{\eqref{eq:240403.1648}}{\ell_1,\ell_2,m}
  I_{2(\ell_1+\ell_2-1-m)}({\diffker^{n}_{j_1}}^{\otimes2\ell_1-1-m}
  \otimes{\diffker^{n}_{j_2}}^{\otimes2\ell_2-1-m})
  \abr{\diffker^{n}_{j_1},\diffker^{n}_{j_2}}^m
  \label{eq:240403.1648}
\end{align}
with some positive integer $\pfqtanconst{\ell_1,\ell_2}{m}$,
we obtain another decomposition
% {\rb producing another decomposition:}
% 
\begin{align}
  G_n^{\prime(\ell_1,\ell_2)\dagger}&=
  % n^{2H(\ell_1+\ell_2)-1} 
  % \sum_{j_1,j_2=1}^{n-1}
  % a(X_\tnjo) a(X_\tnjt)
  % % \nn\\&\hspace{50pt}
  % I_{2\ell_1-1}({\diffker^{n}_{j_1}}^{\otimes2\ell_1-1})
  % I_{2\ell_2-1}({\diffker^{n}_{j_2}}^{\otimes2\ell_2-1})
  % \abr{\diffker^{n}_{j_1}, \diffker^{n}_{j_2}}
  % \\&=
  % n^{2H(\ell_1+\ell_2)-1} 
  % \sum_{j_1,j_2=1}^{n-1}
  % a(X_\tnjo) a(X_\tnjt)
  % \nn\\&\hspace{50pt}
  % \sum_{m=0}^{(2\ell_1-1)\wedge(2\ell_2-1)}
  % \pfqtanconst{\ell_1,\ell_2}{m}\times
  % I_{2(\ell_1+\ell_2-1-m)}({\diffker^{n}_{j_1}}^{\otimes2\ell_1-1-m}
  % \otimes{\diffker^{n}_{j_2}}^{\otimes2\ell_2-1-m})
  % \abr{\diffker^{n}_{j_1},\diffker^{n}_{j_2}}^m
  % \abr{\diffker^{n}_{j_1}, \diffker^{n}_{j_2}}
  % \\&=
  % \sum_{m=0}^{(2\ell_1-1)\wedge(2\ell_2-1)}
  % \pfqtanconst{\ell_1,\ell_2}{m}\times
  % n^{2H(\ell_1+\ell_2)-1} 
  % \sum_{j_1,j_2=1}^{n-1}
  % a(X_\tnjo) a(X_\tnjt)
  % \nn\\&\hspace{50pt}
  % I_{2(\ell_1+\ell_2-1-m)}({\diffker^{n}_{j_1}}^{\otimes2\ell_1-1-m}
  % \otimes{\diffker^{n}_{j_2}}^{\otimes2\ell_2-1-m})
  % \abr{\diffker^{n}_{j_1},\diffker^{n}_{j_2}}^{m+1}
  % \\&=
  \sum_{m=0}^{(2\ell_1-1)\wedge(2\ell_2-1)}
  \pfqtanconst{\ell_1,\ell_2}{m}\times
  G_n^{(\ell_1,\ell_2;m)\dagger\dagger},
  \label{eq:240607.1238}
\end{align}
where we define
\begin{align}
  G_n^{(\ell_1,\ell_2;m)\dagger\dagger}&=
  n^{2H(\ell_1+\ell_2)-1} 
  \sum_{j_1,j_2=1}^{n-1}
  a(X_\tnjo) a(X_\tnjt)
  \nn\\&\hspace{50pt}\times
  I_{2(\ell_1+\ell_2-1-m)}({\diffker^{n}_{j_1}}^{\otimes2\ell_1-1-m}
  \otimes{\diffker^{n}_{j_2}}^{\otimes2\ell_2-1-m})
  \abr{\diffker^{n}_{j_1},\diffker^{n}_{j_2}}^{m+1}.
  \label{eq:240607.1334}
\end{align}

Gathering the relations
\eqref{eq:240607.1234}, \eqref{eq:240403.1754}, \eqref{eq:240607.1235},
\eqref{eq:240607.1237} and \eqref{eq:240607.1238},
we obtain 
\begin{align*}
  \abr{DM_n,u_n}&=
  % \sum_{\ell_1,\ell_2=1}^k
  % \babr{DM_n^{\prime(\ell_1)}, u_n^{(\ell_2)}}
  % +O_M(n^{-H})
  % \\&=
  % \sum_{\ell_1,\ell_2=1}^k
  % % \babr{DM_n^{\prime(\ell_1)}, u_n^{(\ell_2)}}
  % \rbr{G_n^{(\ell_1,\ell_2;2)} + G_n^{(\ell_1,\ell_2;1)}}
  % +O_M(n^{-H})
  % \\&=
  % \sum_{\ell_1,\ell_2=1}^k
  % G_n^{(\ell_1,\ell_2;1)}
  % +O_M(n^{-H})
  % \\&=
  % \sum_{\ell_1,\ell_2=1}^k
  % % G_n^{(\ell_1,\ell_2;1)}  G_n^{(\ell_1,\ell_2;1)}&=
  % 2\ell_1\times
  % \unlconst{\ell_1} \unlconst{\ell_2} \times
  % G_n^{\prime(\ell_1,\ell_2)\dagger}
  % +O_M(n^{-H})
  % \\&=
  \sum_{\ell_1,\ell_2=1}^k
  2\ell_1\times
  \unlconst{\ell_1} \unlconst{\ell_2} \times
  % G_n^{\prime(\ell_1,\ell_2)\dagger} G_n^{\prime(\ell_1,\ell_2)\dagger}&=
  \sum_{m=0}^{(2\ell_1-1)\wedge(2\ell_2-1)}
  \pfqtanconst{\ell_1,\ell_2}{m}\times
  G_n^{(\ell_1,\ell_2;m)\dagger\dagger}
  +O_M(n^{-H}).
  % \\&=
  % \sum_{\ell_1,\ell_2=1}^k
  % \sum_{m=0}^{(2\ell_1-1)\wedge(2\ell_2-1)}
  % 2\ell_1\times
  % \unlconst{\ell_1} \unlconst{\ell_2} \times
  % \pfqtanconst{\ell_1,\ell_2}{m}\times
  % G_n^{(\ell_1,\ell_2;m)\dagger\dagger}
  % +O_M(n^{-H})
\end{align*}
Introducing another notation 
\begin{align}
  \gtwoconst{\ell_1,\ell_2,m}&=
  % C^{(\ell_1,\ell_2,m)}&=
  2\ell_1\times
  \unlconst{\ell_1} \unlconst{\ell_2} \times
  \pfqtanconst{\ell_1,\ell_2}{m},
  \label{eq:240607.1332}
\end{align}
we have
\begin{align*}
  \abr{DM_n,u_n}&=
  \sum_{\ell_1,\ell_2=1}^k
  \sum_{m=0}^{(2\ell_1-1)\wedge(2\ell_2-1)}
  \gtwoconst{\ell_1,\ell_2,m}\times
  % 2\ell_1\times
  % \unlconst{\ell_1} \unlconst{\ell_2} \times
  % \pfqtanconst{\ell_1,\ell_2}{m}\times
  G_n^{(\ell_1,\ell_2;m)\dagger\dagger}
  +O_M(n^{-H}).
\end{align*}

To classify the order of the functionals $G_n^{(\ell_1,\ell_2;m)\dagger\dagger}$,
we set 
\begin{align*}
  \Lambda^{(k)}&= \bcbr{(\ell_1,\ell_2,m):
  \ell_1,\ell_2\in[k], 
  m=0,...,(2\ell_1-1)\wedge(2\ell_2-1)}
  \\
  \Lambda^{(k)}_0 &= \bcbr{(\ell_1,\ell_2,m)\in\Lambda^{(k)}:
  \ell_1+\ell_2-1-m=0}
  \\
  \Lambda^{(k)}_+ &= \bcbr{(\ell_1,\ell_2,m)\in\Lambda^{(k)}:
  \ell_1+\ell_2-1-m>0}.
\end{align*}
When 
$(\ell_1,\ell_2,m)\in\Lambda^{(k)}_0$,
the chaos factor in 
$G_n^{(\ell_1,\ell_2;m)\dagger\dagger}$ disappears,
and
when 
$(\ell_1,\ell_2,m)\in\Lambda^{(k)}_+$, 
$G_n^{(\ell_1,\ell_2;m)\dagger\dagger}$ has a chaos factor.
Notice also that
$(\ell_1,\ell_2,m)\in\Lambda^{(k)}_0$ is equivalent to 
the case where
$\ell_1=\ell_2$ and $m=2\ell_1-1$;
$(\ell_1,\ell_2,m)\in\Lambda^{(k)}_+$ is to 
$\ell_1\neq\ell_2$ or, 
$\ell_1=\ell_2$ and $m<2\ell_1-1$.

If $(\ell_1,\ell_2,m)\in\Lambda^{(k)}_+$, then
the exponent for the factor 
$n^{-(2H(\ell_1+\ell_2)-1)} G_n^{(\ell_1,\ell_2;m)\dagger\dagger}$ is 
$1-2H(m+1) + (-\half -H\times2(\ell_1+\ell_2-1-m))=
\half-2H(\ell_1+\ell_2)$, and hence 
we have 
\begin{align*}
  G_n^{(\ell_1,\ell_2;m)\dagger\dagger}&=
  n^{2H(\ell_1+\ell_2)-1} \times O_M(n^{\half-2H(\ell_1+\ell_2)})=
  O_M(n^{-\half}).
\end{align*}
Otherwise, that is, if $(\ell_1,\ell_2,m)\in\Lambda^{(k)}_0$, it holds that
\begin{align}
  G_n^{(\ell_1,\ell_2;m)\dagger\dagger}(=
  G_n^{(\ell_1,\ell_1;2\ell_1-1)\dagger\dagger})&=
  n^{4H\ell_1-1} 
  \sum_{j_1,j_2=1}^{n-1}
  a(X_\tnjo) a(X_\tnjt)
  % \nn\\&\hspace{50pt}\times
  \abr{\diffker^{n}_{j_1},\diffker^{n}_{j_2}}^{2\ell_1}
  =O_M(1),
  % =O_M(n^0)
  \label{eq:240617.2143}
\end{align}
since the exponent for the factor 
$n^{-(4H\ell_1-1)}  G_n^{(\ell_1,\ell_1;2\ell_1-1)\dagger\dagger}$ is 
$1-2H\times 2\ell_1=1-4H\ell_1$.

Collecting the above arguments, we obtain the following lemma:
% Collecting the above arguments, we obtain 
% the following decomposition of $\abr{DM_n,u_n}$:
\begin{lemma}\label{lemma:240523.2129}
  The functional $\abr{DM_n,u_n}$ decomposes as follows:
  % \begin{align*}
  %   \abr{DM_n,u_n}&=
  %   \sum_{\ell_1=1}^k
  %   % \sum_{(\ell_1,\ell_2,m)\in\Lambda^{(k)}_0}
  %   \gtwoconst{\ell_1,\ell_1,2\ell_1-1}\times
  %   G_n^{(\ell_1,\ell_1,2\ell_1-1)\dagger\dagger}
  %   \\&\quad+
  %   \sum_{(\ell_1,\ell_2,m)\in\Lambda^{(k)}_+}
  %   \gtwoconst{\ell_1,\ell_2,m}\times
  %   G_n^{(\ell_1,\ell_2;m)\dagger\dagger}
  %   +O_M(n^{-H})
  % \end{align*}
  \begin{align}
    \abr{DM_n,u_n}&=
    \sum_{(\ell_1,\ell_2,m)\in\Lambda^{(k)}}
    % \sum_{\ell_1,\ell_2=1}^k
    % \sum_{m=0}^{(2\ell_1-1)\wedge(2\ell_2-1)}
    \gtwoconst{\ell_1,\ell_2,m}\times
    G_n^{(\ell_1,\ell_2;m)\dagger\dagger}
    +O_M(n^{-H})
    \label{eq:240607.1752}
  \end{align}
  where 
  $\gtwoconst{\ell_1,\ell_2,m}$ and 
  $G_n^{(\ell_1,\ell_2;m)\dagger\dagger}$ are defined at 
  \eqref{eq:240607.1332} and \eqref{eq:240607.1334}, respectively.

  Furthermore, the following estimates hold:
  \begin{itemize}
    \item if $(\ell_1,\ell_2,m)\in\Lambda^{(k)}_+$, then 
      $G_n^{(\ell_1,\ell_2;m)\dagger\dagger}=
      O_M(n^{-\half}).$
    \item If $(\ell_1,\ell_2,m)\in\Lambda^{(k)}_0$, then 
      $G_n^{(\ell_1,\ell_2;m)\dagger\dagger}=
      O_M(1).$
  \end{itemize}
\end{lemma}

\vspsm
%%%%%%%%%%%%%%%%%%%%%%%%%%%%%%%%%%%%
% \paragraph{Limit of $G_n^{(\ell_1,\ell_2;m)\dagger\dagger}$ for 
% $(\ell_1,\ell_2,m)\in\Lambda^{(k)}_0$}
%%%%%%%%%%%%%%%%%%%%%%%%%%%%%%%%%%%%
For $(\ell_1,\ell_2,m)\in\Lambda^{(k)}_0$,
the functional 
$G_n^{(\ell_1,\ell_2;m)\dagger\dagger}(=
\GnTwoMain)$
is of $O_M(1)$,
as Lemma \ref{lemma:240523.2129} shows.
% Since $(\ell_1,\ell_2,m)\in\Lambda^{(k)}_0$
% means $\ell_1=\ell_2$ and $m=2\ell_1-1$,
% we write 
% $\GnTwoMain$ for $G_n^{(\ell_1,\ell_2;m)\dagger\dagger}$ in this case.
In fact they converge in $\bbD^{k,p}$ to some limits,
which are related to the summands of the asymptotic variance of $Z_n$.
The following lemma specifies the limit functionals.
\begin{lemma}\label{lemma:240611.2313}
  For $\ell_1=1,...,k$,
    the functional 
    $\GnTwoMain$
    % $G_n^{(\ell_1)\dagger}$
    decomposes as follows:
  \begin{align*}
    \GnTwoMain&=
    % G_n^{(\ell_1)\dagger}&=
    \qtrhoconst{\ell_1}\times
    % C^{\widehat\rho}_{\ell_1}
    \int_0^1 a(X_t)^2dt+
    \hat O_M(n^{(-H)\vee(-3/4)})
    % \hat O_M(n^{-H}) + O_M(n^{-3/4})
  \end{align*}
  with 
  \begin{align*}
    \qtrhoconst{\ell_1}=
    % C^{\widehat\rho}_{\ell_1} = 
    \sum_{i\in\bbZ}
    \widehat\rho(i)^{2\ell_1}.
  \end{align*}
\end{lemma}
The proof of this lemma is given in Section \ref{sec:240611.2314}. %{sec:240607.1650}.
Setting
\begin{align}
  C_{G_\infty} &=   
  \sum_{\ell_1=1}^k 
  \gtwoconst{\ell_1,\ell_1,2\ell_1-1} \times 
  \qtrhoconst{\ell_1}
  \label{eq:240628.1937}
  \\
  G_\infty&= 
  C_{G_\infty} \int_0^1 a(X_t)^2dt,
  % C_{G_\infty} \int_0^1 {\myred (A^{(u)}_t)}^2dt
  \label{eq:240630.1538}
\end{align}
we obtain the following decomposition of $\abr{DM_n,u_n}$
as the summary of the above arguments.
\begin{lemma}\label{lemma:240523.2220}
  The following expansion for $\abr{DM_n,u_n}$ holds:
  \begin{align*}
    \abr{DM_n,u_n}&=
    G_\infty
    +% \\&\quad+
    \sum_{(\ell_1,\ell_2,m)\in\Lambda^{(k)}_+}
    \gtwoconst{\ell_1,\ell_2,m}\times
    G_n^{(\ell_1,\ell_2;m)\dagger\dagger}
    +\hat O_M(n^{(-H)\vee(-3/4)})
  \end{align*}
  Furthermore, 
  the following estimates hold:
  \begin{align*}
    G_n^{(\ell_1,\ell_2;m)\dagger\dagger} &= O_M(n^{-\half})
  \end{align*}
  for $(\ell_1,\ell_2,m)\in\Lambda^{(k)}_+$.
\end{lemma}

% \newpage
\subsubsection{Decomposition of the functional $(D_{u_n})^2 M_n$}
\label{sec:240617.2018}
Using the decomposition \eqref{eq:240607.1752} of $D_{u_n}M_n$ 
(Lemma \ref{lemma:240523.2129}),
$(D_{u_n})^2M_n$ is written as 
\begin{align*}
  (D_{u_n})^2M_n&=
  \sum_{(\ell_1,\ell_2,m)\in\Lambda^{(k)}}
  \gtwoconst{\ell_1,\ell_2,m}\times
  D_{u_n}G_n^{(\ell_1,\ell_2;m)\dagger\dagger}
  +O_M(n^{-H}).
  \\&=
  \sum_{(\ell_1,\ell_2,m)\in\Lambda^{(k)}}
  \sum_{\ell_3=1}^k
  \gtwoconst{\ell_1,\ell_2,m}\times\unlconst{\ell_3}
  D_{u_n^{(\ell_3)\dagger}}G_n^{(\ell_1,\ell_2;m)\dagger\dagger}
  +O_M(n^{-H}).
\end{align*}

For $(\ell_1,\ell_2,m)\in\Lambda^{(k)}_0$,
the functional 
$D_{u_n^{(\ell_3)\dagger}}G_n^{(\ell_1,\ell_2;m)\dagger\dagger}$
is written as 
\begin{align*}
  &D_{u_n^{(\ell_3)\dagger}}G_n^{(\ell_1,\ell_2;m)\dagger\dagger}
  \\&=
  % \Bigg\langle
  %   D\Bcbr{n^{4H\ell_1-1} 
  %   \sum_{j_1,j_2=1}^{n-1}
  %   a(X_\tnjo) a(X_\tnjt)
  %   \abr{\diffker^{n}_{j_1},\diffker^{n}_{j_2}}^{2\ell_1}},
  %   % \\&\hspace{70pt}
  %   n^{2H\ell_3-\half}
  %   \sum_{j_3=1}^{n-1} 
  %   a(X_\tnjs)
  %   I_{2\ell_3-1}({\diffker^{n}_{j_3}}^{\otimes2\ell_3-1})
  %   \diffker^{n}_{j_3}
  % \Bigg\rangle
  % \\&=
  % n^{4H\ell_1-1} n^{2H\ell_3-\half}
  % \sum_{j_1,j_2=1}^{n-1}
  % \sum_{j_3=1}^{n-1} 
  % \abr{D\cbr{a(X_\tnjo) a(X_\tnjt)},\diffker^{n}_{j_3}}
  % a(X_\tnjs)
  % I_{2\ell_3-1}({\diffker^{n}_{j_3}}^{\otimes2\ell_3-1})
  % \abr{\diffker^{n}_{j_1},\diffker^{n}_{j_2}}^{2\ell_1}
  % \\&=
  % n^{4H\ell_1+2H\ell_3-\frac32}
  % % n^{4H\ell_1-1} n^{2H\ell_3-\half}
  % \sum_{j_1,j_2,j_3=1}^{n-1}
  % \abr{D\cbr{a(X_\tnjo) a(X_\tnjt)},\diffker^{n}_{j_3}}
  % a(X_\tnjs)
  % I_{2\ell_3-1}({\diffker^{n}_{j_3}}^{\otimes2\ell_3-1})
  % \abr{\diffker^{n}_{j_1},\diffker^{n}_{j_2}}^{2\ell_1}
  % \\&=
  n^{2H\rbr{2\ell_1+\ell_3-1}-\frac32}
  \sum_{j_1,j_2,j_3=1}^{n-1}
  n^{2H}
  \abr{D\cbr{a(X_\tnjo) a(X_\tnjt)},\diffker^{n}_{j_3}}
  a(X_\tnjs)
  I_{2\ell_3-1}({\diffker^{n}_{j_3}}^{\otimes2\ell_3-1})
  \abr{\diffker^{n}_{j_1},\diffker^{n}_{j_2}}^{2\ell_1}
\end{align*}
The exponent corresponding to the factor 
$n^{-(2H\rbr{2\ell_1+\ell_3-1}-\frac32)}\times
D_{u_n^{(\ell_3)\dagger}}G_n^{(\ell_1,\ell_2;m)\dagger\dagger}$ is 
$(1-2H\times2\ell_1) + (1+(-\half-H(2\ell_3-1)))=
\frac32-4H\ell_1 -2H\ell_3 + H$.
Hence we have 
\begin{align*}
  &D_{u_n^{(\ell_3)\dagger}}G_n^{(\ell_1,\ell_2;m)\dagger\dagger}
  =n^{2H\rbr{2\ell_1+\ell_3-1}-\frac32}\times
  O_M(n^{\frac32-4H\ell_1 -2H\ell_3 + H})
  =
  O_M(n^{-H}).
\end{align*}

For $(\ell_1,\ell_2,m)\in\Lambda^{(k)}_+$,
writing 
$q_{\ell_1,m}:=2\ell_1-1-m$, $q_{\ell_2,m}:=2\ell_2-1-m$,
we decompose
the functional 
$D_{u_n^{(\ell_3)\dagger}}G_n^{(\ell_1,\ell_2;m)\dagger\dagger}$
as 
\begin{align*}
  &D_{u_n^{(\ell_3)\dagger}}G_n^{(\ell_1,\ell_2;m)\dagger\dagger}
  % =\abr{DG_n^{(\ell_1,\ell_2;m)\dagger\dagger},
  % u_n^{(\ell_3)\dagger}}
  % \\&=
  % \Bigg\langle
  %   D\{
  %     n^{2H(\ell_1+\ell_2)-1} 
  %     \sum_{j_1,j_2=1}^{n-1}
  %     a(X_\tnjo) a(X_\tnjt)
  %     % \nn\\&\hspace{50pt}\times
  %     I_{q_{\ell_1,m}+q_{\ell_2,m}}({\diffker^{n}_{j_1}}^{\otimes q_{\ell_1,m}}
  %     \otimes{\diffker^{n}_{j_2}}^{\otimes q_{\ell_2,m}})
  %     \abr{\diffker^{n}_{j_1},\diffker^{n}_{j_2}}^{m+1}
  %   \},
  %   \nn\\&\hspace{50pt}
  %   n^{2H\ell_3-\half}
  %   \sum_{j_3=1}^{n-1} 
  %   a(X_\tnjs)
  %   I_{2\ell_3-1}({\diffker^{n}_{j_3}}^{\otimes2\ell_3-1})
  %   \diffker^{n}_{j_3}
  % \Bigg\rangle
  % % \\&=
  % % n^{2H(\ell_1+\ell_2)-1} n^{2H\ell_3-\half}
  % % \sum_{j_1,j_2=1}^{n-1}\sum_{j_3=1}^{n-1} 
  % % a(X_\tnjs)
  % %   \Bigg\langle
  % %   D\{
  % %     a(X_\tnjo) a(X_\tnjt)
  % %     % \nn\\&\hspace{50pt}\times
  % %     I_{q_{\ell_1,m}+q_{\ell_2,m}}({\diffker^{n}_{j_1}}^{\otimes q_{\ell_1,m}}
  % %     \otimes{\diffker^{n}_{j_2}}^{\otimes q_{\ell_2,m}})
  % %   \},
  % %   \nn\\&\hspace{50pt}
  % %   \diffker^{n}_{j_3}
  % % \Bigg\rangle
  % % I_{2\ell_3-1}({\diffker^{n}_{j_3}}^{\otimes2\ell_3-1})
  % % \abr{\diffker^{n}_{j_1},\diffker^{n}_{j_2}}^{m+1}
  \\&=
  n^{2H(\ell_1+\ell_2+\ell_3)-\frac32}
  \sum_{j_1,j_2,j_3=1}^{n-1}
  a(X_\tnjs)
    \Bigg\langle
    D\{
      a(X_\tnjo) a(X_\tnjt)
      % \nn\\&\hspace{50pt}\times
      I_{q_{\ell_1,m}+q_{\ell_2,m}}({\diffker^{n}_{j_1}}^{\otimes q_{\ell_1,m}}
      \otimes{\diffker^{n}_{j_2}}^{\otimes q_{\ell_2,m}})
    \},
    \diffker^{n}_{j_3}
  \Bigg\rangle
  \nn\\&\hspace{50pt}\times
  I_{2\ell_3-1}({\diffker^{n}_{j_3}}^{\otimes2\ell_3-1})
  \abr{\diffker^{n}_{j_1},\diffker^{n}_{j_2}}^{m+1}
  \\&=
  \cali_n^{(\ell,m;1)}+\cali_n^{(\ell,m;2)},
\end{align*}
where we define 
\begin{align}
  %%%%% 1項目 %%%%%
  \cali_n^{(\ell,m;1)}&=
  n^{2H(\ell_1+\ell_2+\ell_3)-\frac32}
  \sum_{j_1,j_2,j_3=1}^{n-1}
  a(X_\tnjo) a(X_\tnjt) a(X_\tnjs) 
  \nn\\*&\hspace{30pt}\times
  \abr{
  D\rbr{I_{q_{\ell_1,m}+q_{\ell_2,m}}({\diffker^{n}_{j_1}}^{\otimes q_{\ell_1,m}}
  \otimes{\diffker^{n}_{j_2}}^{\otimes q_{\ell_2,m}})},
  \diffker^{n}_{j_3}}
  I_{2\ell_3-1}({\diffker^{n}_{j_3}}^{\otimes2\ell_3-1})
  \abr{\diffker^{n}_{j_1},\diffker^{n}_{j_2}}^{m+1}
  \label{eq:240607.1801}
  \\
  %%%%% 2項目 %%%%%
  \cali_n^{(\ell,m;2)}&=
  n^{2H(\ell_1+\ell_2+\ell_3-1)-\frac32}
  \sum_{j_1,j_2,j_3=1}^{n-1}
  a(X_\tnjs)\times
  n^{2H}
  \abr{D\cbr{a(X_\tnjo) a(X_\tnjt)},\diffker^{n}_{j_3}}
  \nn\\&\hspace{50pt}\times
  I_{q_{\ell_1,m}+q_{\ell_2,m}}({\diffker^{n}_{j_1}}^{\otimes q_{\ell_1,m}}
      \otimes{\diffker^{n}_{j_2}}^{\otimes q_{\ell_2,m}})
  I_{2\ell_3-1}({\diffker^{n}_{j_3}}^{\otimes2\ell_3-1})
  \abr{\diffker^{n}_{j_1},\diffker^{n}_{j_2}}^{m+1}.
  \nn
\end{align}
The exponent for the factor 
$n^{-(2H(\ell_1+\ell_2+\ell_3-1)-\frac32)}\times\cali_n^{(\ell,m;2)}$ is 
\begin{align*}
  &\brbr{(1-2H(m+1))+(-\half-H(q_{\ell_1,m}+q_{\ell_2,m}))}+
  1+(-\half-H(2\ell_3-1))
  \\&=1-2H(\ell_1+\ell_2+\ell_3)+H,
\end{align*}
which indicates that
\begin{align*}
  \cali_n^{(\ell,m;2)}=
  n^{2H(\ell_1+\ell_2+\ell_3-1)-\frac32}\times
  O_M(n^{1-2H(\ell_1+\ell_2+\ell_3)+H})=
  O_M(n^{-H-\half}).
\end{align*}

Thus we have as an intermediate step the following expansion:
\begin{align}
  (D_{u_n})^2M_n&=
  % \sum_{(\ell_1,\ell_2,m)\in\Lambda^{(k)}}
  % \sum_{\ell_3=1}^k
  % \gtwoconst{\ell_1,\ell_2,m}\times\unlconst{\ell_3}
  % D_{u_n^{(\ell_3)\dagger}}G_n^{(\ell_1,\ell_2;m)\dagger\dagger}
  % +O_M(n^{-H}).
  % \\&=
  \sum_{(\ell_1,\ell_2,m)\in\Lambda^{(k)}_0}
  \sum_{\ell_3=1}^k
  \gtwoconst{\ell_1,\ell_2,m}\times\unlconst{\ell_3}
  D_{u_n^{(\ell_3)\dagger}}G_n^{(\ell_1,\ell_2;m)\dagger\dagger}
  \nn\\&\quad+
  \sum_{(\ell_1,\ell_2,m)\in\Lambda^{(k)}_+}
  \sum_{\ell_3=1}^k
  \gtwoconst{\ell_1,\ell_2,m}\times\unlconst{\ell_3}
  D_{u_n^{(\ell_3)\dagger}}G_n^{(\ell_1,\ell_2;m)\dagger\dagger}
  +O_M(n^{-H}).
  \nn\\&=
  % \sum_{(\ell_1,\ell_2,m)\in\Lambda^{(k)}_+}
  % \sum_{\ell_3=1}^k
  % \gtwoconst{\ell_1,\ell_2,m}\times\unlconst{\ell_3}
  % % D_{u_n^{(\ell_3)\dagger}}G_n^{(\ell_1,\ell_2;m)\dagger\dagger}
  % \rbr{\cali_n^{(\ell,m;1)}+\cali_n^{(\ell,m;2)}}
  % +O_M(n^{-H}).
  % \\&=
  \sum_{(\ell_1,\ell_2,m)\in\Lambda^{(k)}_+}
  \sum_{\ell_3=1}^k
  \gtwoconst{\ell_1,\ell_2,m}\times\unlconst{\ell_3}
  \cali_n^{(\ell,m;1)}
  +O_M(n^{-H})
  \label{eq:240608.1202}
\end{align}
% {\rb 等号一つで書いても良いかも．}

We further decompose 
the functional $\cali_n^{(\ell,m;1)}$ defined at \eqref{eq:240607.1801}.
{\it%\gb 
Assume that
$q_{\ell_1,m}, q_{\ell_2,m}>0$ for a while.}
Since we have 
\begin{align*}
  &\abr{D\rbr{
    I_{q_{\ell_1,m}+q_{\ell_2,m}}
    \brbr{{\diffker^{n}_{j_1}}^{\otimes q_{\ell_1,m}}
    \otimes{\diffker^{n}_{j_2}}^{\otimes q_{\ell_2,m}}}},
  \diffker^{n}_{j_3}}
  \\&=
  q_{\ell_1,m}~
  I_{q_{\ell_1,m}+q_{\ell_2,m}-1}
    \brbr{{\diffker^{n}_{j_1}}^{\otimes q_{\ell_1,m}-1}
    \otimes{\diffker^{n}_{j_2}}^{\otimes q_{\ell_2,m}}}
  \abr{\diffker^{n}_{j_1},\diffker^{n}_{j_3}}
  \\&\quad+
  q_{\ell_2,m}~
  I_{q_{\ell_1,m}+q_{\ell_2,m}-1}
    \brbr{{\diffker^{n}_{j_1}}^{\otimes q_{\ell_1,m}}
    \otimes{\diffker^{n}_{j_2}}^{\otimes q_{\ell_2,m}-1}}
  \abr{\diffker^{n}_{j_2},\diffker^{n}_{j_3}},
\end{align*}
the functional $\cali_n^{(\ell,m;1)}$ decomposes as
\begin{align*}
  %%%%% 1項目 %%%%%
  \cali_n^{(\ell,m;1)}&=
  % n^{2(\ell_1+\ell_2+\ell_3) H-\frac32}
  % \sum_{j_1,j_2,j_3\in[n-1]}
  % a(X_{\tnjo}) a(X_{\tnjt}) a(X_\tnjs)
  % \\&\hspace{50pt}\times{\blumy
  % \abr{D\rbr{
  %   I_{q_{\ell_1,m}+q_{\ell_2,m}}
  %   \brbr{{\diffker^{n}_{j_1}}^{\otimes q_{\ell_1,m}}
  %   \otimes{\diffker^{n}_{j_2}}^{\otimes q_{\ell_2,m}}}},
  % \diffker^{n}_{j_3}}}
  % \\&\hspace{50pt}\times
  % I_{2\ell_3-1}({\diffker^{n}_{j_3}}^{\otimes2\ell_3-1})
  % \abr{\diffker^{n}_{j_1},\diffker^{n}_{j_2}}^{m+1}
  % \\&=
  % n^{2(\ell_1+\ell_2+\ell_3) H-\frac32}
  % \sum_{j_1,j_2,j_3\in[n-1]}
  % a(X_{\tnjo}) a(X_{\tnjt}) a(X_\tnjs)
  % \\&\hspace{50pt}\times{\blumy
  % q_{\ell_1,m}~
  % I_{q_{\ell_1,m}+q_{\ell_2,m}-1}
  %   \brbr{{\diffker^{n}_{j_1}}^{\otimes q_{\ell_1,m}-1}
  %   \otimes{\diffker^{n}_{j_2}}^{\otimes q_{\ell_2,m}}}
  % \abr{\diffker^{n}_{j_1},\diffker^{n}_{j_3}}}
  % \\&\hspace{50pt}\times
  % I_{2\ell_3-1}({\diffker^{n}_{j_3}}^{\otimes2\ell_3-1})
  % \abr{\diffker^{n}_{j_1},\diffker^{n}_{j_2}}^{m+1}
  % \\&\quad+
  % n^{2(\ell_1+\ell_2+\ell_3) H-\frac32}
  % \sum_{j_1,j_2,j_3\in[n-1]}
  % a(X_{\tnjo}) a(X_{\tnjt}) a(X_\tnjs)
  % \\&\hspace{50pt}\times{\blumy
  % q_{\ell_2,m}~
  % I_{q_{\ell_1,m}+q_{\ell_2,m}-1}
  %   \brbr{{\diffker^{n}_{j_1}}^{\otimes q_{\ell_1,m}}
  %   \otimes{\diffker^{n}_{j_2}}^{\otimes q_{\ell_2,m}-1}}
  % \abr{\diffker^{n}_{j_2},\diffker^{n}_{j_3}}}
  % \\&\hspace{50pt}\times
  % I_{2\ell_3-1}({\diffker^{n}_{j_3}}^{\otimes2\ell_3-1})
  % \abr{\diffker^{n}_{j_1},\diffker^{n}_{j_2}}^{m+1}
  % \\&=:
  q_{\ell_1,m}\times\cali_n^{(\ell,m;1,1)\dagger}+
  q_{\ell_2,m}\times\cali_n^{(\ell,m;1,2)\dagger}
\end{align*}
with
\begin{align}
  \cali_n^{(\ell,m;1,1)\dagger}&=
  n^{2(\ell_1+\ell_2+\ell_3) H-\frac32}
  \sum_{j_1,j_2,j_3\in[n-1]}
  a(X_{\tnjo}) a(X_{\tnjt}) a(X_\tnjs)
  \nn\\&\hspace{50pt}\times
  I_{q_{\ell_1,m}+q_{\ell_2,m}-1}
  \brbr{{\diffker^{n}_{j_1}}^{\otimes q_{\ell_1,m}-1}
  \otimes{\diffker^{n}_{j_2}}^{\otimes q_{\ell_2,m}}}
  I_{2\ell_3-1}({\diffker^{n}_{j_3}}^{\otimes2\ell_3-1})
  \nn\\&\hspace{50pt}\times
  \abr{\diffker^{n}_{j_1},\diffker^{n}_{j_2}}^{m+1}
  \abr{\diffker^{n}_{j_1},\diffker^{n}_{j_3}}
  \nn\\
  \cali_n^{(\ell,m;1,2)\dagger}&=
  n^{2(\ell_1+\ell_2+\ell_3) H-\frac32}
  \sum_{j_1,j_2,j_3\in[n-1]}
  a(X_{\tnjo}) a(X_{\tnjt}) a(X_\tnjs)
  \nn\\&\hspace{50pt}\times
  I_{q_{\ell_1,m}+q_{\ell_2,m}-1}
  \brbr{{\diffker^{n}_{j_1}}^{\otimes q_{\ell_1,m}}
  \otimes{\diffker^{n}_{j_2}}^{\otimes q_{\ell_2,m}-1}}
  I_{2\ell_3-1}({\diffker^{n}_{j_3}}^{\otimes2\ell_3-1})
  \nn\\&\hspace{50pt}\times
  \abr{\diffker^{n}_{j_1},\diffker^{n}_{j_2}}^{m+1}
  \abr{\diffker^{n}_{j_2},\diffker^{n}_{j_3}}.
  \label{eq:240607.1944}
\end{align}

First we consider the functional $\cali_n^{(\ell,m;1,1)\dagger}$.
By the product formula, we have
\begin{align*}
  &I_{q_{\ell_1,m}+q_{\ell_2,m}-1}
  \brbr{{\diffker^{n}_{j_1}}^{\otimes q_{\ell_1,m}-1}
  \otimes{\diffker^{n}_{j_2}}^{\otimes q_{\ell_2,m}}}
  I_{2\ell_3-1}({\diffker^{n}_{j_3}}^{\otimes2\ell_3-1})
  \\&=
  \sum_{\pi_1,\pi_2:(*)}
  \pfqtorConstOne
  % {\redmy C_{\ell,m,\pi}^{(1)}}
  I_{q_{\ell_1,m}+q_{\ell_2,m}+2\ell_3-2-2(\pi_1+\pi_2)}
  \brbr{{\diffker^{n}_{j_1}}^{\otimes q_{\ell_1,m}-1-\pi_1}
  \otimes{\diffker^{n}_{j_2}}^{\otimes q_{\ell_2,m}-\pi_2}
  \otimes{\diffker^{n}_{j_3}}^{\otimes2\ell_3-1-\pi_1-\pi_2}}
  \\&\qquad\times 
  \abr{\diffker^{n}_{j_1},\diffker^{n}_{j_3}}^{\pi_1}
  \abr{\diffker^{n}_{j_2},\diffker^{n}_{j_3}}^{\pi_2},
\end{align*}
where the summation is taken over all $\pi_1,\pi_2$ such that 
$q_{\ell_1,m}-1\geq\pi_1\geq0$, $q_{\ell_2,m}\geq\pi_2\geq0$ and 
$2\ell_3-1\geq\pi_1+\pi_2(\geq0)$.
(For such $(\pi_1,\pi_2)$, the constant $\pfqtorConstOne$ is a positive integer.)
Then
\begin{align}
  \cali_n^{(\ell,m;1,1)\dagger}&=
  \sum_{\pi_1,\pi_2:(*)}
  \pfqtorConstOne\times
  \cali_n^{(\ell,m;1,1,\pi)\dagger\dagger},
  \label{eq:240608.1231}
\end{align}
where we define
\begin{align*}
  \cali_n^{(\ell,m;1,1,\pi)\dagger\dagger}&=
  n^{2(\ell_1+\ell_2+\ell_3) H-\frac32}
  \sum_{j_1,j_2,j_3\in[n-1]}
  a(X_{\tnjo}) a(X_{\tnjt}) a(X_\tnjs)
  \\*&\hspace{20pt}\times
  I_{q_{\ell_1,m}+q_{\ell_2,m}+2\ell_3-2-2(\pi_1+\pi_2)}
  \brbr{{\diffker^{n}_{j_1}}^{\otimes q_{\ell_1,m}-1-\pi_1}
  \otimes{\diffker^{n}_{j_2}}^{\otimes q_{\ell_2,m}-\pi_2}
  \otimes{\diffker^{n}_{j_3}}^{\otimes2\ell_3-1-\pi_1-\pi_2}}
  \\*&\hspace{50pt}\times
  \abr{\diffker^{n}_{j_1},\diffker^{n}_{j_2}}^{m+1}
  \abr{\diffker^{n}_{j_1},\diffker^{n}_{j_3}}^{\pi_1+1}
  \abr{\diffker^{n}_{j_2},\diffker^{n}_{j_3}}^{\pi_2}.
\end{align*}
By the argument of exponents, 
we can estimate the order of $\cali_n^{(\ell,m;1,1,\pi)\dagger\dagger}$
as follows:
\begin{itemize}
  \item If $q_{\ell_1,m}+q_{\ell_2,m}+2\ell_3-2-2(\pi_1+\pi_2)>0$,
  it holds that
  $\cali_n^{(\ell,m;1,1,\pi)\dagger\dagger}=O_M(n^{-1})$.

  \item If $q_{\ell_1,m}+q_{\ell_2,m}+2\ell_3-2-2(\pi_1+\pi_2)=0$,
  we have
  $\cali_n^{(\ell,m;1,1,\pi)\dagger\dagger}=O_M(n^{-\half})$.
\end{itemize}
Furthermore, we can observe the following relations: 
\begin{itemize}
  \item If
  $q_{\ell_1,m}+q_{\ell_2,m}\neq2\ell_3$, then
  $q_{\ell_1,m}+q_{\ell_2,m}+2\ell_3-2-2(\pi_1+\pi_2)>0$
  holds for any $\pi_1,\pi_2$ satisfying $(*)$.

  \item If
  $q_{\ell_1,m}+q_{\ell_2,m}=2\ell_3$, then
  the condition that $\pi_1=q_{\ell_1,m}-1$ and $\pi_2=q_{\ell_2,m}$ is
  equivalent to 
  $q_{\ell_1,m}+q_{\ell_2,m}+2\ell_3-2-2(\pi_1+\pi_2)=0$
  for any $\pi_1,\pi_2$ satisfying $(*)$.
  We denote $(\pi_1,\pi_2)$ in this case by $\pi_{\ell,m}^{(1)}$.
\end{itemize}

Here we introduce the following notation for later convenience.
% For $(\ell_1,\ell_2,\ell_3,m)$ satisfying 
For $(\ell_1,\ell_2,m)\in\Lambda^{(k)}_+$ and $\ell_3\in[k]$ satisfying 
$q_{\ell_1,m}+q_{\ell_2,m}=2\ell_3$,
we define $\caliqtor$ by 
\begin{align}
  \caliqtor&=
  n^{2(\ell_1+\ell_2+\ell_3) H-\frac32}
  \sum_{j_1,j_2,j_3\in[n-1]}
  a(X_{\tnjo}) a(X_{\tnjt}) a(X_\tnjs)
  \nn\\*&\hspace{50pt}\times
  \abr{\diffker^{n}_{j_1},\diffker^{n}_{j_2}}^{m+1}
  \abr{\diffker^{n}_{j_1},\diffker^{n}_{j_3}}^{q_{\ell_1,m}}
  \abr{\diffker^{n}_{j_2},\diffker^{n}_{j_3}}^{q_{\ell_2,m}}.
  \label{eq:240610.1145}
\end{align}

Then we can summarize the expansion of $\cali_n^{(\ell,m;1,1)\dagger}$:

\noindent
(i) If $q_{\ell_1,m}+q_{\ell_2,m}\neq2\ell_3$, then
\begin{align*}
  \cali_n^{(\ell,m;1,1)\dagger}&=
  \sum_{\pi_1,\pi_2:(*)}
  \pfqtorConstOne\times
  \cali_n^{(\ell,m;1,1,\pi)\dagger\dagger}
  =O_M(n^{-1})
\end{align*}

\noindent
(ii) Assume that $q_{\ell_1,m}+q_{\ell_2,m}=2\ell_3$.
We can observe that 
\begin{align*}
  \cali_n^{(\ell,m;1,1,\pi_{\ell,m}^{(1)})\dagger\dagger}=\caliqtor
  =O_M(n^{-\half}).
\end{align*}
Thus by the decomposition of 
$\cali_n^{(\ell,m;1,1)\dagger}$ \eqref{eq:240608.1231}, 
we have 
\begin{align*}
  \cali_n^{(\ell,m;1,1)\dagger}&=
  % \sum_{\pi_1,\pi_2:(*)}
  % \pfqtorConstOne\times
  % \cali_n^{(\ell,m;1,1,\pi)\dagger\dagger}
  % \\&=
  \pfqtorConstOne[\pi_{\ell,m}^{(1)}]\times
  \cali_n^{(\ell,m;1,1,\pi_{\ell,m}^{(1)})\dagger\dagger}
  +
  \sum_{\substack{\pi_1,\pi_2:(*)\\\pi\neq\pi_{\ell,m}^{(1)}}}
  \pfqtorConstOne\times
  \cali_n^{(\ell,m;1,1,\pi)\dagger\dagger}
  \\&=
  \pfqtorConstOne[\pi_{\ell,m}^{(1)}]\times
  \caliqtor
  % \cali_n^{(\ell,m;1,1,\pi_{\ell,m}^{(1)})\dagger\dagger}
  +O_M(n^{-1}).
\end{align*}

% When $q_{\ell_1,m}+q_{\ell_2,m}=2\ell_3$, and also
% if the condition that 
% $\pi_1=q_{\ell_1,m}-1$ and $\pi_2=q_{\ell_2,m}$ is satisfied,
% the functional $\cali_n^{(\ell,m;1,1,\pi)\dagger\dagger}$ is written as 
% \begin{align*}
%   \cali_n^{(\ell,m;1,1,\pi)\dagger\dagger}&=
%   n^{2(\ell_1+\ell_2+\ell_3) H-\frac32}
%   \sum_{j_1,j_2,j_3\in[n-1]}
%   a(X_{\tnjo}) a(X_{\tnjt}) a(X_\tnjs)
%   \\*&\hspace{50pt}\times
%   \abr{\diffker^{n}_{j_1},\diffker^{n}_{j_2}}^{m+1}
%   \abr{\diffker^{n}_{j_1},\diffker^{n}_{j_3}}^{q_{\ell_1,m}}
%   \abr{\diffker^{n}_{j_2},\diffker^{n}_{j_3}}^{q_{\ell_2,m}}.
% \end{align*}

A similar argument for $\cali_n^{(\ell,m;1,2)\dagger}$ shows the following expansion:

\noindent
(i) If $q_{\ell_1,m}+q_{\ell_2,m}\neq2\ell_3$, then
\begin{align*}
  \cali_n^{(\ell,m;1,2)\dagger}&=O_M(n^{-1})
\end{align*}

\noindent
(ii) If $q_{\ell_1,m}+q_{\ell_2,m}=2\ell_3$, then
we can write
\begin{align*}
  \cali_n^{(\ell,m;1,2)\dagger}&=
  \pfqtorConstTwo[\pi_{\ell,m}^{(2)}]\times
  \caliqtor
  +O_M(n^{-1})
\end{align*}
with some positive constant
$\pfqtorConstTwo[\pi_{\ell,m}^{(2)}]$ 
appearing from the product formula
for the two chaos factors in \eqref{eq:240607.1944}.
Here 
$\pi_{\ell,m}^{(2)}=(q_{\ell_1,m},q_{\ell_2,m}-1)$.

Thus we obtain 

\noindent
(i) if $q_{\ell_1,m}+q_{\ell_2,m}\neq2\ell_3$,
\begin{align}
  \cali_n^{(\ell,m;1)}&=
  q_{\ell_1,m}\times\cali_n^{(\ell,m;1,1)\dagger}+
  q_{\ell_2,m}\times\cali_n^{(\ell,m;1,2)\dagger}
  =O_M(n^{-1})
  \label{eq:240608.1203}
\end{align}

\noindent
(ii) if $q_{\ell_1,m}+q_{\ell_2,m}=2\ell_3$, 
\begin{align}
  \cali_n^{(\ell,m;1)}&=
  (q_{\ell_1,m}\times
  \pfqtorConstOne[\pi_{\ell,m}^{(1)}]+
  % {\redmy C_{\ell,m,\pi_{\ell,m}^{(1)}}^{(1)}}+
  q_{\ell_2,m}\times
  % {\redmy C_{\ell,m,\pi_{\ell,m}^{(2)}}^{(2)}})\times
  \pfqtorConstTwo[\pi_{\ell,m}^{(2)}])\times
  \caliqtor
  +O_M(n^{-1}).
  \label{eq:240608.1204}
\end{align}

Similar arguments work for the case where 
$q_{\ell_1,m}>0, q_{\ell_2,m}=0$ or
$q_{\ell_1,m}=0, q_{\ell_2,m}>0$.
(Note that $q_{\ell_1,m}+q_{\ell_2,m}>0$ holds for 
$(\ell_1,\ell_2,m)\in\Lambda^{(k)}_+$.)

\noindent (i) 
If $q_{\ell_1,m}+q_{\ell_2,m}\neq2\ell_3$, then
it holds that 
\begin{align}
  \cali_n^{(\ell,m;1)}=O_M(n^{-1}).
  \label{eq:240608.1205}
\end{align}

\noindent (ii)
When $q_{\ell_1,m}+q_{\ell_2,m}=2\ell_3$, 
if $q_{\ell_1,m}>0, q_{\ell_2,m}=0$, then 
\begin{align}
  \cali_n^{(\ell,m;1)}&=
  q_{\ell_1,m}\times
  \pfqtorConstOne[\pi_{\ell,m}^{(1)}]\times
  \caliqtor
  +O_M(n^{-1}),
  \label{eq:240608.1206}
\end{align}
and if $q_{\ell_1,m}=0, q_{\ell_2,m}>0$, then
\begin{align}
  \cali_n^{(\ell,m;1)}&=
  q_{\ell_2,m}\times
  \pfqtorConstTwo[\pi_{\ell,m}^{(2)}]\times
  \caliqtor
  +O_M(n^{-1}).
  \label{eq:240608.1207}
\end{align}

To classify the cases for $(\ell_1,\ell_2,\ell_3,m)$, we define the following sets:
\begin{itemize}
  \item $\Lambda^{(k;3)}=
  \cbr{(\ell,m) \text{ such that }
  % \cbr{(\ell_1,\ell_2,\ell_3,m) \text{ such that }
  (\ell_1,\ell_2,m)\in\Lambda^{(k)}_+,
  \ell_3=1,...,k}$
  
  \item $\Lambda^{(k;3)}_{-}=
  \cbr{(\ell,m)\in\Lambda^{(k;3)}\mid 
  % \cbr{(\ell_1,\ell_2,\ell_3,m)\in\Lambda^{(k;3)}\mid 
  2\ell_3\neq q_{\ell_1,m}+q_{\ell_2,m}}$

  \item $\Lambda^{(k;3)}_{\sharp}=
  \cbr{(\ell,m)\in\Lambda^{(k;3)}\mid 
  % \cbr{(\ell_1,\ell_2,\ell_3,m)\in\Lambda^{(k;3)}\mid 
  2\ell_3=q_{\ell_1,m}+q_{\ell_2,m}}$

  \begin{itemize}
    \item $\Lambda^{(k;3)}_{\sharp,0}=
    \cbr{(\ell,m)\in\Lambda^{(k;3)}_{\sharp}\mid 
    % \cbr{(\ell_1,\ell_2,\ell_3,m)\in\Lambda^{(k;3)}_{\sharp}\mid 
    q_{\ell_1,m}>0,q_{\ell_2,m}>0}$
    % {\rb (記号は暫定; $\Lambda^{(k;3)}_{\sharp,3}$の方がいいか？)}

    \item $\Lambda^{(k;3)}_{\sharp,1}=
    \cbr{(\ell,m)\in\Lambda^{(k;3)}_{\sharp}\mid 
    % \cbr{(\ell_1,\ell_2,\ell_3,m)\in\Lambda^{(k;3)}_{\sharp}\mid 
    q_{\ell_1,m}>0,q_{\ell_2,m}=0}$

    \item $\Lambda^{(k;3)}_{\sharp,2}=
    \cbr{(\ell,m)\in\Lambda^{(k;3)}_{\sharp}\mid 
    % \cbr{(\ell_1,\ell_2,\ell_3,m)\in\Lambda^{(k;3)}_{\sharp}\mid 
    q_{\ell_1,m}=0,q_{\ell_2,m}>0}$
  \end{itemize}
\end{itemize}
Obviously, we have 
\begin{align*}
  \Lambda^{(k;3)}=\Lambda^{(k;3)}_{-}\sqcup\Lambda^{(k;3)}_{\sharp},\quad
  \Lambda^{(k;3)}_{\sharp}=
  \Lambda^{(k;3)}_{\sharp,0}\sqcup
  \Lambda^{(k;3)}_{\sharp,1}\sqcup\Lambda^{(k;3)}_{\sharp,2}
\end{align*}

Define constants as follows:
\begin{align}
  \qtorconst&=
  % C^{(\ell,m;\sharp)}&=
  \begin{cases}
    \gtwoconst{\ell_1,\ell_2,m}%\times
    \unlconst{\ell_3}\times
    (q_{\ell_1,m}\times
    \pfqtorConstOne[\pi_{\ell,m}^{(1)}]+
    q_{\ell_2,m}\times
    \pfqtorConstTwo[\pi_{\ell,m}^{(2)}]) & 
    \tifsm (\ell,m)\in\Lambda^{(k;3)}_{\sharp,0}
    % \tifsm (\ell_1,\ell_2,\ell_3,m)\in\Lambda^{(k;3)}_{\sharp,0}
    \\
    \gtwoconst{\ell_1,\ell_2,m}%\times
    \unlconst{\ell_3}\times
    q_{\ell_1,m}\times
    \pfqtorConstOne[\pi_{\ell,m}^{(1)}]& 
    \tifsm (\ell,m)\in\Lambda^{(k;3)}_{\sharp,1}
    \\
    \gtwoconst{\ell_1,\ell_2,m}%\times
    \unlconst{\ell_3}\times
    q_{\ell_2,m}\times
    \pfqtorConstTwo[\pi_{\ell,m}^{(2)}]&
    \tifsm (\ell,m)\in\Lambda^{(k;3)}_{\sharp,2}
  \end{cases}
  \label{eq:240630.1831}
\end{align}
By summing up the expansions 
\eqref{eq:240608.1202}, 
\eqref{eq:240608.1203},
\eqref{eq:240608.1205},
\eqref{eq:240608.1204},
\eqref{eq:240608.1206} and 
\eqref{eq:240608.1207},
we obtain the following decomposition of $(D_{u_n})^2M_n$.
\begin{proposition}
  The functional $(D_{u_n})^2M_n$ is expanded as 
  \begin{align}
    (D_{u_n})^2M_n&=
    \sum_{(\ell_1,\ell_2,\ell_3,m)\in\Lambda^{(k;3)}_{\sharp}}
    \qtorconst\times\caliqtor
    +O_M(n^{-H}),
    \label{eq:240610.1143}
  \end{align}
  and it holds that
  $\caliqtor=O_M(n^{-\half})$ 
  for $(\ell_1,\ell_2,\ell_3,m)\in\Lambda^{(k;3)}_{\sharp}$.

  In particular, $(D_{u_n})^2M_n=O_M(n^{-\half})$.
\end{proposition}
% {\rb 証明を書くかは上の最初の補題たちを埋めてから考える．}

\subsubsection{About the functionals $N_n$ and $D_{u_n} N_n$}
% \subsubsection{About the functionals $N_n$ and $\abr{DN_n, u_n}$}
% 
By the definition \eqref{eq:240606.1422} of $N_n$,
we can write 
\begin{align}
  N_n&= 
  N_n^{(1)} + N_n^{(2)}  
  + \sum_{\ell=1}^k \cbr{N_n^{(3,\ell)} + N_n^{(4,\ell)}}
  % \\&\quad
  +\hat O_M(n^{(H-1)\vee(\half-H)})
  \label{eq:240629.2411}
\end{align}
% \begin{align*}
%   N_n&=
%   n\times\cnstVar{2k}{0}
%   (S_n^{(0;1;1)\dagger}
%   -S_n^{(\infty;1)\dagger})
%   +n\sum_{\ell=1}^k 
%   \cbr{S_n^{(2,2;\ell,1)}+S_n^{(2,2;\ell,2)}}
%   % +N_n'
%   +\hat O_M(n^{(H-1)\vee(\half-H)})
%   % \label{eq:240606.1422}
% \end{align*}
with 
\begin{align}
  % N_n^{(1)}&={\redmy
  % % n\times\cnstVar{2k}{0} S_n^{(0;1;1)\dagger}=
  % 2^{-1} \times\cnstVar{2k}{0}
  % \sum_{j=1}^{n-1}
  % a'(X_\tnj) V^{[1]}_\tnj \times 
  % \cbr{-B\brbr{\bbone^-_{2n,2j+1}} + B\brbr{\bbone^+_{2n,2j}}},}
  % \\
  N_n^{(1)}&=
  % n\times\cnstVar{2k}{0} S_n^{(0;1;1)\dagger}=
  2^{-1} \times\cnstVar{2k}{0}
  \sum_{j=1}^{n-1}
  a'(X_\tnj) V^{[1]}_\tnj \times 
  \cbr{-I_1\brbr{\bbone^-_{2n,2j+1}} + I_1\brbr{\bbone^+_{2n,2j}}}
  \label{eq:240629.2412}
  \\
  N_n^{(2)}&=
  % -n\times\cnstVar{2k}{0} S_n^{(\infty;1)\dagger}=
  -2^{-1}\times\cnstVar{2k}{0}
  \brbr{a(X_0)+a(X_1)}
  \label{eq:240629.2413}
  \\
  N_n^{(3,\ell)}&=
  % n S_n^{(2,2;\ell,1)}=
  n^{2H\ell} 
  \sum_{j=1}^{n-1} 
  A^{(2,2;\ell)}_{n,j}\times
  I_{2\ell+1}((\diffker^{n}_{j})^{\otimes 2\ell-1}
  \otimes(\bbone^n_j)^{\otimes2})
  \label{eq:240629.2414}
  \\
  N_n^{(4,\ell)}&=
  % n S_n^{(2,2;\ell,2)}=
  n^{2H\ell} 
  \sum_{j=1}^{n-1} 
  A^{(2,2;\ell)}_{n,j}\times
  I_{2\ell+1}((\diffker^{n}_{j})^{\otimes 2\ell-1}
  \otimes(\bbone^n_{j+1})^{\otimes2}).
  \label{eq:240629.2415}
\end{align}

The order of $N_n$ has been estimated as $N_n=O_M(1)$ in Proposition \ref{prop:240522.2303}.
From the following proposition,
we can see that $D_{u_n} N_n$ is negligible in the subsequent argument. 
\begin{proposition}
  \label{prop:240610.2337} %一旦コメントアウト
  % 
  % \item [(i)] The functional $N_n$ is estimated as 
  % \begin{align*}
  %   N_n=O_M(1).
  % \end{align*}
  % 
  % \item [(ii)] 
  The functional $D_{u_n} N_n$ is estimated as 
  % \begin{align*}
  %   (D_{u_n})^2 N_n=
  %   \hat O_M(n^{-H}) + O_M(n^{-3/4}).
  % \end{align*}
  \begin{align*}
    D_{u_n} N_n&=
    \hat O_M(n^{(H-1)\vee(\half-H)}).
  \end{align*}
\end{proposition}

\begin{proof}
Since the functional $D_{u_n^{(\ell)\dagger}}N_n^{(2)}$ can be written as 
\begin{align*}
  D_{u_n^{(\ell)\dagger}}N_n^{(2)}&=
  % \abr{DN_n^{(2)},u_n^{(\ell)\dagger}}
  % \\&=
  % \abr{DN_n^{(2)},
  % n^{2H\ell-\half}
  % \sum_{j=1}^{n-1} 
  % a(X_\tnj)
  % I_{2\ell-1}({\diffker^{n}_{j}}^{\otimes2\ell-1})
  % \diffker^{n}_{j}}
  % \\&=
  % n^{2H\ell-\half}
  % \sum_{j=1}^{n-1} 
  % a(X_\tnj)
  % \abr{DN_n^{(2)},
  % \diffker^{n}_{j}}
  % I_{2\ell-1}({\diffker^{n}_{j}}^{\otimes2\ell-1})
  % \\&=
  n^{2H(\ell-1)-\half}
  \sum_{j=1}^{n-1} 
  a(X_\tnj)
  n^{2H}
  \abr{DN_n^{(2)},
  \diffker^{n}_{j}}
  I_{2\ell-1}({\diffker^{n}_{j}}^{\otimes2\ell-1}),
\end{align*}
by the estimate using exponent, we have
$D_{u_n^{(\ell)\dagger}}N_n^{(2)}
% =O_M(n^{2H(\ell-1)-\half+1-\half-H(2\ell-1)})
=O_M(n^{-H})$, and hence
\begin{align*}
  D_{u_n}N_n^{(2)}=O_M(n^{-H}).
\end{align*}

Writing $\graphNnThree$ for
the weighted graph corresponding to the functional
$N_n^{(3,\ell)}$,
we have $\barq(\graphNnThree)=2\ell+1$,
that is an odd number.
To estimate the functional 
$D_{u_n^{(\ell_2)}}N_n^{(3,\ell_1)}$,
we can use Proposition \ref{prop:240617.2009} 
with $q=2\ell_2$.
Since $q=2\ell_2$ is an even number and 
$q\neq\barq(\graphNnThree)$,
the case (ii) of Proposition \ref{prop:240617.2009} applies to 
$D_{u_n^{(\ell_2)}}N_n^{(3,\ell_1)}$.
Hence we obtain 
\begin{align*}
  D_{u_n^{(\ell_2)}}N_n^{(3,\ell_1)}=O_M(n^{-\half})
\end{align*}
since  
$N_n^{(3,\ell_1)}=O_M(1)$.
Notice that the weighted graph representing $N_n^{(4,\ell)}$ is 
the same as $\graphNnThree$, 
and we also have 
\begin{align*}
  D_{u_n^{(\ell_2)}}N_n^{(4,\ell_1)}=O_M(n^{-\half}).
\end{align*}

Write $N_n^{(1)}=N_n^{(1,1)}+N_n^{(1,2)}$ with
\begin{align*}
  N_n^{(1,1)}&=
  -2^{-1} \times\cnstVar{2k}{0}
  \sum_{j=1}^{n-1}
  a'(X_\tnj) V^{[1]}_\tnj \times 
  I_1\brbr{\bbone^-_{2n,2j+1}},
  \\
  N_n^{(1,2)}&=
  2^{-1} \times\cnstVar{2k}{0}
  \sum_{j=1}^{n-1}
  a'(X_\tnj) V^{[1]}_\tnj \times 
  I_1\brbr{\bbone^+_{2n,2j}}.
\end{align*}
and denote the weighted graph representing the functionals 
$N_n^{(1,1)}$ and $N_n^{(1,2)}$ by $\graphNnOne$.
The graph $\graphNnOne$ is written as 
$\graphNnOne=
(\cbr{v_0},0, \cbr{(v_0,1)\mapsto1,(v_0,2)\mapsto0})$,
and $\barq(\graphNnOne)=1$.
Thus we can again apply Proposition \ref{prop:240617.2009} (ii)
to obtain 
\begin{align*}
  D_{u_n^{(\ell)}}N_n^{(1,1)},
  D_{u_n^{(\ell)}}N_n^{(1,2)}=O_M(n^{-\half})
\end{align*}
for $\ell=1,...,k$,
since we have 
$e(\graphNnOne)=1+(-1)=0$.
% since we have $N_n^{(1,1)}, N_n^{(1,2)}=O_M(1)$.

\end{proof}

% \newpage

% \newpage
\subsection{The Limit of random symbols}
\label{sec:240628.1845}
Consider the random symbols 
$\mS^{(2,0)}_{0,n}$,
$\mS^{(3,0)}_{n}$,
$\mS^{(1,0)}_{n}$ and
$\mS^{(2,0)}_{1,n}$
(defined at 
\eqref{eq:240629.2341},
\eqref{eq:240629.2342},
\eqref{eq:240629.2343} and
\eqref{eq:240629.2344},
respectively)
for $u_n, r_n, N_n, G_\infty$
in our context.
Here we identify the random symbols 
$\mS_0^{(2,0)}$, $\mS^{(3,0)}$, $\mS^{(1,0)}$ and 
$\mS_1^{(2,0)}$
that will satisfy the condition {\bf[D]} (iii) in Section \ref{sec:240626.1819}.

\subsubsection*{The Limit of $\mS^{(2,0)}_{0,n}$}
To find the limit random symbol of $\mS^{(2,0)}_{0,n}$,
we first use the IBP formula to obtain another random symbol
$\bar\mS^{(2,0)}_{0,n}$ such that
\begin{align*}
  E\sbr{\Psi(\sfz)\mS^{(2,0)}_{0,n}(\sfi\sfz)}=
  E\sbr{\Psi(\sfz)\bar\mS^{(2,0)}_{0,n}(\sfi\sfz)}.
\end{align*}
Recall that the symbol $\mS^{(2,0)}_{0,n}(\sfi\sfz)$ is defined by 
\begin{align*}
  \mS^{(2,0)}_{0,n}(\sfi\sfz)=\half\qtan_n[\sfi\sfz]^2=
  \half
  r_n^{-1}
  \brbr{\abr{DM_n,u_n}-G_\infty} [\sfi\sfz]^2.
\end{align*}
By Lemma \ref{lemma:240523.2220},
the coefficient of $\mS^{(2,0)}_{0,n}(\sfi\sfz)$ is written as 
\begin{align}
  &\half n^{\half}\rbr{\abr{DM_n,u_n}-G_\infty}
  =% \nn\\&=
  \half 
  \sum_{(\ell_1,\ell_2,m)\in\Lambda^{(k)}_+}
  \gtwoconst{\ell_1,\ell_2,m}\times
  n^{\half}
  G_n^{(\ell_1,\ell_2;m)\dagger\dagger}
  +\qtanresid
  % +\hat O_M(n^{(\half-H)\vee(-\frac14)})
  \label{eq:240610.1718}
\end{align}
with some functional
$\qtanresid=\hat O_M(n^{(\half-H)\vee(-\frac14)})$.
Thus we only need to consider the limit of the symbols 
% Thus we only need to consider the {\rb weak limit} of the symbols 
$n^{\half}
G_n^{(\ell_1,\ell_2;m)\dagger\dagger}[\sfi\sfz]^2$
for $(\ell_1,\ell_2,m)\in\Lambda^{(k)}_+$.

We can write
% \begin{align}
%   E\sbr{\Psi(\sfz)\times n^{\half} G_n^{(\ell_1,\ell_2;m)\dagger\dagger}
%   [\sfi\sfz]^2}.
%   \label{eq:240610.1729}
% \end{align}
\begin{align*}
  &E\sbr{\Psi(\sfz)\times n^{\half} G_n^{(\ell_1,\ell_2;m)\dagger\dagger}
  [\sfi\sfz]^2}
  \\&=
  % E\Big[\Psi(\sfz)\times n^{\half} 
  % n^{2H(\ell_1+\ell_2)-1} \Big\{
  % \sum_{j_1,j_2=1}^{n-1}
  % a(X_\tnjo) a(X_\tnjt)
  % \nn\\&\hspace{50pt}\times
  % I_{2(\ell_1+\ell_2-1-m)}({\diffker^{n}_{j_1}}^{\otimes2\ell_1-1-m}
  % \otimes{\diffker^{n}_{j_2}}^{\otimes2\ell_2-1-m})
  % \abr{\diffker^{n}_{j_1},\diffker^{n}_{j_2}}^{m+1}\Big\}
  % [\sfi\sfz]^2\Big]
  % \\&=
  n^{2H(\ell_1+\ell_2)-\half} 
  \sum_{j_1,j_2=1}^{n-1}
  E\Big[\Psi(\sfz)\times  
  % \Big\{
  a(X_\tnjo) a(X_\tnjt)
  % \nn\\&\hspace{50pt}\times
  % I_{2(\ell_1+\ell_2-1-m)}({\diffker^{n}_{j_1}}^{\otimes2\ell_1-1-m}
  % \otimes{\diffker^{n}_{j_2}}^{\otimes2\ell_2-1-m})
  I_{q_{\ell_1,m}+q_{\ell_2,m}}({\diffker^{n}_{j_1}}^{\otimes q_{\ell_1,m}}
  \otimes{\diffker^{n}_{j_2}}^{\otimes q_{\ell_2,m}})%\Big\}
  \Big]
  \nn\\&\hspace{50pt}\times
  \abr{\diffker^{n}_{j_1},\diffker^{n}_{j_2}}^{m+1}
  [\sfi\sfz]^2.
\end{align*}
Here recall that we write 
$q_{\ell_1,m}:=2\ell_1-1-m$ and $q_{\ell_2,m}:=2\ell_2-1-m$.
{\it Assume that $q_{\ell_1,m}\geq1$ for a while.}
% We can write \eqref{eq:240610.1729} explicitly as 
% The above expectation is written explicitly as 
% \begin{align*}
%   &E\sbr{\Psi(\sfz)\times n^{\half} G_n^{(\ell_1,\ell_2;m)\dagger\dagger}
%   [\sfi\sfz]^2}
%   \\&=
%   % E\Big[\Psi(\sfz)\times n^{\half} 
%   % n^{2H(\ell_1+\ell_2)-1} \Big\{
%   % \sum_{j_1,j_2=1}^{n-1}
%   % a(X_\tnjo) a(X_\tnjt)
%   % \nn\\&\hspace{50pt}\times
%   % I_{2(\ell_1+\ell_2-1-m)}({\diffker^{n}_{j_1}}^{\otimes2\ell_1-1-m}
%   % \otimes{\diffker^{n}_{j_2}}^{\otimes2\ell_2-1-m})
%   % \abr{\diffker^{n}_{j_1},\diffker^{n}_{j_2}}^{m+1}\Big\}
%   % [\sfi\sfz]^2\Big]
%   % \\&=
%   n^{2H(\ell_1+\ell_2)-\half} 
%   \sum_{j_1,j_2=1}^{n-1}
%   E\Big[\Psi(\sfz)\times  
%   % \Big\{
%   a(X_\tnjo) a(X_\tnjt)
%   % \nn\\&\hspace{50pt}\times
%   % I_{2(\ell_1+\ell_2-1-m)}({\diffker^{n}_{j_1}}^{\otimes2\ell_1-1-m}
%   % \otimes{\diffker^{n}_{j_2}}^{\otimes2\ell_2-1-m})
%   I_{q_{\ell_1,m}+q_{\ell_2,m}}({\diffker^{n}_{j_1}}^{\otimes q_{\ell_1,m}}
%   \otimes{\diffker^{n}_{j_2}}^{\otimes q_{\ell_2,m}})%\Big\}
%   \Big]
%   \nn\\&\hspace{50pt}\times
%   \abr{\diffker^{n}_{j_1},\diffker^{n}_{j_2}}^{m+1}
%   [\sfi\sfz]^2
% \end{align*}
% 
By the IBP formula, the above expectation in the summand is decomposed as 
\begin{align*}
  &E\sbr{\Psi(\sfz)\times
  a(X_\tnjo) a(X_\tnjt)
  I_{q_{\ell_1,m}+q_{\ell_2,m}}({\diffker^{n}_{j_1}}^{\otimes q_{\ell_1,m}}
  \otimes{\diffker^{n}_{j_2}}^{\otimes q_{\ell_2,m}})}[\sfi\sfz]^2
  \\&=
  \half
  E\Big[\Psi(\sfz)
  \abr{DG_\infty,\diffker^{n}_{j_1}}\times
  a(X_{\tnjo}) a(X_{\tnjt})
  I_{q_{\ell_1,m}+q_{\ell_2,m}-1}({\diffker^{n}_{j_1}}^{\otimes q_{\ell_1,m}-1}
  \otimes{\diffker^{n}_{j_2}}^{\otimes q_{\ell_2,m}})
  \Big] [\sfi\sfz]^4
  \\&\quad+
  E\Big[\Psi(\sfz)
  \abr{D\rbr{a(X_{\tnjo}) a(X_{\tnjt})},\diffker^{n}_{j_1}}
  I_{q_{\ell_1,m}+q_{\ell_2,m}-1}({\diffker^{n}_{j_1}}^{\otimes q_{\ell_1,m}-1}
  \otimes{\diffker^{n}_{j_2}}^{\otimes q_{\ell_2,m}})
  \Big] [\sfi\sfz]^2
\end{align*}
Hence we have 
\begin{align}
  &E\sbr{\Psi(\sfz)
  n^{\half} G_n^{(\ell_1,\ell_2;m)\dagger\dagger}
  [\sfi\sfz]^2}
  % \nn\\&=
  % n^{2(\ell_1+\ell_2) H-\half}
  % \sum_{j_1,j_2\in[n-1]}
  % E\Big[
  % \Psi(\sfz) a(X_{\tnjo}) a(X_{\tnjt})
  % I_{q_{\ell_1,m}+q_{\ell_2,m}}({\diffker^{n}_{j_1}}^{\otimes q_{\ell_1,m}}
  % \otimes{\diffker^{n}_{j_2}}^{\otimes q_{\ell_2,m}})
  % \Big]
  % \abr{\diffker^{n}_{j_1},\diffker^{n}_{j_2}}^{m+1}
  % [\sfi\sfz]^2
  \nn\\&=
  E\Big[\Psi(\sfz) \times\cali^{(\ell_1,\ell_2,m;1)}_{\qtan,n}
  [\sfi\sfz]^4 \Big]+
  E\Big[\Psi(\sfz) \times\cali^{(\ell_1,\ell_2,m;2)}_{\qtan,n}
  [\sfi\sfz]^2 \Big]
  \label{eq:240425.1623}
\end{align}
with 
\begin{align*}
  \cali^{(\ell_1,\ell_2,m;1)}_{\qtan,n}&=
  \half n^{2(\ell_1+\ell_2) H-\half-2H}
  \sum_{j_1,j_2\in[n-1]}
  n^{2H}\abr{DG_\infty,\diffker^{n}_{j_1}}\times 
  a(X_{\tnjo}) a(X_{\tnjt})
  \\&\hspace{50pt}\times
  I_{q_{\ell_1,m}+q_{\ell_2,m}-1}({\diffker^{n}_{j_1}}^{\otimes q_{\ell_1,m}-1}
  \otimes{\diffker^{n}_{j_2}}^{\otimes q_{\ell_2,m}})
  \abr{\diffker^{n}_{j_1},\diffker^{n}_{j_2}}^{m+1}
  \\
  \cali^{(\ell_1,\ell_2,m;2)}_{\qtan,n}&=
  n^{2(\ell_1+\ell_2) H-\half-2H}
  \sum_{j_1,j_2\in[n-1]}
  n^{2H}\abr{D\rbr{a(X_{\tnjo}) a(X_{\tnjt})},\diffker^{n}_{j_1}}
  \\&\hspace{50pt}\times
  I_{q_{\ell_1,m}+q_{\ell_2,m}-1}({\diffker^{n}_{j_1}}^{\otimes q_{\ell_1,m}-1}
  \otimes{\diffker^{n}_{j_2}}^{\otimes q_{\ell_2,m}})
  \abr{\diffker^{n}_{j_1},\diffker^{n}_{j_2}}^{m+1}.
\end{align*}
For $(\ell_1,\ell_2,m)\in\Lambda^{(k)}_+$,
it holds that 
$q_{\ell_1,m}+q_{\ell_2,m}-1=
2(\ell_1+\ell_2-1-m)-1\geq1$.
Since the exponents for the functionals
$n^{-(2(\ell_1+\ell_2) H-\half-2H)}\cali^{(\ell_1,\ell_2,m;1)}_{\qtan,n}$ and 
$n^{-(2(\ell_1+\ell_2) H-\half-2H)}\cali^{(\ell_1,\ell_2,m;2)}_{\qtan,n}$ are both 
\begin{align*}
  1-2H(m+1) + (-\half -H (q_{\ell_1,m}+q_{\ell_2,m}-1))=
  \half -2H(\ell_1+\ell_2) +H,
\end{align*}
we have 
\begin{align*}
  \cali^{(\ell_1,\ell_2,m;1)}_{\qtan,n},\cali^{(\ell_1,\ell_2,m;2)}_{\qtan,n}&=
  n^{2(\ell_1+\ell_2) H-\half-2H}\times O_M(n^{\half -2H(\ell_1+\ell_2) +H})=
  O_M(n^{-H}).
\end{align*}

In the case where $q_{\ell_1,m}=0$ and $q_{\ell_2,m}>0$,
we can obtain a similar decomposition using IBP formula with respect to $\diffker^{n}_{j_2}$.
 
Hence there exist functionals 
$\cali^{(i)}_{\qtan,n}$ $(i=1,2)$
% satisfying
of $O_M(n^{-H})$
% $\cali^{(1)}_{\qtan,n}=O_M(n^{-H}),
% \cali^{(2)}_{\qtan,n}=O_M(n^{-H})$
such that 
% we can find a random symbol $\bar\mS^{(2,0)}_{0,n}(\sfi\sfz)$ such that 
the equality 
\begin{align*}
  E\sbr{\Psi(\sfz)\mS^{(2,0)}_{0,n}(\sfi\sfz)}=
  E\sbr{\Psi(\sfz)\bar\mS^{(2,0)}_{0,n}(\sfi\sfz)}
\end{align*}
holds with 
\begin{align*}
  \bar\mS^{(2,0)}_{0,n}(\sfi\sfz)=
  \cali^{(1)}_{\qtan,n}[\sfi\sfz]^4 +
  \rbr{\cali^{(2)}_{\qtan,n}+\qtanresid}[\sfi\sfz]^2.
\end{align*}

% {\mygreen
% Hence there exist a random symbol $\bar\mS^{(2,0)}_{0,n}(\sfi\sfz)$
% written as 
% \begin{align*}
%   \bar\mS^{(2,0)}_{0,n}(\sfi\sfz)=
%   \cali^{(1)}_{\qtan,n}[\sfi\sfz]^4 +
%   \rbr{\cali^{(2)}_{\qtan,n}+\qtanresid}[\sfi\sfz]^2
% \end{align*}
% with some functionals 
% $\cali^{(i)}_{\qtan,n}$ $(i=1,2)$ 
% of
% $O_M(n^{-H})$,
% such that 
% the equality 
% \begin{align*}
%   E\sbr{\Psi(\sfz)\mS^{(2,0)}_{0,n}(\sfi\sfz)}=
%   E\sbr{\Psi(\sfz)\bar\mS^{(2,0)}_{0,n}(\sfi\sfz)}
% \end{align*}
% holds.
% }

% {\myred Hence we can find a random symbol $\bar\mS^{(2,0)}_{0,n}(\sfi\sfz)$ such that 
% the equality 
% \begin{align*}
%   E\sbr{\Psi(\sfz)\mS^{(2,0)}_{0,n}(\sfi\sfz)}=
%   E\sbr{\Psi(\sfz)\bar\mS^{(2,0)}_{0,n}(\sfi\sfz)}
% \end{align*}
% holds 
% and there exist functionals 
% $\cali^{(i)}_{\qtan,n}$ $(i=1,2)$
% such that 
% $\cali^{(1)}_{\qtan,n}=O_M(n^{-H}),
% \cali^{(2)}_{\qtan,n}=O_M(n^{-H})$
% % +\hat O_M(n^{(\half-H)\vee(-\frac14)})$
% and 
% \begin{align*}
%   \bar\mS^{(2,0)}_{0,n}(\sfi\sfz)=
%   \cali^{(1)}_{\qtan,n}[\sfi\sfz]^4 +
%   \rbr{\cali^{(2)}_{\qtan,n}+\qtanresid}[\sfi\sfz]^2.
% \end{align*}}
Thus setting 
\begin{align*}
  \mS^{(2,0)}_{0}(\sfi\sfz)=0,
\end{align*}
we obtain the following lemma:
%%%%% $\qtan$のlimit %%%%%
%%%%% 240425-cal より %%%%%
\begin{lemma}
  There exists a sequence of random symbol $\bar\mS^{(2,0)}_{0,n}(\sfi\sfz)$
  such that 
  \begin{align*}
    E\sbr{\Psi(\sfz)\mS^{(2,0)}_{0,n}(\sfi\sfz)}&=
    E\sbr{\Psi(\sfz)\bar\mS^{(2,0)}_{0,n}(\sfi\sfz)}
  \end{align*}
  and 
  $\bar\mS^{(2,0)}_{0,n}(\sfi\sfz)\to\mS^{(2,0)}_{0}(\sfi\sfz)$ 
  in $L^p$ for any $p\geq1$.
\end{lemma}
In other words, the limit of the random symbol $\mS^{(2,0)}_{0,n}$
vanishes.

% \newpage
\subsubsection*{The Limit of $\mS^{(3,0)}_{n}$}
Recall that the random symbol $\mS^{(3,0)}_{n}$ is defined as 
\begin{align*}
  \mS^{(3,0)}_{n}(\sfi\sfz)=
  \frac13\qtor_n[\sfi\sfz]^3=
  \frac13 n^\half(D_{u_n})^2M_n[\sfi\sfz]^3.
\end{align*}
The decomposition of the functional $(D_{u_n})^2M_n$ is given at \eqref{eq:240610.1143},
and its leading term is of $O_M(n^{-\half})$.
Hence the coefficient satisfies 
$\frac13 n^\half(D_{u_n})^2M_n=O_M(1)$.
In fact, 
the following lemma shows that 
the functional $n^\half(D_{u_n})^2M_n$ is convergent in $L^p$:
\begin{lemma}\label{lemma:240611.2018}
  As $n\to\infty$, it holds that 
  \begin{align*}
    n^\half(D_{u_n})^2M_n\to
    \qtorconstTotal\times 
    \int^1_0 a(X_t)^3 dt
    % \qquad\text{in }L^p,
  \end{align*}
  in $L^p$ for any $p\geq1$, where 
  $\qtorconstTotal$ is defined at \eqref{eq:240611.2222}.
\end{lemma}
\noindent
The proof of this lemma is postponed to Section \ref{sec:240611.2231}.

Hence we set 
$\bar\mS^{(3,0)}_{n}(\sfi\sfz)=\mS^{(3,0)}_{n}(\sfi\sfz)$,
and define the random symbol $\mS^{(3,0)}(\sfi\sfz)$ by
\begin{align}\label{eq:240626.1954}
  \mS^{(3,0)}(\sfi\sfz)=
  \frac13 
  \qtorconstTotal\times 
  \int^1_0 a(X_t)^3 dt
  [\sfi\sfz]^3.
\end{align}

% \newpage 
\subsubsection*{The Limit of $\mS^{(1,0)}_{n}$}
% \subsubsection{Limit random symbol for $\mS^{(1,0)}_{n}$}
Recall that $N_n$ is written as \eqref{eq:240629.2411} with 
$N_n^{(1)}$,
$N_n^{(2)}$,
$N_n^{(3,\ell)}$ and
$N_n^{(4,\ell)}$
defined 
\eqref{eq:240629.2412},
\eqref{eq:240629.2413},
\eqref{eq:240629.2414} and 
\eqref{eq:240629.2415}, respectively.
% {\redmy Recall that we have defined 
% \begin{align*}
%   N_n&= 
%   N_n^{(1)} + N_n^{(2)}  
%   + \sum_{\ell=1}^k \cbr{N_n^{(3,\ell)} + N_n^{(4,\ell)}}
%   % \\&\quad
%   +\hat O_M(n^{(H-1)\vee(\half-H)})
% \end{align*}
% with
% \begin{align*}
%   N_n^{(1)}&=
%   % n\times\cnstVar{2k}{0} S_n^{(0;1;1)\dagger}=
%   2^{-1} \times\cnstVar{2k}{0}
%   \sum_{j=1}^{n-1}
%   a'(X_\tnj) V^{[1]}_\tnj \times 
%   \cbr{-B\brbr{\bbone^-_{2n,2j+1}} + B\brbr{\bbone^+_{2n,2j}}},
%   \\
%   N_n^{(2)}&=
%   % -n\times\cnstVar{2k}{0} S_n^{(\infty;1)\dagger}=
%   -2^{-1}\times\cnstVar{2k}{0}
%   \brbr{a(X_0)+a(X_1)}
%   \\
%   N_n^{(3,\ell)}&=
%   % n S_n^{(2,2;\ell,1)}=
%   n^{2H\ell} 
%   \sum_{j=1}^{n-1} 
%   A^{(2,2;\ell)}_{n,j}\times
%   I_{2\ell+1}((\diffker^{n}_{j})^{\otimes 2\ell-1}
%   \otimes(\bbone^n_j)^{\otimes2})
%   \\
%   N_n^{(4,\ell)}&=
%   % n S_n^{(2,2;\ell,2)}=
%   n^{2H\ell} 
%   \sum_{j=1}^{n-1} 
%   A^{(2,2;\ell)}_{n,j}\times
%   I_{2\ell+1}((\diffker^{n}_{j})^{\otimes 2\ell-1}
%   \otimes(\bbone^n_{j+1})^{\otimes2})
% \end{align*}}

First we calculate the limit of $N_n^{(4,\ell)}[\sfi\sfz]$.
Since $2\ell-1\geq1$, we can apply the IBP formula 
with respect to $\diffker^{n}_{j}$ as follows:
\begin{align*}
  &E\sbr{\Psi(\sfz)N_n^{(4,\ell)}[\sfi\sfz]}
  \\&=
  % E\sbr{\Psi(\sfz)\times
  % n^{2H\ell} 
  % \sum_{j=1}^{n-1} 
  % A^{(2,2;\ell)}_{n,j}\times
  % I_{2\ell+1}((\diffker^{n}_{j})^{\otimes 2\ell-1}
  % \otimes(\bbone^n_{j+1})^{\otimes2})
  % [\sfi\sfz]}
  % \\&=
  n^{2H\ell} 
  \sum_{j=1}^{n-1} 
  E\sbr{\Psi(\sfz)\times
  A^{(2,2;\ell)}_{n,j}\times
  I_{2\ell+1}((\diffker^{n}_{j})^{\otimes 2\ell-1}
  \otimes(\bbone^n_{j+1})^{\otimes2})}
  [\sfi\sfz]
  \\&=
  n^{2H\ell} 
  \sum_{j=1}^{n-1} 
  \half
  E\sbr{\Psi(\sfz)\abr{DG_\infty, \diffker^{n}_{j}}
  A^{(2,2;\ell)}_{n,j}\times
  I_{2\ell}({\diffker^{n}_{j}}^{\otimes 2\ell-2}
  \otimes{\bbone^n_{j+1}}^{\otimes2})
  }[\sfi\sfz]^3
  \\&\quad+
  n^{2H\ell} 
  \sum_{j=1}^{n-1} 
  E\sbr{\Psi(\sfz)
  \babr{D{A^{(2,2;\ell)}_{n,j}}, \diffker^{n}_{j}}\times
  I_{2\ell}({\diffker^{n}_{j}}^{\otimes 2\ell-2}
  \otimes{\bbone^n_{j+1}}^{\otimes2})
  }[\sfi\sfz]
  % \\&=
  % E\sbr{\Psi(\sfz)\times
  % \half 
  % n^{2H\ell-2H} 
  % \sum_{j=1}^{n-1} 
  % n^{2H} 
  % \abr{DG_\infty, \diffker^{n}_{j}}
  % A^{(2,2;\ell)}_{n,j}\times
  % I_{2\ell}({\diffker^{n}_{j}}^{\otimes 2\ell-2}
  % \otimes{\bbone^n_{j+1}}^{\otimes2})
  % [\sfi\sfz]^3}
  % \\&\quad+
  % E\sbr{\Psi(\sfz)\times
  % n^{2H\ell-2H} 
  % \sum_{j=1}^{n-1} 
  % n^{2H} 
  % \babr{D{A^{(2,2;\ell)}_{n,j}}, \diffker^{n}_{j}}\times
  % I_{2\ell}({\diffker^{n}_{j}}^{\otimes 2\ell-2}
  % \otimes{\bbone^n_{j+1}}^{\otimes2})
  % [\sfi\sfz]}
\end{align*}
Thus we define 
\begin{align*}
  \bar N_n^{(4,\ell)} (\sfi\sfz)&:=
  \bar N_n^{(4,\ell;1)} [\sfi\sfz]^3+\bar N_n^{(4,\ell;2)} [\sfi\sfz]
\end{align*}
with 
\begin{align*}
  \bar N_n^{(4,\ell;1)}&:=
  \half n^{2H\ell-2H} 
  \sum_{j=1}^{n-1} 
  n^{2H} 
  \abr{DG_\infty, \diffker^{n}_{j}}
  A^{(2,2;\ell)}_{n,j}\times
  I_{2\ell}({\diffker^{n}_{j}}^{\otimes 2\ell-2}
  \otimes{\bbone^n_{j+1}}^{\otimes2})
  \\
  \bar N_n^{(4,\ell;2)}&:=
  n^{2H\ell-2H} 
  \sum_{j=1}^{n-1} 
  n^{2H} 
  \babr{D{A^{(2,2;\ell)}_{n,j}}, \diffker^{n}_{j}}\times
  I_{2\ell}({\diffker^{n}_{j}}^{\otimes 2\ell-2}
  \otimes{\bbone^n_{j+1}}^{\otimes2}),
\end{align*}
then we have 
\begin{align*}
  &E\sbr{\Psi(\sfz)N_n^{(4,\ell)}[\sfi\sfz]}=
  E\sbr{\Psi(\sfz)\bar N_n^{(4,\ell)} (\sfi\sfz)}
\end{align*}

The exponent for the functional 
$ n^{-(2H\ell-2H)} \bar N_n^{(4,\ell;i)}$ $(i=1,2)$ is 
$1+(-1-H(2\ell-1)) = -H(2\ell-1)$,
which means that 
$\bar N_n^{(4,\ell;1)},\bar N_n^{(4,\ell;2)}=
n^{2H\ell-2H} \times O_M(n^{-H(2\ell-1)})=O_M(n^{-H})$.
Hence the weak limit of the random symbol $N_n^{(4,\ell)}[\sfi\sfz]$
is $0$.

By a similar argument, we can show that 
the weak limit of $N_n^{(3,\ell)}[\sfi\sfz]$ is also $0$.

With the following notation
\begin{align*}
  \weifuncNOne=
  2^{-1} \cnstVar{2k}{0}\times
  a'(X_\tnj) V^{[1]}_\tnj,\quad
  \kernelNOne=\bbone^+_{2n,2j}-\bbone^-_{2n,2j+1},
\end{align*}
the functional $N_n^{(1)}$ is written as
\begin{align*}
    N_n^{(1)}&=
    \sum_{j=1}^{n-1}
    \weifuncNOne\times 
    I_1(\kernelNOne).
    % \cbr{-B\brbr{\bbone^-_{2n,2j+1}} + B\brbr{\bbone^+_{2n,2j}}},
\end{align*}
Again by the IBP formula, we have
\begin{align*}
  &E\sbr{\Psi(\sfz)N_n^{(1)}[\sfi\sfz]}
  % \\&=
  % E\sbr{\Psi(\sfz)\times
  % \sum_{j=1}^{n-1}
  % \weifuncNOne\times 
  % I_1(\kernelNOne)
  % % \cbr{-B\brbr{\bbone^-_{2n,2j+1}} + B\brbr{\bbone^+_{2n,2j}}}
  % [\sfi\sfz]}
  % \\&=
  % \sum_{j=1}^{n-1}
  % E\sbr{\Psi(\sfz)\times
  % \weifuncNOne\times 
  % I_1(\kernelNOne)
  % % \cbr{-B\brbr{\bbone^-_{2n,2j+1}} + B\brbr{\bbone^+_{2n,2j}}}
  % }
  % [\sfi\sfz]
  \\&=
  \sum_{j=1}^{n-1}
  E\sbr{\Psi(\sfz)\times\half\babr{DG_\infty,\kernelNOne}
  \weifuncNOne}
  [\sfi\sfz]^3
  +% \\&\quad+
  \sum_{j=1}^{n-1}
  E\sbr{\Psi(\sfz)\times
  \babr{D\weifuncNOne, \kernelNOne}}
  [\sfi\sfz]
  % \\&=
  % E\sbr{\Psi(\sfz)
  % \times\bbcbr{
  % \Brbr{\sum_{j=1}^{n-1}
  % \half\abr{DG_\infty,\kernelNOne} \weifuncNOne}
  % [\sfi\sfz]^3 + 
  % \Brbr{\sum_{j=1}^{n-1}
  % \babr{D\weifuncNOne, \kernelNOne}}
  % [\sfi\sfz]}
  % }
\end{align*}
Define 
\begin{align*}
  \bar N_n^{(1)} (\sfi\sfz)&:=
  \bar N_n^{(1;1)} [\sfi\sfz]^3+\bar N_n^{(1;2)} [\sfi\sfz]
\end{align*}
with 
\begin{align*}
  \bar N_n^{(1;1)}&:=
  \sum_{j=1}^{n-1}
  \half\babr{DG_\infty,\kernelNOne} \weifuncNOne,\tand
  % \\
  \bar N_n^{(1;2)}:=
  \sum_{j=1}^{n-1}
  \babr{D\weifuncNOne, \kernelNOne},
\end{align*}
then we have 
\begin{align*}
  &E\sbr{\Psi(\sfz)N_n^{(1)}[\sfi\sfz]}=
  E\sbr{\Psi(\sfz)\bar N_n^{(1)} (\sfi\sfz)}.
\end{align*}

Since we can show that 
$\babr{DG_\infty,\kernelNOne}=O_M(n^{-2H})$ and 
$\babr{D\weifuncNOne, \kernelNOne}=O_M(n^{-2H})$,
we have 
\begin{align*}
  \bar N_n^{(1;i)}=O_M(n^{1-2H})\tfor i=1,2.
\end{align*}
Hence the weak limit of $N_n^{(1)}[\sfi\sfz]$ is $0$.

Since the weak limit of $N_n^{(2)}[\sfi\sfz]$ is obviously itself,
we set
\begin{align}\label{eq:240626.1955}
  \mS^{(1,0)}(\sfi\sfz)=N_n^{(2)}[\sfi\sfz]
  =-2^{-1}\times\cnstVar{2k}{0}
  \brbr{a(X_0)+a(X_1)}[\sfi\sfz],
\end{align}
and we obtain the following proposition:
\begin{proposition}
  There exists a sequence of random symbols $\bar\mS^{(1,0)}_{n}(\sfi\sfz)$
  such that 
  \begin{align*}
    E\sbr{\Psi(\sfz)\mS^{(1,0)}_{n}(\sfi\sfz)}&=
    E\sbr{\Psi(\sfz)\bar\mS^{(1,0)}_{n}(\sfi\sfz)}
    % \\
    % \bar\mS^{(1,0)}_{n}(\sfi\sfz)&\to 
    % \mS^{(1,0)}(\sfi\sfz)
    % % N_n^{(2)}[\sfi\sfz] 
    % \qquad\text{in }L^p\text{ for any }p\geq1.
  \end{align*}
  and 
  $\bar\mS^{(1,0)}_{n}(\sfi\sfz)\to 
  \mS^{(1,0)}(\sfi\sfz)$
  % N_n^{(2)}[\sfi\sfz] 
  in $L^p$ for any $p\geq1$.
  % where 
  % \begin{align*}
  %   N_n^{(2)}[\sfi\sfz] &=  
  %   - \half\cnstVar{2k}{0} 
  %   \cbr{f_1 (V^{[1]}_1)^{2k} + f_0 (V^{[1]}_0)^{2k}}[\sfi\sfz].
  % \end{align*}
  % {\rb $N_n^{(2)}$の記号は変えた方がいいと思われる}
  % {\rb ``any $p\geq1$''でよいか？}
\end{proposition}

\subsubsection*{The Limit of $\mS^{(2,0)}_{1,n}$}
Since it holds that
$D_{u_n} N_n=\hat O(n^{(H-1)\vee(\half-H)})$
by Proposition \ref{prop:240610.2337},
we define $\mS^{(2,0)}_{1}$ by
\begin{align*}
  \mS^{(2,0)}_{1}(\sfi\sfz)=0.
\end{align*}
% 
% 

% \newpage

% \newpage
\subsection{Proof of Theorem \ref{thm:240620.2214}}
\label{sec:240630.2339}
\subsubsection{Condition {\bf[D]} (ii)}
% \subsubsection{Order estimates of related functionals}
\label{sec:240628.1902}
From
Proposition \ref{prop:240607.1051}, %$u_n$
Lemma \ref{lemma:240523.2220} and  %$\abr{DM_n,u_n}-G_\infty$
Proposition \ref{prop:240522.2303}, %$Z_n$の確率展開, $N_n$
the conditions 
\eqref{220215.1241}, 
\eqref{220215.1242} and 
\eqref{220215.1246} of
{\bf[D]} (ii) 
are verified in our context.
% 
% By Lemma \ref{lemma:240523.2220}, we have 
% \begin{align*}
%   G_n^{(2)}=\abr{DM_n,u_n}-G_\infty&=O_M(n^{-\half}).
% \end{align*}
We check the remainder as follows.

\subsubsection*{{\bf[D]} (ii) \eqref{220215.1243} and \eqref{220215.1244}}

% \begin{proposition}
%   % \item Both the functionals $G_n^{(3)}$ and $D_{u_n}G_n^{(3)}$ are of the order of $O_M(n^{-H})$.
%   The functional $G_n^{(3)}$ is estimated as $G_n^{(3)}=O_M(n^{-H})$.

%   As a corollary, it also holds that $D_{u_n}G_n^{(3)}=O_M(n^{-H})$.
% \end{proposition}

The functional $G_n^{(3)}=D_{u_n}G_\infty$ can be written as 
\begin{align*}
  G_n^{(3)}=
  D_{u_n}G_\infty&=
  \abr{DG_\infty, u_n}
  % \\&=
  % \sum_{\ell=1}^k \abr{DG_\infty, u_n^{(\ell)}}
  =% \\&=
  \sum_{\ell=1}^k \unlconst{\ell} \abr{DG_\infty, u_n^{(\ell)\dagger}}
\end{align*}
with 
\begin{align*}
  \abr{DG_\infty, u_n^{(\ell)\dagger}}&=
  % \abr{DG_\infty, 
  % n^{2H\ell-\half}
  % \sum_{j=1}^{n-1} 
  % a(X_\tnj)
  % I_{2\ell-1}({\diffker^{n}_{j}}^{\otimes2\ell-1})
  % \diffker^{n}_{j}}
  % \\&=
  % n^{2H\ell-\half}
  % \sum_{j=1}^{n-1} 
  % \abr{DG_\infty, \diffker^{n}_{j}}
  % a(X_\tnj)
  % I_{2\ell-1}({\diffker^{n}_{j}}^{\otimes2\ell-1})
  % \\&=
  n^{2H(\ell-1)-\half}
  \sum_{j=1}^{n-1} 
  n^{2H}
  \abr{DG_\infty, \diffker^{n}_{j}}
  a(X_\tnj)
  I_{2\ell-1}({\diffker^{n}_{j}}^{\otimes2\ell-1}).
\end{align*}
%%%%%%%%%%%%%%%%%%%%%%%%%%%%%%%%%%%%
%%%%% [CommentOut:FurtherTodo] %%%%%
% {\rb[この書き方はcoolか？]
%%%%%%%%%%%%%%%%%%%%%%%%%%%%%%%%%%%%
The exponent for the factor 
$n^{-(2H(\ell-1)-\half)}\babr{DG_\infty, u_n^{(\ell)\dagger}}$ is
$1+(-\half-H(2\ell-1))=\half-H(2\ell-1)$.
Hence we have 
\begin{align*}
  \abr{DG_\infty, u_n^{(\ell)\dagger}}=
  n^{2H(\ell-1)-\half}\times O_M(n^{\half-H(2\ell-1)})=
  O_M(n^{-H})
\end{align*}
for any $\ell=1,...,k$,
and 
\begin{align}
  G_n^{(3)}=D_{u_n}G_\infty&=O_M(n^{-H}),
  \label{eq:240610.1039}
\end{align}
which confirms the condition \eqref{220215.1243}. % since $r_n=n^{-\half}$.

By Proposition \ref{prop:240607.1051} and \eqref{eq:240610.1039},
we have
\begin{align}
  D_{u_n}G_n^{(3)}=
  (D_{u_n})^2G_\infty=
  O_M(n^{-H}),
  \label{eq:240610.1110}
\end{align}
and the condition \eqref{220215.1244}.

\begin{remark}
  Writing $G_n^{(3)}$ explicitly and applying the argument of exponent
  (in particular Proposition \ref{prop:240617.2009}),
  we can obtain a sharper estimate for $D_{u_n}G_n^{(3)}$,
  namely 
  $D_{u_n}G_n^{(3)}=O_M(n^{-H-\half})$.
  However we do not need the estimate of
  that order for the proof of the asymptotic expansion.
\end{remark}

\subsubsection*{{\bf[D]} (ii) \eqref{220215.1245}}
% \subsubsection{Order estimate related to $G_n^{(2)}$}
% \begin{proposition}
%   \item [(i)] The functional $G_n^{(2)}$ is estimated as 
%   \begin{align*}
%     G_n^{(2)}=O_M(n^{-\half}).
%   \end{align*}
%   \item [(ii)] The functional $(D_{u_n})^2G_n^{(2)}$ is estimated as 
%   \begin{align*}
%     (D_{u_n})^2G_n^{(2)}=
%     O_M(n^{-H}).
%   \end{align*}
% \end{proposition}

The functional $(D_{u_n})^2G_n^{(2)}$ is written as 
\begin{align*}
  (D_{u_n})^2G_n^{(2)}=
  (D_{u_n})^3M_n-(D_{u_n})^2G_\infty.
\end{align*}
By \eqref{eq:240610.1110}, we have 
$(D_{u_n})^2G_\infty=O_M(n^{-H})$.
For $(D_{u_n})^3M_n$,
recall that we have the decomposition \eqref{eq:240610.1143} of $(D_{u_n})^2M_n$.
In that decomposition,
the principal term is the linear combination of the functionals  
$\caliqtor$ (defined at \eqref{eq:240610.1145}).
When we write $\graphQtorSharp$ for
the weighted graph representing the functional $\caliqtor$,
it is obvious that $\barq(\graphQtorSharp)=0$.
Hence Proposition \ref{prop:240617.2009} (i) applies to 
$D_{u_n^{(\ell_4)}}\caliqtor$ for $\ell_4=1,...,k$,
and we have the estimate
\begin{align*}
  D_{u_n^{(\ell_4)}}\caliqtor=O_M(n^{-\half-H}),
\end{align*}
since we have checked that 
$\caliqtor=O_M(n^{-\half})$ in Section \ref{sec:240617.2018}.
Hence we have 
\begin{align*}
  (D_{u_n})^3M_n&=
  \sum_{(\ell_1,\ell_2,\ell_3,m)\in\Lambda^{(k;3)}_{\sharp}}
  \qtorconst\times D_{u_n}\caliqtor
  +O_M(n^{-H})
  % \\&
  % =O_M(n^{-1})+O_M(n^{-H})
  =O_M(n^{-H}).
\end{align*}

\subsubsection*{{\bf[D]} (ii) \eqref{220215.1247}}
By Proposition \ref{prop:240610.2337},
we have 
$D_{u_n}N_n=\hat O_M(n^{(H-1)\vee(\half-H)})$.
Thus the functional $(D_{u_n})^2N_n$ is estimated as 
$(D_{u_n})^2N_n=\hat O_M(n^{(H-1)\vee(\half-H)})$
thanks to $u_n=O_M(1)$
(Proposition \ref{prop:240607.1051}).

% \newpage

% \newpage
\subsubsection{Condition {\bf[D]} (iv)}
% \subsubsection{Nondegeneracy}
\label{sec:240628.1903}
% \subsection{Condition for nondegeneracy}
%%%%%  %%%%%
%%%%% 240426-cal より %%%%%
% We set {\rb define(?)} $s_n=s_\infty$ as 
% \begin{align*}
%   s_n=s_\infty:=
%   \half C_\Delta\times 
%   % G_\infty^{\dagger}
%   \int_0^1 a(X_t)^2dt
% \end{align*}
% with 
% \begin{align*}
%   C_{\Delta}&:=
%   \Bcbr{\sum_{\ell_1=1}^k 
%   2\ell_1\times   C^{(\ell_1)} \times C^{\widehat\rho}_{\ell_1}}.
% \end{align*}
% {\rb Then we have the following estimate, which 
% ensures the nondegeneracy of the distribution of $M_n$.[書き直す]}
% 
Since Assumption \ref{ass:240626.1518} (iii) validates 
Condition {\bf[D]} (iv) (a),
we will verify {\bf[D]} (iv) (b).

% \begin{proposition}\label{prop:240523.1206}
%   There exists some positive 
%   constant $C_\eqref{prop:240523.1206}>0$ such that
%   % random variable $s_\infty\in\bbD^\infty$ such that
%   \begin{align*}
%     P\sbr{\abr{DM_n,DM_n}\leq 
%     C_\eqref{prop:240523.1206} \int_0^1 a(X_t)^2dt
%     }&=O(n^{-\frac{p}2})
%   \end{align*}
%   for any $p\geq1$.
% \end{proposition}
% \begin{proof}
  Since we have 
  \begin{align*}
    M_n = \sum_{\ell=1}^k 
    M_n^{(\ell)},
    % \delta(u_n^{(\ell)}),
    \tand
    M_n^{(\ell)}=
    \delta(u_n^{(\ell)})=
    % n^{\half}S_n^{(1,\ell)}-\rsdOne{\ell}=
    M_n^{\prime(\ell)}-\rsdOne{\ell}
  \end{align*}
  with $\rsdOne{\ell}=O_M(n^{-H})$
  (see Lemma \ref{lemma:240523.1805} and \eqref{eq:240612.0901}),
  we can write
  \begin{align*}
    \abr{DM_n,DM_n}&=
    \abr{\sum_{\ell=1}^k D\rbr{M_n^{\prime(\ell)} - \rsdOne{\ell}},
    \sum_{\ell=1}^k D\rbr{M_n^{\prime(\ell)} - \rsdOne{\ell}}}
    \\&=
    \abr{\sum_{\ell=1}^k DM_n^{\prime(\ell)},
    \sum_{\ell=1}^k DM_n^{\prime(\ell)}}
    +O_M(n^{-H}).
    % \\&=
    % \sum_{\ell_1,\ell_2=1}^k
    % \abr{DM_n^{\prime(\ell_1)}, DM_n^{\prime(\ell_2)}}
    % +O_M(n^{-H})
  \end{align*}

By the definition of $M_n^{\prime(\ell)}$, we have
\begin{align*}
  DM_n^{\prime(\ell)}&=
  % C^{2k}_{2\ell}\times n^{2\ell H-\half}
  % \sum_{j\in[n-1]} D\cbr{a(X_{\tnj}) \times 
  % I_{2\ell} \nrbr{{\diffker^n_j}^{\otimes 2\ell}}}
  % \\&=
  % C^{2k}_{2\ell}\times n^{2\ell H-\half}
  % \sum_{j\in[n-1]} (Da(X_{\tnj})) \times 
  % I_{2\ell} \nrbr{{\diffker^n_j}^{\otimes 2\ell}}
  % \\&\quad+
  % 2\ell\times C^{2k}_{2\ell}\times n^{2\ell H-\half}
  % \sum_{j\in[n-1]} a(X_{\tnj}) \times 
  % I_{2\ell-1} \nrbr{{\diffker^n_j}^{\otimes 2\ell-1}} \diffker^n_j
  % \\&=
  2\ell\times u_n^{(\ell)} + R_n^{DM';(\ell)}
\end{align*}
with % where we define 
\begin{align*}
  R_n^{DM';(\ell)} &=
  \unlconst{\ell}\times
  % C^{2k}_{2\ell}\times 
  n^{2\ell H-\half}
  \sum_{j\in[n-1]} (Da(X_{\tnj})) \times 
  I_{2\ell} \nrbr{{\diffker^n_j}^{\otimes 2\ell}}.
\end{align*}
Writing
$R_n^{DM'}:=\sum_{\ell=1}^k R_n^{DM';(\ell)}$,
% \begin{align*}
%   \sum_{\ell=1}^k DM_n^{\prime(\ell)}&=
%   % % \sum_{\ell=1}^k \cbr{2\ell\times u_n^{(\ell)} + R_n^{DM';(\ell)}}
%   % % \\&=
%   % \sum_{\ell=1}^k 2\ell\times u_n^{(\ell)}+
%   % \sum_{\ell=1}^k R_n^{DM';(\ell)}
%   % =
%   \sum_{\ell=1}^k 2\ell\times u_n^{(\ell)}+R_n^{DM'}
% \end{align*}
we decompose $\abr{DM_n,DM_n}$ as follows:
\begin{align}
  &\abr{DM_n,DM_n}
  \nn\\&=
  \abr{\sum_{\ell=1}^k DM_n^{\prime(\ell)},
  \sum_{\ell=1}^k DM_n^{\prime(\ell)}}
  +O_M(n^{-H})
  \nn\\&=
  \abr{\sum_{\ell=1}^k 2\ell\times u_n^{(\ell)}+R_n^{DM'},
  \sum_{\ell=1}^k 2\ell\times u_n^{(\ell)}+R_n^{DM'}}
  +O_M(n^{-H})
  \nn\\&=
  % \abr{\sum_{\ell=1}^k 2\ell\times u_n^{(\ell)},
  % \sum_{\ell=1}^k 2\ell\times u_n^{(\ell)}}+
  % 2\abr{\sum_{\ell=1}^k 2\ell\times u_n^{(\ell)},R_n^{DM'}}+
  % \abr{R_n^{DM'},R_n^{DM'}}
  % +O_M(n^{-H})
  % \\&=
  \abr{\sum_{\ell=1}^k 2\ell\times u_n^{(\ell)},
  \sum_{\ell=1}^k 2\ell\times u_n^{(\ell)}}+
  2\mR^{(1)}_n+
  \mR^{(2)}_n
  +O_M(n^{-H}),
  \label{eq:240612.0941}
\end{align}
where we define
\begin{align*}
  \mR^{(1)}_n&:=
  \abr{R_n^{DM'},\sum_{\ell=1}^k 2\ell\times u_n^{(\ell)}},\tand
  % \\
  \mR^{(2)}_n:=
  \abr{R_n^{DM'},R_n^{DM'}}\geq0.
\end{align*}
We have
\begin{align*}
  \mR^{(1)}_n&=
  % \abr{\sum_{\ell=1}^k 2\ell\times u_n^{(\ell)},R_n^{DM'}}
  % =% \\&=
  % \abr{\sum_{\ell=1}^k 2\ell\times u_n^{(\ell)},
  % \sum_{\ell=1}^k R_n^{DM';(\ell)}}
  % \\&=
  \sum_{\ell_1,\ell_2=1}^k
  2\ell_2\abr{R_n^{DM';(\ell_1)}, u_n^{(\ell_2)}}
  =
  \sum_{\ell_1,\ell_2=1}^k
  2\ell_2
  G_n^{(\ell_1,\ell_2;2)}=O_M(n^{-H}),
\end{align*}
where $G_n^{(\ell_1,\ell_2;2)}$ is defined at \eqref{eq:240403.1754} and 
estimated at \eqref{eq:240607.1235}.

The principal term of \eqref{eq:240612.0941} is decomposed as 
\begin{align*}
  &\abr{\sum_{\ell=1}^k 2\ell\times u_n^{(\ell)},
  \sum_{\ell=1}^k 2\ell\times u_n^{(\ell)}}
  \\&=
  \sum_{\ell_1=1}^k\sum_{\ell_2=1}^k
  2\ell_1\times 2\ell_2\times
  \abr{ u_n^{(\ell_1)},u_n^{(\ell_2)}}
  \\&=
  \sum_{\ell_1=1}^k\sum_{\ell_2=1}^k
  2\ell_2\times G_n^{(\ell_1,\ell_2;1)}
  \\&=
  \sum_{\ell_1=1}^k\sum_{\ell_2=1}^k
  2\ell_2\times 
  \sum_{m=0}^{(2\ell_1-1)\wedge(2\ell_2-1)}
  \gtwoconst{\ell_1,\ell_2,m}\times
  G_n^{(\ell_1,\ell_2;m)\dagger\dagger}
  \\&=
  \sum_{\ell_1=1}^k
  2\ell_1\times 
  \gtwoconst{\ell_1,\ell_1,2\ell_1-1}\times
  G_n^{(\ell_1,\ell_1;2\ell_1-1)\dagger\dagger}
  +O_M(n^{-\half})
  \\&=
  \sum_{\ell_1=1}^k
  2\ell_1\times 
  \gtwoconst{\ell_1,\ell_1,2\ell_1-1}\times
  \qtrhoconst{\ell_1}\times
  \int_0^1 a(X_t)^2dt
  % G_n^{(\ell_1,\ell_1;2\ell_1-1)\dagger\dagger}
  +O_M(n^{-\half}),
\end{align*}
where we used the definition \eqref{eq:240403.1754} of $G_n^{(\ell_1,\ell_2;1)}$,
the relations \eqref{eq:240607.1237}, \eqref{eq:240607.1238} and
\eqref{eq:240607.1332},
the estimates in Lemma \ref{lemma:240523.2129}, and 
the convergence proved in Lemma \ref{lemma:240611.2313}.

Setting 
\begin{align*}
  C_{\Delta}&:=
  \sum_{\ell_1=1}^k
  2\ell_1\times 
  \gtwoconst{\ell_1,\ell_1,2\ell_1-1}\times
  \qtrhoconst{\ell_1},\tand 
  G_\infty^{\dagger}:=
  \int_0^1 a(X_t)^2dt
\end{align*}
for short, 
we have
\begin{align*}
  &\abr{DM_n,DM_n}
  =% \nn\\&=
  % \abr{\sum_{\ell=1}^k DM_n^{\prime(\ell)},
  % \sum_{\ell=1}^k DM_n^{\prime(\ell)}}
  % +O_M(n^{-H})
  % \nn\\&=
  % \abr{\sum_{\ell=1}^k 2\ell\times u_n^{(\ell)}+R_n^{DM'},
  % \sum_{\ell=1}^k 2\ell\times u_n^{(\ell)}+R_n^{DM'}}
  % +O_M(n^{-H})
  % \nn\\&=
  % \abr{\sum_{\ell=1}^k 2\ell\times u_n^{(\ell)},
  % \sum_{\ell=1}^k 2\ell\times u_n^{(\ell)}}+
  % 2\mR^{(1)}_n+
  % \mR^{(2)}_n
  % +O_M(n^{-H})
  % \\&=
  C_{\Delta}\times
  G_\infty^{\dagger}
  % \int_0^1 a(X_t)^2dt
  +\mR^{(2)}_n
  +O_M(n^{-\half}).
\end{align*}
Note that $C_{\Delta}>0$.

Denote the above term of $O_M(n^{-\half})$ by $\mR^{(0)}_n$, and 
let 
\begin{align*}
  % s_n=
  s_\infty:=
  \half C_\Delta\times G_\infty^{\dagger}.
\end{align*}
% {\redmy that is 
% $C_\eqref{prop:240523.1206}:=\half C_\Delta$.}
Note that we have 
% $G_\infty$ is defined at \eqref{eq:240630.1538} as
$G_\infty = C_{G_\infty} G_\infty^{\dagger}$ 
with some positive constant $C_{G_\infty}$ (defined at \eqref{eq:240628.1937}),
% $G_\infty = C_{G_\infty} \int_0^1 a(X_t)^2dt,$
and recall that we have assumed 
$(G_\infty)^{-1}\in L^{\infty-}=\cap_{p>1}L^p$
in Assumption \ref{ass:240626.1518} (iii).
Then we have
\begin{align*}
  P\sbr{\abr{DM_n,DM_n}\leq s_\infty}%&=
  % P\sbr{C_{\Delta}\times G_\infty^{\dagger}+\mR^{(2)}_n+\mR^{(0)}_n
  % \leq s_\infty}
  &\leq% \\&\leq
  % P\sbr{C_{\Delta}\times G_\infty^{\dagger}+\mR^{(0)}_n\leq s_\infty}
  % =% \\&=
  % P\sbr{2s_\infty+\mR^{(0)}_n\leqs_\infty}=
  P\sbr{\mR^{(0)}_n\leq-s_\infty}
  \hspace{100pt}
  (\because \mR^{(2)}_n\geq0)
  \\&\leq 
  % P\sbr{\abs{\mR^{(0)}_n}\geq s_\infty}=
  P\sbr{\babs{\mR^{(0)}_n}(s_\infty)^{-1}\geq1}
  \\&\leq
  \norm{\babs{\mR^{(0)}_n}(s_\infty)^{-1}}_p^p
  \\&\leq %\leq% 
  \cbr{\bnorm{\mR^{(0)}_n}_{2p}\times 
  \norm{s_\infty^{-1}}_{2p}}^p
  =O(n^{-\frac{p}2})
\end{align*}
for any $p\geq1$.
Thus,
Condition {\bf[D]} (iv) (b) holds with $s_n:=s_\infty$.

% \end{proof}

\section{Technical lemmas}\label{sec:240626.1842}
% \subfile{7-lemmas.tex}

In the subsequent sections, we use the
following lemma, which is %taken from \cite{2024Yamagishi-asymptotic}, and 
a collection of basic properties of the solution of 
SDE \eqref{eq:240626.1939}.
\begin{lemma}[Lemma 6.1 of \cite{2024Yamagishi-asymptotic}]%\label{230924.1800}
  \label{lemma:240617.2154}
  Suppose that $X_t$ is the solution of SDE \eqref{eq:240626.1939}
  under Assumption \ref{ass:240626.1518} (i).
  Let $T>0$ and
  $f\in C^\infty(\bbR)$ with bounded derivatives of any order.
  \item[(i)]
  Let $r_0<r_1\in[0,T]$.
  Writing $Y^{(0)}_t=t$ and $Y^{(1)}_t=B_t$,
  consider the $k$-th $(k\in\bbN)$ iterated integral 
  \begin{align*}
    F_{r_0,r_1}^{(i)}&=
    \int^{r_1}_{r_0}
    \int^{t^{(1)}}_{r_0}\cdots\rbr{
    \int^{t^{(k-1)}}_{r_0}
    f(X_{t^{(k)}})dY^{(i_k)}_{t^{(k)}}}
    \cdots dY^{(i_2)}_{t^{(2)}}dY^{(i_1)}_{t^{(1)}}
  \end{align*}
  for $(i_j)_{j=1}^k\in\cbr{0,1}^k$.
  Then for any $N\in\bbZ_{\geq0}$ and $p\geq1$,
  the following estimate holds with any $\epsilon>0$:
  \begin{align*}
    \norm{F_{r_0,r_1}^{(i)}}_{N,p}&\leq
    C \abs{r_0-r_1}^{k_0+k_1H-\epsilon\bbone_\cbr{k_1>0}},
  \end{align*}
  where 
  $k_\ell=\sharp\cbr{1\leq j\leq k \mid i_j=\ell}$ for $\ell=0,1$, and
  the constant $C$ is independent of $r_0$ and $r_1$.

  \item[(ii)] %\ref{230807.1543}
  Let $r_0<r_1\in[0,T]$, $N\in\bbZ_{\geq0}$ and $p>1$.
  For any $\epsilon>0$,
  it holds that 
  \begin{align*}
    \norm{f(X_{r_1})-f(X_{r_0})}_{N,p} \leq C \abs{r_1-r_0}^{H-\epsilon}
  \end{align*}
  with the constant $C$ independent of $r_0$ and $r_1$.

  \item[(iii)]
  The difference of the integral of $f(X_t)$ %with respect to $dt$ 
  and its Riemann sum is estimated as
  \begin{align}\label{230924.1812}
    \int^T_0 f(X_t)dt - 
    \frac{T}{n}\sum_{j=0}^{n-1} f(X_{t_{j}})
    =O_M(n^{-1}),
  \end{align}
  where $t_j=jT/n$.
\end{lemma}

\subsection{Lemmas related to stochastic expansion}
\label{sec:240626.1850}
\subsubsection{Decomposition of the second order difference of fSDE}
\label{sec:240606.1759}
First we decompose the difference of the process $X$ satisfying SDE \eqref{eq:240626.1939}.
We write 
\begin{align*}
  \diff{n}{j}X:=X_\tnj-X_\tnjm
\end{align*}
for $j=1,...,n$,
where $\tnj=j/n$.
Recall that $V^{[i;k]}$ is the $k$-th derivative of $V^{[i]}:\bbR\to\bbR$, and
% {\rb $V^{[i]}$の微分のnotationについて触れたい．}
we write 
$V^{[i]}_t=V^{[i]}(X_t)$ %for $i=1,2$.
and $V^{[i;k]}_t=V^{[i;k]}(X_t)$.

The following lemma shows two different ways of the decomposition of
$\diff{n}{j}X$.
Approximating the integrals by the value at the left and
right ends of the interval,
we obtain the decompositions \eqref{eq:240605.1541} and 
\eqref{eq:240605.1542}, respectively.
\begin{lemma}\label{lemma:240523.1445}
  For $n\geq2, j=1,...,n-1$, the difference of $X$ can be written as follows:
  \begin{align}
    \diff{n}{j}X &= 
    V^{[1]}_\tnj \int^\tnj_\tnjm  dB_t + 
    V^{[2]}_\tnj \int^\tnj_\tnjm  dt 
    \nn\\&\quad+
    V^{[1,1]}_\tnj V^{[1]}_\tnj \int^\tnj_\tnjm dB_t \int^t_\tnj dB_s 
    \nn\\&\quad+
    V^{[1,1]}_\tnj V^{[2]}_\tnj \int^\tnj_\tnjm dB_t \int^t_\tnj ds +
    V^{[2,1]}_\tnj V^{[1]}_\tnj \int^\tnj_\tnjm dt \int^t_\tnj dB_s 
    \nn\\&\quad+ 
    \hat O(n^{-2}\vee n^{-3H})
    \label{eq:240605.1541}
    \\
    \diff{n}{j+1}X &= 
    V^{[1]}_\tnj \int^\tnjp_\tnj  dB_t + 
    V^{[2]}_\tnj \int^\tnjp_\tnj  dt 
    \nn\\&\quad+
    V^{[1,1]}_\tnj V^{[1]}_\tnj \int^\tnjp_\tnj dB_t \int^t_\tnj dB_s 
    \nn\\&\quad+
    V^{[1,1]}_\tnj V^{[2]}_\tnj \int^\tnjp_\tnj dB_t \int^t_\tnj ds +
    V^{[2,1]}_\tnj V^{[1]}_\tnj \int^\tnjp_\tnj dt \int^t_\tnj dB_s 
    \nn\\&\quad+ 
    \hat O(n^{-2}\vee n^{-3H})
    \label{eq:240605.1542}
  \end{align}
\end{lemma}
\begin{proof}
  For the expansion \eqref{eq:240605.1541},
  we approximate the functions in the integrals by the values at 
  the right end of the interval (i.e. $\tnj$):
  \begin{align*}
    &\diff{n}{j}X 
    % X_\tnj - X_\tnjm
    % \\&=
    % \int^\tnj_\tnjm dX_t
    \\&=
    \int^\tnj_\tnjm V^{[1]}_t dB_t+
    \int^\tnj_\tnjm V^{[2]}_t dt
    \\&=
    V^{[1]}_\tnj \int^\tnj_\tnjm  dB_t + 
    \int^\tnj_\tnjm \rbr{V^{[1]}_t-V^{[1]}_\tnj} dB_t 
    +% \\&\quad+
    V^{[2]}_\tnj \int^\tnj_\tnjm  dt +
    \int^\tnj_\tnjm \rbr{V^{[2]}_t-V^{[2]}_\tnj} dt 
    % \\&=
    % V^{[1]}_\tnj \int^\tnj_\tnjm  dB_t + 
    % V^{[2]}_\tnj \int^\tnj_\tnjm  dt 
    % \\&\quad+
    % \int^\tnj_\tnjm dB_t %\rbr{V^{[1]}(X_t)-V^{[1]}(X_\tnj)} 
    % \int^t_\tnj V^{[1,1]}_s dX_s
    % +% \\&\quad+
    % \int^\tnj_\tnjm dt %\rbr{V^{[2]}(X_t) - V^{[2]}(X_\tnj)}
    % \int^t_\tnj V^{[2,1]}_s dX_s
    \\&=
    V^{[1]}_\tnj \int^\tnj_\tnjm  dB_t + 
    V^{[2]}_\tnj \int^\tnj_\tnjm  dt 
    \\&\quad+
    \int^\tnj_\tnjm dB_t \int^t_\tnj V^{[1,1]}_s V^{[1]}_s dB_s+ 
    \int^\tnj_\tnjm dB_t \int^t_\tnj V^{[1,1]}_s V^{[2]}_s ds
    \\&\quad+
    \int^\tnj_\tnjm dt \int^t_\tnj V^{[2,1]}_s V^{[1]}_s dB_s + 
    \int^\tnj_\tnjm dt \int^t_\tnj V^{[2,1]}_s V^{[2]}_s ds
    \\&=
    V^{[1]}_\tnj \int^\tnj_\tnjm  dB_t + 
    V^{[2]}_\tnj \int^\tnj_\tnjm  dt 
    \\&\quad+
    V^{[1,1]}_\tnj V^{[1]}_\tnj \int^\tnj_\tnjm dB_t \int^t_\tnj dB_s 
    \\&\quad+
    V^{[1,1]}_\tnj V^{[2]}_\tnj \int^\tnj_\tnjm dB_t \int^t_\tnj ds +
    V^{[2,1]}_\tnj V^{[1]}_\tnj \int^\tnj_\tnjm dt \int^t_\tnj dB_s 
    \\&\quad+ 
    \hat O(n^{-2}\vee n^{-3H})
  \end{align*}
  The estimate of the residual terms (i.e. $\hat O(n^{-2}\vee n^{-3H})$)
  follows from Lemma \ref{lemma:240617.2154} (i).
  
  For $\diff{n}{j+1}X$, we can prove the expansion \eqref{eq:240605.1542}
  in a similar manner,
  but with approximation by the values at the left end.
\end{proof}

\begin{proof}[Proof of Lemma \ref{lemma:240523.1438}]
  Since we have defined 
  \begin{align*}
    \bbone^n_j=\bbone_{[\tnjm,\tnj]},\quad
    \diffker^n_j=\bbone^n_{j+1}-\bbone^n_{j},
  \end{align*}
  we have 
  \begin{align*}
    B(\bbone^n_j)=B_\tnj-B_\tnjm,\quad
    B(\diffker^n_j)=
    (B_\tnjp-B_\tnj)-(B_\tnj-B_\tnjm).
  \end{align*}
  Recall that we have introduced 
  % We introduce 
  the following notation for short:
  \begin{align*}
    \fvoo_t = \dfvo_t\fvo_t,\quad
    \fvot_t = \dfvo_t\fvt_t,\quad
    \fvto_t = \dfvt_t\fvo_t.%,\quad
  \end{align*}
  By Lemma \ref{lemma:240523.1445},
  the second order difference of $X$ %the solution $X$ to SDE \ref{eq:240626.1939}
  decomposes as follows:
  \begin{align*}
    \secDiff{n}{j}X &= 
    \diff{n}{j+1}X-\diff{n}{j}X
    % (X_\tnjp - X_\tnj) - (X_\tnj - X_\tnjm)
    \\&=
    \Bigg\{
    \fvo_\tnj \int^\tnjp_\tnj  dB_t + 
    \fvt_\tnj \int^\tnjp_\tnj  dt
    +% \\&\quad+
    \fvoo_\tnj \int^\tnjp_\tnj dB_t \int^t_\tnj dB_s 
    % \dfvo_\tnj \fvo_\tnj \int^\tnjp_\tnj dB_t \int^t_\tnj dB_s 
    \\&\quad+
    \fvot_\tnj \int^\tnjp_\tnj dB_t \int^t_\tnj ds +
    % \dfvo_\tnj \fvt_\tnj \int^\tnjp_\tnj dB_t \int^t_\tnj ds +
    \fvto_\tnj \int^\tnjp_\tnj dt \int^t_\tnj dB_s\Bigg\}
    % \dfvt_\tnj \fvo_\tnj \int^\tnjp_\tnj dt \int^t_\tnj dB_s\Bigg\}
    \\&\quad-\Bigg\{
    \fvo_\tnj \int^\tnj_\tnjm  dB_t + 
    \fvt_\tnj \int^\tnj_\tnjm  dt
    +% \\&\quad+
    \fvoo_\tnj \int^\tnj_\tnjm dB_t \int^t_\tnj dB_s 
    % \dfvo_\tnj \fvo_\tnj \int^\tnj_\tnjm dB_t \int^t_\tnj dB_s 
    \\&\quad+
    \fvot_\tnj \int^\tnj_\tnjm dB_t \int^t_\tnj ds +
    % \dfvo_\tnj \fvt_\tnj \int^\tnj_\tnjm dB_t \int^t_\tnj ds +
    \fvto_\tnj \int^\tnj_\tnjm dt \int^t_\tnj dB_s \Bigg\}
    % \dfvt_\tnj \fvo_\tnj \int^\tnj_\tnjm dt \int^t_\tnj dB_s \Bigg\}
    \\&\quad+ 
    \hat O(n^{-2}\vee n^{-3H})
    \\&=
    \fvo_\tnj \times\rbr{\int^\tnjp_\tnj  dB_t - \int^\tnj_\tnjm  dB_t}+
    \fvt_\tnj \times\rbr{\int^\tnjp_\tnj  dt -  \int^\tnj_\tnjm  dt}
    \\&\quad+
    \fvoo_\tnj \times\rbr{
      \int^\tnjp_\tnj dB_t \int^t_\tnj dB_s - 
      \int^\tnj_\tnjm dB_t \int^t_\tnj dB_s}
    \\&\quad+
    \fvot_\tnj \times\rbr{
      \int^\tnjp_\tnj dB_t \int^t_\tnj ds - 
      \int^\tnj_\tnjm dB_t \int^t_\tnj ds}
    \\&\quad+
    \fvto_\tnj \times\rbr{
      \int^\tnjp_\tnj dt \int^t_\tnj dB_s - 
      \int^\tnj_\tnjm dt \int^t_\tnj dB_s}
    \\&\quad+
    \hat O(n^{-2}\vee n^{-3H})
  \end{align*}

  For the second term,
  notice that 
  \begin{align*}
    \int^\tnjp_\tnj  dt -  \int^\tnj_\tnjm  dt=0.
  \end{align*}
  Also for the third term, we have 
  \begin{align*}
    \int^\tnjp_\tnj dB_t \int^t_\tnj dB_s - 
    \int^\tnj_\tnjm dB_t \int^t_\tnj dB_s
    &=
    \half\cbr{B(\bbone^n_{j+1})^2 + B(\bbone^n_{j})^2}
    % \\&=
    % \half\cbr{
    %   \abr{{\bbone^n_{j+1}},{\bbone^n_{j+1}}}+
    %   I_2({\bbone^n_{j+1}}^{\otimes2}) + 
    %   \abr{{\bbone^n_{j}},{\bbone^n_{j}}}+
    % I_2({\bbone^n_{j}}^{\otimes2})}.
    \\&=
    \half\cbr{
    2\times n^{-2H}+
    I_2({\bbone^n_{j+1}}^{\otimes2}) + 
    I_2({\bbone^n_{j}}^{\otimes2})}.
  \end{align*}
  
  Since we have defined 
  \begin{align*}
    \bbone^+_{n,j}=\bbone^n_{j}(t) (t-\tnjm) n \tand
    \bbone^-_{n,j}=\bbone^n_{j}(t) (\tnj-t) n,
  \end{align*}
  for the fourth and fifth terms, we have 
  \begin{align*}
    \int^\tnjp_\tnj dB_t \int^t_\tnj ds - 
    \int^\tnj_\tnjm dB_t \int^t_\tnj ds
    &=
    n^{-1}B(\bbone^+_{n,j+1}) +
    n^{-1}B(\bbone^-_{n,j})
    \\
    \int^\tnjp_\tnj dt \int^t_\tnj dB_s - 
    \int^\tnj_\tnjm dt \int^t_\tnj dB_s
    &=
    \int^\tnjp_\tnj dB_s \int^\tnjp_s dt 
    +
    \int^\tnj_\tnjm dB_s \int^s_\tnjm dt 
    \\&=
    n^{-1}B(\bbone^-_{n,j+1})
    +n^{-1}B(\bbone^+_{n,j}).
  \end{align*}

  Therefore,
  we have 
  \begin{align*}
    \secDiff{n}{j}X &= 
    \fvo_\tnj \times B(\diffker^n_j)+
    \\&\quad+
    \fvoo_\tnj \times
    \cbr{n^{-2H}+
    \half I_2({\bbone^n_{j+1}}^{\otimes2}) + 
    \half I_2({\bbone^n_{j}}^{\otimes2})}
    \\&\quad+
    \fvot_\tnj \times
    n^{-1}\rbr{B(\bbone^+_{n,j+1})+B(\bbone^-_{n,j})}
    \\&\quad+
    \fvto_\tnj \times
    n^{-1}\rbr{B(\bbone^-_{n,j+1})+B(\bbone^+_{n,j})}
    \\&\quad+
    \hat O(n^{-2}\vee n^{-3H})
  \end{align*}
\end{proof}

% \newpage
\subsubsection{Proof of Lemma \ref{lemma:240523.1807}
(Decomposition of $S_n^{(0)}-S_\infty$)}
\label{sec:240606.1341}
% \subsubsection{Decomposition of $S_n^{(0)}-S_\infty$}
Recall that we have defined
\begin{align*}
  S_n^{(0)}&=
  \cnstVar{2k}{0}\times
  n^{-1}
  \sum_{j=1}^{n-1} 
  a(X_\tnj)\tand
  % \\
  S_\infty= 
  \cnstVar{2k}{0}\times
  \int^1_0 a(X_t) dt.
  % \int^1_0 f(X_t) (V^{[1]}(X_t))^{2k} dt.
  % \int^1_0 f_t (V^{[1]}_t)^{2k} dt.
\end{align*}
For brevity, we write 
\begin{align*}
  S_n^{(0)\dagger}:=
  (\cnstVar{2k}{0})^{-1} S_n^{(0)},
  \tand
  S_\infty^\dagger:=
  (\cnstVar{2k}{0})^{-1} S_\infty.
\end{align*}

The functional $S_\infty^\dagger$ decomposes as follows
\begin{align*}
  S_\infty^\dagger&=
  % (\cnstVar{2k}{0})^{-1} S_\infty&= 
  % \int^1_0 a(X_t) dt
  % \\&=
  \sum_{j=1}^{n-1}
  \int^{\frac{2j+1}{2n}}_{\frac{2j-1}{2n}} a(X_t) dt +
  \int^\frac{1}{2n}_{0} a(X_t) dt +
  \int^1_{1-\frac{1}{2n}} a(X_t) dt 
  % \\&=
  % \sum_{j=1}^{n-1}
  % \int^{\frac{2j+1}{2n}}_{\frac{2j-1}{2n}} a(X_t) dt +
  % \\&\quad+
  % a(X_0) \int^\frac{1}{2n}_{0} dt +
  % a(X_1) \int^1_{1-\frac{1}{2n}} dt
  % \\&\quad+
  % \int^\frac{1}{2n}_{0} \cbr{a(X_t) - a(X_0)} dt +
  % \int^1_{1-\frac{1}{2n}} \cbr{a(X_t) - a(X_1)} dt
  \\&=
  {S_n^{(\infty;0)\dagger}+S_n^{(\infty;1)\dagger}+S_n^{(\infty;2)\dagger}},
\end{align*}
where we define 
\begin{align*}
  S_n^{(\infty;0)\dagger}&=
  \sum_{j=1}^{n-1}
  \int^{\frac{2j+1}{2n}}_{\frac{2j-1}{2n}} a(X_t) dt
  \\
  S_n^{(\infty;1)\dagger}&=
  % a(X_1) \int^1_{1-\frac{1}{2n}}  dt +
  % a(X_0) \int^\frac{1}{2n}_{0} dt
  % =
  \frac{1}{2n}%\times
  \brbr{a(X_0)+a(X_1)}
  \\
  S_n^{(\infty;2)\dagger}&=
  \int^\frac{1}{2n}_{0} \cbr{a(X_t) - a(X_0)} dt +
  \int^1_{1-\frac{1}{2n}} \cbr{a(X_t) - a(X_1)} dt.
  % =\hat O(n^{-H-1})
\end{align*}
We can observe that 
$S_n^{(\infty;2)\dagger}
% =
% \int^\frac{1}{2n}_{0} \cbr{a(X_t) - a(X_0)} dt +
% \int^1_{1-\frac{1}{2n}} \cbr{a(X_t) - a(X_1)} dt
=\hat O(n^{-H-1})$
by the estimate of Lemma \ref{lemma:240617.2154} (ii).

We denote the difference 
$S_n^{(0)\dagger}-
% $(\cnstVar{2k}{0})^{-1} S_n^{(0)}-
S_n^{(\infty;0)\dagger}$
by $S_n^{(0;1)\dagger}$,
which is decomposed as follows:
\begin{align*}
  S_n^{(0;1)\dagger}&=
  % (\cnstVar{2k}{0})^{-1} S_n^{(0)}-S_n^{(\infty;0)\dagger}
  % \\&=
  n^{-1} \sum_{j=1}^{n-1} a(X_\tnj) - 
  \sum_{j=1}^{n-1}
  \int^{\frac{2j+1}{2n}}_{\frac{2j-1}{2n}} a(X_t) dt
  \\&=
  \sum_{j=1}^{n-1}
  \int^{\frac{2j+1}{2n}}_{\frac{2j-1}{2n}} 
  \cbr{a(X_\tnj)-a(X_t)} dt.
\end{align*}
For $t\in[\frac{2j-1}{2n},\frac{2j+1}{2n}]$,
the difference ${a(X_\tnj)-a(X_t)}$ is expanded as 
\begin{align*}
  a(X_\tnj)-a(X_t)&=
  % \int^\tnj_t da(X_t)
  % =
  % \int^\tnj_t a'(X_t) dX_t
  % \nn\\&=
  \int^\tnj_t a'(X_s) (V^{[1]}_s dB_s + V^{[2]}_s ds)
  % \nn\\&=
  % \int^\tnj_t a'_t V^{[1]}_t dB_t + 
  % \int^\tnj_t a'_t V^{[2]}_t dt
  \nn\\&=
  a'(X_\tnj) V^{[1]}_\tnj (B_\tnj-B_t) + 
  % a'_\tnj V^{[1]}_\tnj (B_\tnj-B_t) + 
  \int^\tnj_t \cbr{a'(X_s) V^{[1]}_s - a'(X_\tnj) V^{[1]}_\tnj} dB_s 
  \nn\\&\quad+ 
  a'(X_\tnj) V^{[2]}_\tnj (\tnj-t) + 
  % a'_\tnj V^{[2]}_\tnj (\tnj-t) + 
  \int^\tnj_t \cbr{a'(X_s) V^{[2]}_s - a'(X_\tnj) V^{[2]}_\tnj} ds.
  % \int^\tnj_t \cbr{a'_t V^{[2]}_t - a'_\tnj V^{[2]}_\tnj} dt
\end{align*}
Hence 
$S_n^{(0;1)\dagger}$ can be written as 
% \begin{align*}
$S_n^{(0;1)\dagger}=\sum_{i=1}^{4} S_n^{(0;1;i)\dagger}$
% \end{align*}
with 
\begin{align*}
  S_n^{(0;1;1)\dagger}&=
  \sum_{j=1}^{n-1}
  \int^{\frac{2j+1}{2n}}_{\frac{2j-1}{2n}} 
  a'(X_\tnj) V^{[1]}_\tnj (B_\tnj-B_t) dt
  % \\&=
  % \sum_{j=1}^{n-1}
  % a'(X_\tnj) V^{[1]}_\tnj
  % \int^{\frac{2j+1}{2n}}_{\frac{2j-1}{2n}} 
  %  (B_\tnj-B_t) dt
  \\
  S_n^{(0;1;2)\dagger}&=
  \sum_{j=1}^{n-1}
  \int^{\frac{2j+1}{2n}}_{\frac{2j-1}{2n}} 
  a'(X_\tnj) V^{[2]}_\tnj (\tnj-t) dt
  \\
  S_n^{(0;1;3)\dagger}&=
  \sum_{j=1}^{n-1}
  \int^{\frac{2j+1}{2n}}_{\frac{2j-1}{2n}}
  \int^\tnj_t \cbr{a'(X_s) V^{[1]}_s - a'(X_\tnj) V^{[1]}_\tnj} dB_s dt
  \\
  S_n^{(0;1;4)\dagger}&=
  \sum_{j=1}^{n-1}
  \int^{\frac{2j+1}{2n}}_{\frac{2j-1}{2n}} 
  \int^\tnj_t \cbr{a'(X_s) V^{[2]}_s - a'(X_\tnj) V^{[2]}_\tnj} ds dt.
\end{align*}
By an argument similar to the proof of Lemma \ref{lemma:240523.1445},
we can show that 
$S_n^{(0;1;3)\dagger}=\hat O_M(n^{-2H})$ and
$S_n^{(0;1;4)\dagger}=\hat O_M(n^{-H-1})$.
Also we have $S_n^{(0;1;2)\dagger}=0$ since 
\begin{align*}
  \int^{\frac{2j+1}{2n}}_{\frac{2j-1}{2n}} 
  % a'(X_\tnj) V^{[2]}_\tnj
  (\tnj-t) dt=0.
\end{align*}

The following integral can be written as follows:
\begin{align*}
  \int^{\frac{2j+1}{2n}}_{\frac{2j-1}{2n}} (B_\tnj-B_t) dt
  &=% \\&=
  % \int^{\frac{2j+1}{2n}}_{\frac{2j-1}{2n}}
  % \int^\tnj_t dB_s dt
  % \\&=
  % \int^{\frac{2j+1}{2n}}_{\tnj}
  % \int^\tnj_t dB_s dt
  % +
  % \int^{\tnj}_{\frac{2j-1}{2n}}
  % \int^\tnj_t dB_s dt
  % \\&=
  % -\int^{\frac{2j+1}{2n}}_{\tnj}
  % \int^t_\tnj dB_s dt
  % +
  % \int^{\tnj}_{\frac{2j-1}{2n}}
  % \int^\tnj_t dB_s dt
  % \\&=
  -\int^{\frac{2j+1}{2n}}_{\tnj}
  \int^{\frac{2j+1}{2n}}_s  dt dB_s
  +
  \int^{\tnj}_{\frac{2j-1}{2n}}
  \int^s_{\frac{2j-1}{2n}}  dt dB_s
  \\&=
  % -\int^{\frac{2j+1}{2n}}_{\tnj}
  % \cbr{{\frac{2j+1}{2n}}-s} dB_s
  % +
  % \int^{\tnj}_{\frac{2j-1}{2n}}
  % \cbr{s-{\frac{2j-1}{2n}}} dB_s
  % \\&=
  2^{-1}% \half\times 
  n^{-1}\times \cbr{
  -B\brbr{\bbone^-_{2n,2j+1}} + B\brbr{\bbone^+_{2n,2j}}},
\end{align*}
where we recall that we have defined at \eqref{eq:240606.1219}
\begin{align*}
  \bbone^+_{n,j}(t)=\bbone^n_{j}(t) (t-\tnjm) n \tand
  \bbone^-_{n,j}(t)=\bbone^n_{j}(t) (\tnj-t) n.
\end{align*}
Thus we can write $S_n^{(0;1;1)\dagger}$ as follows:
\begin{align*}
  S_n^{(0;1;1)\dagger}&=
  % \sum_{j=1}^{n-1}
  % \int^{\frac{2j+1}{2n}}_{\frac{2j-1}{2n}} 
  % a'(X_\tnj) V^{[1]}_\tnj (B_\tnj-B_t) dt
  % \\&=
  % \sum_{j=1}^{n-1}
  % a'(X_\tnj) V^{[1]}_\tnj
  % \int^{\frac{2j+1}{2n}}_{\frac{2j-1}{2n}} (B_\tnj-B_t) dt
  % \\&=
  2^{-1} n^{-1}
  \sum_{j=1}^{n-1}
  a'(X_\tnj) V^{[1]}_\tnj \times 
  \cbr{-B\brbr{\bbone^-_{2n,2j+1}} + B\brbr{\bbone^+_{2n,2j}}},
  % \label{eq:240606.1236}
\end{align*}
By the argument of exponent, we can show that 
$S_n^{(0;1;1)\dagger}
% =n^{-1}O_M(n^{1+(-1)})
=O_M(n^{-1})$.

Therefore, we obtain 
\begin{align*}
  S_n^{(0)}-S_\infty&=
  \cnstVar{2k}{0}\times
  (S_n^{(0)\dagger}-
  % ((\cnstVar{2k}{0})^{-1} S_n^{(0)}-
  \ncbr{S_n^{(\infty;0)\dagger}+S_n^{(\infty;1)\dagger}+S_n^{(\infty;2)\dagger}})
  % \\&=
  % \cnstVar{2k}{0}\times
  % ((\cnstVar{2k}{0})^{-1} S_n^{(0)}-S_n^{(\infty;0)\dagger}
  % -\ncbr{S_n^{(\infty;1)\dagger}+S_n^{(\infty;2)\dagger}})
  % \\&=
  % \cnstVar{2k}{0}\times
  % (S_n^{(0;1)\dagger}
  % -\ncbr{S_n^{(\infty;1)\dagger}+S_n^{(\infty;2)\dagger}})
  \\&=
  \cnstVar{2k}{0}\times
  (\sum_{i=1}^{4} S_n^{(0;1;i)\dagger}
  -\ncbr{S_n^{(\infty;1)\dagger}+S_n^{(\infty;2)\dagger}})
  % \\&=
  % \cnstVar{2k}{0}\times
  % (S_n^{(0;1;1)\dagger}+
  % S_n^{(0;1;3)\dagger}+S_n^{(0;1;4)\dagger}
  % -\ncbr{S_n^{(\infty;1)\dagger}+S_n^{(\infty;2)\dagger}})
  % \\&=
  % \cnstVar{2k}{0}\times
  % (S_n^{(0;1;1)\dagger}
  % -S_n^{(\infty;1)\dagger})+
  % \hat O_M(n^{-2H})+\hat O_M(n^{-H-1})+\hat O(n^{-H-1})
  \\&=
  \cnstVar{2k}{0}\times
  (S_n^{(0;1;1)\dagger}
  -S_n^{(\infty;1)\dagger})+
  \hat O_M(n^{-2H}).
\end{align*}
\rightline{\qed}

% \newpage

% \newpage
\subsection{Lemmas related to limits of functionals}
\label{sec:240607.1650}
Let 
\begin{align}
  \widehat\rho(j)&=
  \half\cbr{-\abs{j+2}^{2H}
  +4\abs{j+1}^{2H}
  -6\abs{j}^{2H}
  +4\abs{j-1}^{2H}
  -\abs{j-2}^{2H}}.
\end{align}
When we write 
$\diffker^1_j=\bbone_{[j,j+1]}-\bbone_{[j-1,j]}$ for $j\in\bbZ$,
we have 
\begin{align*}
  \abr{\diffker^1_{0},\diffker^1_j}=\widehat\rho(j),
\end{align*}
since 
\begin{align*}
  \abr{\bbone_{[0,1]},\bbone_{[j,j+1]}}=
  \half\cbr{\abs{j+1}^{2H}-2\abs{j}^{2H}+\abs{j-1}^{2H}}.
\end{align*}
% \begin{align}
%   \phi(j)&:=
%   \abr{\bbone_{[0,1]},\bbone_{[j,j+1]}}
%   \nn\\&=
%   \half\cbr{\abs{j+1}^{2H}-\abs{j}^{2H}-
%   (\abs{j}^{2H}-\abs{j-1}^{2H})}
%   \nn\\&=
%   \half\cbr{\abs{j+1}^{2H}-2\abs{j}^{2H}+\abs{j-1}^{2H}}
%   \label{eq:240620.1504}
% \end{align}
(Here the inner product is considered on the Hilbert space associated to 
the fBm on $\bbR$.) %with the Hurst parameter $H$.

The following lemma gives the basic properties of $\widehat\rho$.
\begin{lemma}\label{lemma:240620.1743}
(i) The following relation holds for $j_1,j_2\in[n-1]$:
\begin{align}
  \babr{\diffker^{n}_{j_1},\diffker^{n}_{j_2}}=
  n^{-2H}\widehat\rho(j_2-j_1).
  \label{eq:240617.2141}
\end{align}

\item[(ii)] It holds that 
$\widehat\rho(j) = O(\abs{j}^{2H-4})$ as $\abs{j}\to\infty$.
In other words,
there exists $C>0$ such that
$\abs{\widehat\rho(j)} \leq C(\abs{j}\vee1)^{2H-4}$ for $j\in\bbZ$.
\end{lemma}

\begin{proof}
(i)
It follows from the self-similarity property of the fBm.

\item{(ii)}
This follows from %Lemma 3.3 and 
Lemma 3.4 of \cite{mishura2023asymptotic}.
\end{proof}

\subsubsection{Proof of Lemma \ref{lemma:240611.2313}}\label{sec:240611.2314}

\begin{proof}
From \eqref{eq:240617.2143} and \eqref{eq:240617.2141},
we have
\begin{align*}
  \GnTwoMain&=
  % n^{4H\ell_1-1} 
  % \sum_{j_1,j_2=1}^{n-1}
  % a(X_\tnjo) a(X_\tnjt)
  % \abr{\diffker^{n}_{j_1},\diffker^{n}_{j_2}}^{2\ell_1}
  % \\&=
  % n^{4H\ell_1-1} 
  % \sum_{j_1,j_2=1}^{n-1}
  % a(X_\tnjo) a(X_\tnjt)
  % \rbr{n^{-2H}\widehat\rho(j_2-j_1)}^{2\ell_1}
  % \\&=
  % n^{4H\ell_1-1} n^{-2H2\ell_1}
  % \sum_{j_1,j_2=1}^{n-1}
  % a(X_\tnjo) a(X_\tnjt)
  % \widehat\rho(j_2-j_1)^{2\ell_1}
  % \\&=
  n^{-1} 
  \sum_{j_1,j_2=1}^{n-1}
  a(X_\tnjo) a(X_\tnjt)
  \widehat\rho(j_2-j_1)^{2\ell_1}.
\end{align*}
We write $\cali_n^{(\ell_1)\dagger}:=\GnTwoMain$ for short.
We define the following functionals:  
\begin{align*}
  % G_n^{(\ell_1)\dagger}&:=
  % n^{-1} 
  % \sum_{j_1,j_2\in[n-1]}
  % a(X_\tnjo) a(X_\tnjt)
  % \widehat\rho(j_2-j_1)^{2\ell_1}
  % \quad\text{(参考)}
  % \\
  \cali_n^{(\ell_1;1)\dagger}&=
  n^{-1} 
  \sum_{j_1,j_2\in[n-1]}
  a(X_\tnjo)^2
  \widehat\rho(j_2-j_1)^{2\ell_1}
  % \\&=
  % n^{-1} 
  % \sum_{j_1\in[n-1]}
  % a(X_\tnjo)^2
  % \sum_{j_2\in[n-1]}
  % \widehat\rho(j_2-j_1)^{2\ell_1}
  \\
  \cali_n^{(\ell_1;2)\dagger}&=
  n^{-1} 
  \sum_{j_1\in[n-1]}
  a(X_\tnjo)^2
  \sum_{i\in\bbZ}
  \widehat\rho(i)^{2\ell_1}
  \\
  \cali_n^{(\ell_1;\infty)\dagger}&=
  \int^1_0 a(X_t)^2 dt \times
  \sum_{i\in\bbZ}
  \widehat\rho(i)^{2\ell_1}
\end{align*}

Then the difference of 
$\cali_n^{(\ell_1)\dagger}$ and $\cali_n^{(\ell_1;1)\dagger}$ 
can be written as 
\begin{align*}
  &\cali_n^{(\ell_1)\dagger}-\cali_n^{(\ell_1;1)\dagger}
  =% \\&=
  % n^{-1} 
  % \sum_{j_1,j_2\in[n-1]}
  % a(X_\tnjo) a(X_\tnjt)
  % \widehat\rho(j_2-j_1)^{2\ell_1}
  % -
  % n^{-1} 
  % \sum_{j_1\in[n-1]}
  % a(X_\tnjo)^2
  % \sum_{j_2\in[n-1]}
  % \widehat\rho(j_2-j_1)^{2\ell_1}
  % \\&=
  n^{-1} 
  \sum_{j_1,j_2\in[n-1]}
  a(X_\tnjo) \rbr{a(X_\tnjt)-a(X_\tnjo)}
  \widehat\rho(j_2-j_1)^{2\ell_1},
\end{align*}
which is estimated using Lemma \ref{lemma:240617.2154} (ii) as 
% SDEの解(の関数値)の差についての評価より，
\begin{align}
  % &
  \norm{\cali_n^{(\ell_1)\dagger}-\cali_n^{(\ell_1;1)\dagger}}_{k,p}
  &\leq% \\&\leq
  n^{-1} 
  \sum_{j_1,j_2\in[n-1]}
  \norm{a(X_\tnjo)}_{k,2p} \norm{a(X_\tnjt)-a(X_\tnjo)}_{k,2p}
  \widehat\rho(j_2-j_1)^{2\ell_1}
  \nn\\&\simleq
  n^{-1} 
  \sum_{j_1,j_2\in[n-1]}
  % \norm{a(X_\tnjt)-a(X_\tnjo)}_{k,2p}
  \rbr{\abs{j_2-j_1}/n}^{\beta}
  \widehat\rho(j_2-j_1)^{2\ell_1}
  \nn\\&=
  n^{-1-\beta} 
  \sum_{j_1,j_2\in[n-1]}
  \abs{j_2-j_1}^{\beta}
  \widehat\rho(j_2-j_1)^{2\ell_1}
  =O(n^{-\beta})
  \label{eq:240617.2221}
\end{align}
for any $\beta\in(\half,H)$.
% {\rb $\hat O_M$で書き直す．}

For $\cali_n^{(\ell_1;1)\dagger}-\cali_n^{(\ell_1;2)\dagger}$,
we first consider the following difference:
\begin{align*}
  &n^{-1} 
  \sum_{j_1\in[n-1]}
  \Bcbr{
  \sum_{i\in\bbZ}
  \widehat\rho(i)^{2\ell_1}
  -
  \sum_{j_2\in[n-1]}
  \widehat\rho(j_2-j_1)^{2\ell_1}}
  \\&=
  n^{-1} 
  \sum_{j_1\in[n-1]}
  \Bcbr{
  \sum_{i\in\bbZ}
  \widehat\rho(i)^{2\ell_1}
  -
  \sum_{i\in[n-1]_{j_1}}
  \widehat\rho(i)^{2\ell_1}}
  \\&=
  n^{-1} 
  \sum_{j_1\in[n-1]}
  \sum_{i\in\bbZ\setminus [n-1]_{j_1}}
  \widehat\rho(i)^{2\ell_1},
\end{align*}
where we write 
$[n-1]_{j_1}={1-j_1,...,(n-1)-j_1}$.

Let $\epsilon\in(0,1)$.
When $n^\epsilon<j_1<n-n^\epsilon$, it holds that
$-j_1<-n^\epsilon$ and $n^\epsilon<n-j_1$, and
% $j_1<n-n^\epsilon$
% $n^\epsilon<n-j_1$
if $i\in\bbZ\setminus [n-1]_{j_1}$, then 
we have $\abs{i}>n^\epsilon$.
Hence, 
\begin{align*}
  % &
  n^{-1} 
  \sum_{n^\epsilon<j_1<n-n^\epsilon}
  % \Bcbr{
  % \sum_{i\in\bbZ}
  % \widehat\rho(i)^{2\ell_1}
  % -
  % \sum_{j_2\in[n-1]}
  % \widehat\rho(j_2-j_1)^{2\ell_1}}
  \sum_{i\in\bbZ\setminus [n-1]_{j_1}}
  \widehat\rho(i)^{2\ell_1}
  &\leq% \\&\leq
  n^{-1} 
  \sum_{n^\epsilon<j_1<n-n^\epsilon}
  \sum_{\abs{i}>n^\epsilon}
  \widehat\rho(i)^{2\ell_1}
  \\*&=
  n^{-1} \times O(n^1) \times
  O(n^{\epsilon((2H-4)2\ell_1+1)})
  \\*&=O(n^{\epsilon((2H-4)2\ell_1+1)})
  \leq O(n^{-3\epsilon}).
\end{align*}
Using the above estimate, we can show
\begin{align*}
  &\norm{\cali_n^{(\ell_1;1)\dagger}-\cali_n^{(\ell_1;2)\dagger}}_{k,p}
  \\&\leq
  \sup_{j_1}\norm{a(X_\tnjo)^2}_{k,p}
  n^{-1} 
  \sum_{j_1\in[n-1]}
  \abs{\sum_{j_2\in[n-1]}
  \widehat\rho(j_2-j_1)^{2\ell_1}
  -
  \sum_{i\in\bbZ}
  \widehat\rho(i)^{2\ell_1}}
  \\&=
  \sup_{j_1}\norm{a(X_\tnjo)^2}_{k,p}
  n^{-1} 
  \sum_{n^\epsilon<j_1<n-n^\epsilon}
  \abs{\sum_{j_2\in[n-1]}
  \widehat\rho(j_2-j_1)^{2\ell_1}
  -
  \sum_{i\in\bbZ}
  \widehat\rho(i)^{2\ell_1}}
  \\&\quad+
  \sup_{j_1}\norm{a(X_\tnjo)^2}_{k,p}
  n^{-1} 
  \sum_{\substack{1\leq j_1\leq n^\epsilon\\n-n^\epsilon\leq j_1\leq n-1}}
  \abs{\sum_{j_2\in[n-1]}
  \widehat\rho(j_2-j_1)^{2\ell_1}
  -
  \sum_{i\in\bbZ}
  \widehat\rho(i)^{2\ell_1}}
  \\&=
  % O(n^{\epsilon((2H-4)2\ell_1+1)})+
  % O(n^{-1+\epsilon})
  % \\&=
  O(n^{-3\epsilon})+
  O(n^{-1+\epsilon}).
\end{align*}
Setting $\frac14=\epsilon$,
we have
\begin{align}
  &\norm{\cali_n^{(\ell_1;1)\dagger}-\cali_n^{(\ell_1;2)\dagger}}_{k,p}=
  O(n^{-\frac34}).
  \label{eq:240617.2222}
\end{align}

Since we have 
\begin{align*}
  \cali_n^{(\ell_1;2)\dagger}-\cali_n^{(\ell_1;\infty)\dagger}
  &=O_M(n^{-1})
\end{align*}
by Lemma \ref{lemma:240617.2154} (iii),
we obtain from \eqref{eq:240617.2221} and \eqref{eq:240617.2222}
the following estimate
\begin{align*}
  \GnTwoMain-
  \cali_n^{(\ell_1;\infty)\dagger}=
  \cali_n^{(\ell_1)\dagger}-
  \cali_n^{(\ell_1;\infty)\dagger}&=
  \hat O_M(n^{-H})+O_M(n^{-\frac34}).
\end{align*}

\end{proof}

\subsubsection{Proof of Lemma \ref{lemma:240611.2018}}\label{sec:240611.2231}
\begin{proof}%[Proof of Lemma \ref{lemma:240611.2018}]
  Recall that we have 
  \begin{align*}
    (D_{u_n})^2M_n&=
    \sum_{(\ell_1,\ell_2,\ell_3,m)\in\Lambda^{(k;3)}_{\sharp}}
    \qtorconst \times\caliqtor
    +O_M(n^{-H})
  \end{align*}
  with some positive integer $\qtorconst$ defined at \eqref{eq:240630.1831} and 
\begin{align*}
  \caliqtor&=
  n^{2(\ell_1+\ell_2+\ell_3) H-\frac32}
  \sum_{j_1,j_2,j_3\in[n-1]}
  a(X_{\tnjo}) a(X_{\tnjt}) a(X_\tnjs)
  \nn\\*&\hspace{50pt}\times
  \abr{\diffker^{n}_{j_1},\diffker^{n}_{j_2}}^{m+1}
  \abr{\diffker^{n}_{j_1},\diffker^{n}_{j_3}}^{q_{\ell_1,m}}
  \abr{\diffker^{n}_{j_2},\diffker^{n}_{j_3}}^{q_{\ell_2,m}}.
  % \label{eq:240610.1145}
\end{align*}
Note that at least one of 
$q_{\ell_1,m}$ and $q_{\ell_2,m}$ is positive.
For brevity, we write
$\cali_{\qtor,n}^{(\ell,m)}:=
n^{\half}\times\caliqtor$.
%
% {\redbf By the self-similarity property of fBm,
% we have 
% \begin{align*}
%   \babr{\diffker^{n}_{j_1},\diffker^{n}_{j_2}}=
%   n^{-2H}\widehat\rho(j_2-j_1),
% \end{align*}}
Using the relation \eqref{eq:240617.2141},
we can write 
\begin{align*}
  \cali_{\qtor,n}^{(\ell,m)}&=
  n^{-1}
  % n^{2(\ell_1+\ell_2+\ell_3) H-1}
  % n^{-2H(m+1 + q_{\ell_1,m} + q_{\ell_2,m})}
  % \\*&\hspace{50pt}\times 
  \sum_{j_1,j_2,j_3\in[n-1]}
  a(X_{\tnjo}) a(X_{\tnjt}) a(X_\tnjs)
  \widehat\rho(j_2-j_1)^{m+1}
  \widehat\rho(j_3-j_1)^{q_{\ell_1,m}}
  \widehat\rho(j_3-j_2)^{q_{\ell_2,m}},
\end{align*}
where we have used 
$2(\ell_1+\ell_2+\ell_3)=
2(m+1 + q_{\ell_1,m} + q_{\ell_2,m})$
since 
$(\ell_1,\ell_2,\ell_3,m)\in\Lambda^{(k;3)}_{\sharp}
=
\cbr{(\ell_1,\ell_2,\ell_3,m)\in\Lambda^{(k;3)}\mid 
2\ell_3=q_{\ell_1,m}+q_{\ell_2,m}}$.

Assume that $q_{\ell_1,m}>0$ for a while.
Decompose the functional $\cali_{\qtor,n}^{(\ell,m)}$ as
\begin{align}
  \cali_{\qtor,n}^{(\ell,m)}&=
  \cali_{\qtor,n}^{(\ell,m;1)}+
  \cali_{\qtor,n}^{(\ell,m;2)}+
  \cali_{\qtor,n}^{(\ell,m;3)}
  \label{eq:240611.2151}
\end{align}
with 
\begin{align*}
  \cali_{\qtor,n}^{(\ell,m;1)}&=
  n^{-1}
  \sum_{j_1,j_2,j_3\in[n-1]}
  a(X_{\tnjo})^3
  % a(X_{\tnjo}) a(X_{\tnjt}) a(X_\tnjs)
  \widehat\rho(j_2-j_1)^{m+1}
  \widehat\rho(j_3-j_1)^{q_{\ell_1,m}}
  \widehat\rho(j_3-j_2)^{q_{\ell_2,m}}
  \\
  \cali_{\qtor,n}^{(\ell,m;2)}&=
  n^{-1}
  \sum_{j_1,j_2,j_3\in[n-1]}
  a(X_{\tnjo})^2
  (a(X_\tnjs)-a(X_{\tnjo}))
  \widehat\rho(j_2-j_1)^{m+1}
  \widehat\rho(j_3-j_1)^{q_{\ell_1,m}}
  \widehat\rho(j_3-j_2)^{q_{\ell_2,m}}
  \\
  \cali_{\qtor,n}^{(\ell,m;3)}&=
  n^{-1}
  \sum_{j_1,j_2,j_3\in[n-1]}
  a(X_{\tnjo}) a(X_\tnjs)
  (a(X_\tnjt)-a(X_{\tnjo}))
  \widehat\rho(j_2-j_1)^{m+1}
  \widehat\rho(j_3-j_1)^{q_{\ell_1,m}}
  \widehat\rho(j_3-j_2)^{q_{\ell_2,m}}.
\end{align*}
The Malliavin norm of $\cali_{\qtor,n}^{(\ell,m;2)}$ is estimated as follows:
\begin{align}
  &\norm{\cali_{\qtor,n}^{(\ell,m;2)}}_{k,p}
  % \\&\leq
  % n^{-1}
  % \sum_{j_1,j_2,j_3\in[n-1]}
  % \norm{a(X_{\tnjo})^2
  % (a(X_\tnjs)-a(X_{\tnjo}))}_{k,p}
  % \abs{\widehat\rho(j_2-j_1)}^{m+1}
  % \abs{\widehat\rho(j_3-j_1)}^{q_{\ell_1,m}}
  % \abs{\widehat\rho(j_3-j_2)}^{q_{\ell_2,m}}
  \nn\\&\simleq
  n^{-1}
  \sum_{j_1,j_2,j_3\in[n-1]}
  \norm{a(X_{\tnjo})}_{k,3p}^2
  \norm{(a(X_\tnjs)-a(X_{\tnjo}))}_{k,3p}
  \nn\\&\hspace{40pt}
  \abs{\widehat\rho(j_2-j_1)}^{m+1}
  \abs{\widehat\rho(j_3-j_1)}^{q_{\ell_1,m}}
  \abs{\widehat\rho(j_3-j_2)}^{q_{\ell_2,m}}
  \nn\\&\simleq
  n^{-1}
  \sum_{j_1,j_2,j_3\in[n-1]}
  O(1)
  \rbr{\abs{j_3-j_1}/n}^{\beta}
  % \norm{(a(X_\tnjs)-a(X_{\tnjo}))}_{k,3p}
  \abs{\widehat\rho(j_2-j_1)}^{m+1}
  \abs{\widehat\rho(j_3-j_1)}^{q_{\ell_1,m}}
  \abs{\widehat\rho(j_3-j_2)}^{q_{\ell_2,m}}
  \label{eq:240611.2051}
  \\&\simleq
  n^{-1-\beta}
  \sum_{j_1\in[n-1]}
  \sum_{j_3\in[n-1]}
  (1\vee\abs{j_3-j_1})^{\beta+(2H-4)q_{\ell_1,m}}
  \sum_{j_2\in[n-1]}
  (1\vee\abs{j_2-j_1})^{(2H-4)(m+1)}
  \nn\\&=
  O(n^{-\beta})
  \label{eq:240611.2154}
\end{align}
with any $\beta\in(\half,H)$.
At \eqref{eq:240611.2051}, we used the estimate by Lemma \ref{lemma:240617.2154} (ii).
Similarly, 
the Malliavin norm of $\cali_{\qtor,n}^{(\ell,m;3)}$ is estimated as 
\begin{align}
  \norm{\cali_{\qtor,n}^{(\ell,m;3)}}_{k,p}
  =O(n^{-\beta}).
  \label{eq:240611.2155}
\end{align}

Define the following functional:
\begin{align*}
  % \cali_{\qtor,n}^{(\ell,m;1)}&=
  % n^{-1}
  % \sum_{j_1,j_2,j_3\in[n-1]}
  % a(X_{\tnjo})^3
  % \widehat\rho(j_2-j_1)^{m+1}
  % \widehat\rho(j_3-j_1)^{q_{\ell_1,m}}
  % \widehat\rho(j_3-j_2)^{q_{\ell_2,m}}
  % \\&=
  % n^{-1}
  % \sum_{j_1\in[n-1]}
  % a(X_{\tnjo})^3
  % \sum_{j_2,j_3\in[n-1]}
  % \widehat\rho(j_2-j_1)^{m+1}
  % \widehat\rho(j_3-j_1)^{q_{\ell_1,m}}
  % \widehat\rho(j_3-j_2)^{q_{\ell_2,m}}
  % \\
  \cali_{\qtor,n}^{(\ell,m;\infty)}&=
  n^{-1}
  \sum_{j_1\in[n-1]}
  a(X_{\tnjo})^3
  \sum_{i_1,i_2\in\bbZ}
  \widehat\rho(i_1)^{m+1}
  \widehat\rho(i_2)^{q_{\ell_1,m}}
  \widehat\rho(i_2-i_1)^{q_{\ell_2,m}}.
\end{align*}
Also we write 
$i_1:=j_2-j_1$, $i_2:=j_3-j_1$ and 
$[n-1]_{j}:=[n-1]-j=\cbr{1-j,...,n-1-j}$.
Then we can write 
\begin{align*}
  \cali_{\qtor,n}^{(\ell,m;1)}&=
  % n^{-1}
  % \sum_{j_1\in[n-1]}
  % a(X_{\tnjo})^3
  % \sum_{j_2,j_3\in[n-1]}
  % \widehat\rho(j_2-j_1)^{m+1}
  % \widehat\rho(j_3-j_1)^{q_{\ell_1,m}}
  % \widehat\rho(j_3-j_2)^{q_{\ell_2,m}}
  % \\&=
  n^{-1}
  \sum_{j_1\in[n-1]}
  a(X_{\tnjo})^3
  \sum_{i_1,i_2\in[n-1]_{j_1}}
  % \sum_{j_2,j_3\in[n-1]}
  \widehat\rho(i_1)^{m+1}
  \widehat\rho(i_2)^{q_{\ell_1,m}}
  \widehat\rho(i_2-i_1)^{q_{\ell_2,m}}
\end{align*}
and hence 
% \begin{align*}
%   \cali_{\qtor,n}^{(\ell,m;\infty)}-\cali_{\qtor,n}^{(\ell,m;1)}&=
%   n^{-1}
%   \sum_{j_1\in[n-1]}
%   a(X_{\tnjo})^3
%   \sum_{(i_1,i_2)\in\bbZ^2\setminus([n-1]_{j_1})^2}
%   \widehat\rho(i_1)^{m+1}
%   \widehat\rho(i_2)^{q_{\ell_1,m}}
%   \widehat\rho(i_2-i_1)^{q_{\ell_2,m}}
% \end{align*}
we have
\begin{align}
  &\norm{\cali_{\qtor,n}^{(\ell,m;\infty)}-\cali_{\qtor,n}^{(\ell,m;1)}}_{p}
  \nn\\&\leq
  % n^{-1}
  % \sum_{j_1\in[n-1]}
  % \norm{a(X_{\tnjo})^3}_{p}
  % \sum_{(i_1,i_2)\in\bbZ^2\setminus([n-1]_{j_1})^2}
  % \widehat\rho(i_1)^{m+1}
  % \widehat\rho(i_2)^{q_{\ell_1,m}}
  % \widehat\rho(i_2-i_1)^{q_{\ell_2,m}}
  % \nn\\&\leq
  \sup_t\norm{a(X_{t})^3}_{p}\times
  n^{-1}
  \sum_{j_1\in[n-1]}
  \sum_{(i_1,i_2)\in\bbZ^2\setminus([n-1]_{j_1})^2}
  \abs{\widehat\rho(i_1)}^{m+1}
  \abs{\widehat\rho(i_2)}^{q_{\ell_1,m}}
  \abs{\widehat\rho(i_2-i_1)}^{q_{\ell_2,m}}.
  \label{eq:240611.2146}
\end{align}

Let us define 
\begin{align*}
  \Delta^\qtor_{n,j_1}:=
  \sum_{(i_1,i_2)\in\bbZ^2\setminus([n-1]_{j_1})^2}
  \abs{\widehat\rho(i_1)}^{m+1}
  \abs{\widehat\rho(i_2)}^{q_{\ell_1,m}}
  \abs{\widehat\rho(i_2-i_1)}^{q_{\ell_2,m}}.
\end{align*}
By Lemma \ref{lemma:240620.1743} (ii), % ($\widehat\rho$のdecayの速さについての補題)},
we have 
\begin{align*}
  \sup_{j} 
  \abs{\widehat\rho(j)}^{q_{\ell_2,m}}<\infty,\tand
  \sum_{i_1\in\bbZ}
  \abs{\widehat\rho(i_1)}^{m+1},
  \sum_{i_2\in\bbZ}
  \abs{\widehat\rho(i_2)}^{q_{\ell_1,m}}<\infty,
\end{align*}
Here we used the assumption $q_{\ell_1,m}>0$.
Hence,
\begin{align*}
  \sup_{j_1=1,...,n-1}
  \Delta^\qtor_{n,j_1}&\leq
  \sum_{(i_1,i_2)\in\bbZ^2}
  \abs{\widehat\rho(i_1)}^{m+1}
  \abs{\widehat\rho(i_2)}^{q_{\ell_1,m}}
  \abs{\widehat\rho(i_2-i_1)}^{q_{\ell_2,m}}
  \\&\leq
  \sup_{i_1,i_2}
  \abs{\widehat\rho(i_2-i_1)}^{q_{\ell_2,m}}
  \sum_{i_1\in\bbZ}
  \abs{\widehat\rho(i_1)}^{m+1}
  \sum_{i_2\in\bbZ}
  \abs{\widehat\rho(i_2)}^{q_{\ell_1,m}}
  <\infty.
\end{align*}

Consider the case where $n^\epsilon\leq j_1\leq n-n^\epsilon$
with some $0<\epsilon<1$, which will be specified later.
We can bound $\Delta^\qtor_{n,j_1}$ as 
\begin{align*}
  \Delta^\qtor_{n,j_1}
  % &\overset{\text{def}}{=}
  % \sum_{(i_1,i_2)\in\bbZ^2\setminus([n-1]_{j_1})^2}
  % \widehat\rho(i_1)^{m+1}
  % \widehat\rho(i_2)^{q_{\ell_1,m}}
  % \widehat\rho(i_2-i_1)^{q_{\ell_2,m}}
  % \\&\leq
  &\leq
  \sum_{(i_1,i_2)\in(\bbZ\setminus[n-1]_{j_1})\times\bbZ}
  \abs{\widehat\rho(i_1)}^{m+1}
  \abs{\widehat\rho(i_2)}^{q_{\ell_1,m}}
  \abs{\widehat\rho(i_2-i_1)}^{q_{\ell_2,m}}
  \\&\quad+
  \sum_{(i_1,i_2)\in\bbZ\times(\bbZ\setminus[n-1]_{j_1})}
  \abs{\widehat\rho(i_1)}^{m+1}
  \abs{\widehat\rho(i_2)}^{q_{\ell_1,m}}
  \abs{\widehat\rho(i_2-i_1)}^{q_{\ell_2,m}}
%   \\&=:
%   \Delta^{\qtor;(1)}_{n,j_1}+\Delta^{\qtor;(2)}_{n,j_1}
% \end{align*}
% 
% \begin{align*}
%   \Delta^{\qtor;(1)}_{n,j_1}&\overset{\text{def}}{=}
%   \sum_{i_1\in(\bbZ\setminus[n-1]_{j_1})}
%   \sum_{i_2\in\bbZ}
%   \widehat\rho(i_1)^{m+1}
%   \widehat\rho(i_2)^{q_{\ell_1,m}}
%   \widehat\rho(i_2-i_1)^{q_{\ell_2,m}}
  \\*&\leq
  \sup_{j} 
  \abs{\widehat\rho(j)}^{q_{\ell_2,m}}\times
  \sum_{i_1\in\bbZ\setminus[n-1]_{j_1}}
  \abs{\widehat\rho(i_1)}^{m+1}
  \sum_{i_2\in\bbZ}
  \abs{\widehat\rho(i_2)}^{q_{\ell_1,m}}
% \end{align*}
% 
% \begin{align*}
%   \Delta^{\qtor;(2)}_{n,j_1}&\overset{\text{def}}{=}
  % \sum_{(i_1,i_2)\in\bbZ\times(\bbZ\setminus[n-1]_{j_1})}
  % \widehat\rho(i_1)^{m+1}
  % \widehat\rho(i_2)^{q_{\ell_1,m}}
  % \widehat\rho(i_2-i_1)^{q_{\ell_2,m}}
  % \\*&\leq 
\\*&\quad+
  \sup_j
  \abs{\widehat\rho(j)}^{q_{\ell_2,m}} \times
  \sum_{i_1\in\bbZ}
  \abs{\widehat\rho(i_1)}^{m+1}
  \sum_{i_2\in\bbZ\setminus[n-1]_{j_1}}
  \abs{\widehat\rho(i_2)}^{q_{\ell_1,m}}
\end{align*}
By the assumption
$n^\epsilon\leq j_1\leq n-n^\epsilon$,
we have 
\begin{align*}
  \bbZ\setminus[n-1]_{j_1}&=
  \cbr{i\in\bbZ\mid i\leq-j_1 \torsm n-j_1\leq i}
  \subset 
  \cbr{i\in\bbZ\mid \abs{i}\geq n^\epsilon}
\end{align*}
and hence 
\begin{align*}
  \sum_{i_1\in\bbZ\setminus[n-1]_{j_1}}
  \abs{\widehat\rho(i_1)}^{m+1}
  \leq
  \sum_{i_1\in\cbr{i\in\bbZ\mid \abs{i}\geq n^\epsilon}}
  \abs{\widehat\rho(i_1)}^{m+1}
  =
  O(n^{\epsilon((2H-4)(m+1)+1)}).
\end{align*}
Similarly, 
\begin{align*}
  \sum_{i_2\in\bbZ\setminus[n-1]_{j_1}}
  \abs{\widehat\rho(i_2)}^{q_{\ell_1,m}}
  % \sum_{i_1\in\bbZ\setminus[n-1]_{j_1}}
  % \widehat\rho(i_1)^{m+1}
  \leq
  \sum_{i_2\in\cbr{i\in\bbZ\mid \abs{i}\geq n^\epsilon}}
  \abs{\widehat\rho(i_2)}^{q_{\ell_1,m}}
  % \widehat\rho(i_1)^{m+1}
  =
  O(n^{\epsilon((2H-4)q_{\ell_1,m}+1)}),
\end{align*}
since we are assuming that $q_{\ell_1,m}>0$.
We obtain 
\begin{align*}
  \Delta^\qtor_{n,j_1}&\leq
  O(n^{\epsilon((2H-4)(m+1)+1)})+
  O(n^{\epsilon((2H-4)q_{\ell_1,m}+1)})
  \leq O(n^{\epsilon((2H-4)+1)}),
\end{align*}
and therefore
\begin{align*}
  n^{-1}
  \sum_{j_1\in[n-1]}
  \Delta^\qtor_{n,j_1}&=
  n^{-1}
  \sum_{j_1<n^\epsilon\torsm n-n^\epsilon<j_1}
  \Delta^\qtor_{n,j_1}
  +
  n^{-1}
  \sum_{n^\epsilon\leq j_1\leq n-n^\epsilon}
  \Delta^\qtor_{n,j_1}
  \\&\leq
  n^{-1} O(n^\epsilon)+
  n^{-1} O(n^1) \times
  O(n^{\epsilon((2H-4)+1)})
  \\&=
  O(n^{(\epsilon-1)})+
  O(n^{\epsilon(2H-3)})
  =
  O(n^{(\epsilon-1)\vee(-\epsilon)}).
\end{align*}
By setting $\epsilon=\half$ and \eqref{eq:240611.2146},
it holds that 
\begin{align}
  &\norm{\cali_{\qtor,n}^{(\ell,m;\infty)}-\cali_{\qtor,n}^{(\ell,m;1)}}_{p}
  =O(n^{-\half})
  \label{eq:240611.2152}
\end{align}

Define the functional $\cali_{\qtor}^{(\ell,m;\infty)}$ by 
\begin{align}
  \cali_{\qtor}^{(\ell,m;\infty)}&=
  \int^1_0 a(X_{t})^3 dt \times
  \sum_{i_1,i_2\in\bbZ}
  \widehat\rho(i_1)^{m+1}
  \widehat\rho(i_2)^{q_{\ell_1,m}}
  \widehat\rho(i_2-i_1)^{q_{\ell_2,m}}.
  \label{eq:240604.1707}
\end{align}
Then we have
\begin{align}
  &\cali_{\qtor,n}^{(\ell,m;\infty)}-
  \cali_{\qtor}^{(\ell,m;\infty)}
  \nn\\&=
  \Bcbr{\sum_{i_1,i_2\in\bbZ}
  \widehat\rho(i_1)^{m+1}
  \widehat\rho(i_2)^{q_{\ell_1,m}}
  \widehat\rho(i_2-i_1)^{q_{\ell_2,m}}}
  \times
  \Bcbr{n^{-1}\sum_{j_1\in[n-1]} a(X_{\tnjo})^3
  -\int^1_0 a(X_{t})^3 dt}
  \nn\\&
  =O_M(n^{-1}),
  \label{eq:240611.2153}
\end{align}
where we used Lemma \ref{lemma:240617.2154} (iii).

By
\eqref{eq:240611.2151},
\eqref{eq:240611.2154},
\eqref{eq:240611.2155},
\eqref{eq:240611.2152} and \eqref{eq:240611.2153},
we obtain 
\begin{align}
  % n^{\half}\times\caliqtor-
  % \cali_{\qtor}^{(\ell,m;\infty)}=
  \cali_{\qtor,n}^{(\ell,m)}-
  \cali_{\qtor}^{(\ell,m;\infty)}&=
  O_M(n^{-\half})
  \label{eq:240604.1709}
\end{align}
for the case 
$q_{\ell_1,m}>0$.
Similar arguments work for the case where
$q_{\ell_1,m}=0$ and $q_{\ell_2,m}>0$,
and we have 
$\cali_{\qtor,n}^{(\ell,m)}-\cali_{\qtor}^{(\ell,m;\infty)}=
O_M(n^{-\half})$
with 
\begin{align}
  \cali_{\qtor}^{(\ell,m;\infty)}&=
  \int^1_0 a(X_{t})^3 dt \times
  \sum_{i_1,i_2\in\bbZ}
  \widehat\rho(i_1)^{m+1}
  \widehat\rho(i_2)^{q_{\ell_2,m}}.
  \label{eq:240604.1735}
\end{align}
In fact, $\cali_{\qtor}^{(\ell,m;\infty)}$ defined by \eqref{eq:240604.1735}
coincides with that by \eqref{eq:240604.1707}
when $q_{\ell_1,m}=0$.

Setting
\begin{align}
  \qtorconstRho
  % C_{\qtor}^{(\ell,m)}
  &=
  \sum_{i_1,i_2\in\bbZ}
  \widehat\rho(i_1)^{m+1}
  \widehat\rho(i_2)^{q_{\ell_1,m}}
  \widehat\rho(i_2-i_1)^{q_{\ell_2,m}}
  \nn\\
  % C_{\qtor}
  \qtorconstTotal
  &=
  \sum_{(\ell_1,\ell_2,\ell_3,m)\in\Lambda^{(k;3)}_{\sharp}}
  % C^{(\ell,m;\sharp)}
  \qtorconst
  \times\qtorconstRho,
  % \times C_{\qtor}^{(\ell,m)},
  \label{eq:240611.2222}
\end{align}
we obtain 
\begin{align}
  n^{\half}\times (D_{u_n})^2M_n&=
  \sum_{(\ell_1,\ell_2,\ell_3,m)\in\Lambda^{(k;3)}_{\sharp}}
  C^{(\ell,m;\sharp)}\times
  \cali_{\qtor,n}^{(\ell,m)}
  +O_M(n^{\half-H})
  \nn\\&=
  \sum_{(\ell_1,\ell_2,\ell_3,m)\in\Lambda^{(k;3)}_{\sharp}}
  C^{(\ell,m;\sharp)}\times
  % \cali_{\qtor,n}^{(\ell,m)}
  \rbr{\cali_{\qtor}^{(\ell,m;\infty)}+O_M(n^{-\half})}
  +O_M(n^{\half-H})
  %
  % \nn\\&=
  % \sum_{(\ell_1,\ell_2,\ell_3,m)\in\Lambda^{(k;3)}_{\sharp}}
  % C^{(\ell,m;\sharp)}\times
  % \cali_{\qtor}^{(\ell,m;\infty)}
  % +O_M(n^{\half-H})
  %
  % \nn\\&=
  % \sum_{(\ell_1,\ell_2,\ell_3,m)\in\Lambda^{(k;3)}_{\sharp}}
  % C^{(\ell,m;\sharp)}\times
  % % \cali_{\qtor}^{(\ell,m;\infty)}
  % C_{\qtor}^{(\ell,m)}\times \int^1_0 a(X_{t})^3 dt 
  % +O_M(n^{\half-H})
  % %
  \nn\\&=
  % \sum_{(\ell_1,\ell_2,\ell_3,m)\in\Lambda^{(k;3)}_{\sharp}}
  % C^{(\ell,m;\sharp)}\times C_{\qtor}^{(\ell,m)}
  % C_{\qtor}
  \qtorconstTotal
  \times \int^1_0 a(X_{t})^3 dt 
  +O_M(n^{\half-H}).
  \nn
\end{align}
\end{proof}

\printbibliography

\end{document}